\documentclass[a4paper,11pt]{amsart}
\usepackage{amsmath,amssymb,amsthm}
\usepackage{tikz-cd}
\usepackage[margin=1in]{geometry}
\usepackage[hidelinks]{hyperref}
\newtheorem{theorem}{Theorem}[section]
\newtheorem{proposition}[theorem]{Proposition}
\newtheorem{lemma}[theorem]{Lemma}
\newtheorem{corollary}[theorem]{Corollary}
\newtheorem{introtheorem}{Theorem}

\theoremstyle{definition}
\newtheorem{definition}[theorem]{Definition}
\newtheorem*{definition*}{Definition}
\newtheorem{example}[theorem]{Example}
\newtheorem{conjecture}[theorem]{Conjecture}
\newtheorem{question}[theorem]{Question}
\newtheorem{hypothesis}[theorem]{Hypothesis}
\theoremstyle{remark}
\newtheorem{remark}[theorem]{Remark}
\numberwithin{equation}{section}
\newcommand{\Q}{\mathbb{Q}}
\newcommand{\Z}{\mathbb{Z}}

\DeclareMathOperator{\Br}{Br}
\DeclareMathOperator{\inv}{inv}
\DeclareMathOperator{\res}{res}
\DeclareMathOperator{\End}{End}
\DeclareMathOperator{\Gal}{Gal}
\DeclareMathOperator{\GL}{GL}
\DeclareMathOperator{\tr}{tr}
\DeclareMathOperator{\Tr}{Tr}
\DeclareMathOperator{\ind}{ind}
\DeclareMathOperator{\diag}{diag}
\DeclareMathOperator{\Frob}{Frob}
\DeclareMathOperator{\Ver}{Ver}
\newcommand{\op}{\mathrm{op}}
\newcommand{\F}{\mathbb{F}}
\newcommand{\Qp}{\mathbb{Q}_p}

\newcommand{\Qb}{\breve{\mathbb{Q}}_p}

\DeclareMathOperator{\Hom}{Hom}
\DeclareMathOperator{\Aut}{Aut}
\newcommand{\R}{\mathbb{R}}
\newcommand{\C}{\mathbb{C}}
\newcommand{\Qbar}{\overline{\mathbb{Q}}}
\newcommand{\OO}{\mathcal{O}}
\DeclareMathOperator{\Ind}{Ind}
\DeclareMathOperator{\Art}{Art}
\DeclareMathOperator{\Nm}{N}
\DeclareMathOperator{\ord}{ord}
\newcommand{\JL}{\mathrm{JL}}
\DeclareMathOperator{\Nrd}{Nrd}
\newcommand{\Bpinf}{B_{p,\infty}}
\DeclareMathOperator{\Trd}{Trd}
\DeclareMathOperator{\SL}{SL}
\hypersetup{pdftitle={Brauer classes of representations of division algebras and their Langlands parameters},pdfauthor={Ryan Xu}}

\begin{document}

\title[Brauer classes of transfers and Langlands parameters]{Brauer classes of representations of division algebras and their Langlands parameters}
\author{Ryan Xu}
\address{Universit\'e Paris-Saclay, Facult\'e des Sciences d'Orsay, 91400 Orsay, France}
\email{ke.xu@universite-paris-saclay.fr}

\begin{abstract}
For every depth-zero supercuspidal representation $\pi$ of $\GL_n(K)$ with algebraic-valued central character, we prove that the Brauer obstructions to realizing its Jacquet--Langlands transfer and its normalized Langlands parameter $\varphi_\pi\otimes|\cdot|^{-(n-1)/2}$ over their common field of rationality $E=\Q(\pi)$ sum to $(n-1)[\Bpinf\otimes_\Q E]$. Here $K/\Qp$ is finite, the transfer is to the units of the central division algebra over $K$ of invariant $1/n$, and $\Bpinf$ is the quaternion algebra over $\Q$ ramified exactly at $p$ and $\infty$. We conjecture this identity without the depth restriction. For arbitrarily large $n=p-1$, we construct simple supercuspidal representations satisfying the same sum formula. Their transfers have finite-order central character and Schur index $n$, whereas their normalized parameters have index $n/2$. Since these transfers are not self-dual, the unequal indices contradict the isomorphism of Schur algebras conjectured by Prasad and Ramakrishnan.
\end{abstract}

\maketitle

\tableofcontents

\section{Introduction}\label{sec:intro}

We begin by asking how the Jacquet--Langlands and local Langlands correspondences change the obstruction to a model over the field of rationality. Let $K$ be a finite extension of $\Q_p$ with residue field $\F_q$, and let $\pi$ be a supercuspidal representation of $\GL_n(K)$ whose central character takes algebraic values. Then $\pi$, its Jacquet--Langlands transfer $\rho$ to the units of the central division algebra over $K$ of invariant $1/n$ and its Langlands parameter $\tau$ normalized by the twist $|\cdot|^{-(n-1)/2}$ have a common field of rationality $E$, a number field. Each of $\rho$ and $\tau$ defines a class in the Brauer group of $E$, the class $o_E(\rho)$, respectively $o_E(\tau)$, of the $E$-algebra spanned by its image (Definition~\ref{def:brauer-class}), and this class vanishes if and only if the representation has a model over $E$. We show that the two classes can differ, with Schur indices $n$ and $n/2$ for unbounded $n$. In every case that we can compute, however, their sum depends only on the rank $n$, the residue characteristic $p$ and the field $E$, and it is $n-1$ times the restriction of the class of the rational quaternion algebra $\Bpinf$ ramified exactly at $p$ and $\infty$,
\begin{equation}\label{eq:intro-sum}
 o_E(\rho)+o_E(\tau)=(n-1)\res_{E/\Q}[\Bpinf]\qquad\text{in }\Br(E).
\end{equation}
We conjecture that this holds for every such $\pi$, and we prove it in
depth zero.

In depth zero, both classes are represented by cyclic algebras for one cyclic extension of $E$, with mutually inverse generators. Inverting the generator negates the class of a cyclic algebra $(L/E,\delta,b)$, and the class is additive in $b$, so the value of the central character at a uniformizer, a factor of $b$ for both classes, cancels from the sum. What remains is the class of a cyclic algebra with $b$ a rational number, and the quaternion class comes from the sign of this number and from its $p$-adic valuation, which the normalization of $\tau$ contributes (\S\ref{subsec:intro-results}).

Write $q=p^f$, let $n\ge1$, and let $D$ be the central division algebra over $K$ of invariant $1/n$. For an irreducible supercuspidal representation $\pi$ of $\GL_n(K)$, write $\rho=\JL(\pi)$ for its Jacquet--Langlands transfer to $D^\times$ and $\varphi_\pi$ for its Langlands parameter.

\begin{definition}\label{def:normalized-parameter}
The \emph{normalized Langlands parameter} of $\pi$ is
\[
 \tau=\varphi_\pi\otimes|\cdot|^{-(n-1)/2}.
\]
We call the twist by $|\cdot|^{-(n-1)/2}$ the \emph{half twist}. 
\end{definition}
Here the reciprocity map sends geometric Frobenius elements to uniformizers, $|s|=q^{-1}$ for a geometric Frobenius element $s$, and half powers use the positive square root of $q$ (\S\ref{subsec:local}).

When the central character of $\pi$ takes algebraic values, the half twist makes the map $\pi\mapsto\tau$ commute with the automorphisms of $\C$ (\S\ref{subsec:intro-classes}), so $\pi$ and $\tau$ then have the same field of rationality. In the normalization of Prasad and Ramakrishnan \cite[\S6]{PR12}, which we recall in \S\ref{subsec:order-eight}, the normalized parameter $\tau$ is the parameter attached to $\rho$.

\subsection{Fields of rationality and Brauer classes}\label{subsec:intro-classes}

\begin{definition}[Suzuki {\cite[Definition~3.0.1]{Su23}}]\label{def:rationality} Let $r$ be a representation of a group on a complex vector space, and write $r^\sigma=r\otimes_{\C,\sigma}\C$ for $\sigma\in\Aut(\C)$. The \emph{field of rationality} $\Q(r)$ is the subfield of $\C$ fixed by all $\sigma$ such that $r^\sigma\cong r$.
\end{definition}

Every field of definition of $r$ contains $\Q(r)$. For irreducible finite-dimensional $r$, the field $\Q(r)$ is generated by the values of the character of $r$.

The Jacquet--Langlands correspondence commutes with the automorphisms of $\C$, as Dat remarks after \cite[Th\'eor\`eme~2.1.2]{Dat07}, and so does the map $\pi\mapsto\tau$ \cite[Fait~2.2.3]{Dat07}. So when the central character of $\pi$ takes algebraic values, the three fields $\Q(\pi)$, $\Q(\rho)$ and $\Q(\tau)$ coincide, and the common field $E$ is a number field (Lemma~\ref{lem:common-rationality-field}). When the central character has finite order, Suzuki shows that the Jacquet--Langlands correspondence and the map $\pi\mapsto\varphi_\pi(-\frac{n-1}2)$ to Suzuki's twisted parameter commute with these automorphisms \cite[Lemma~4.0.3]{Su23}.

But the field of rationality need not be a field of definition. The two-dimensional irreducible representation of the quaternion group of order eight has a rational character, and the $\Q$-span of its image is the algebra of rational Hamilton quaternions, a division algebra, so the representation has no model over $\Q$ (Definition~\ref{def:model-span}, Lemma~\ref{lem:models} and \cite[Remarks~3.0.4 and~3.0.10]{Su23}). Following Suzuki, we measure this obstruction by the class of the algebra spanned by the image.

\begin{definition}[cf.~{\cite[Definition~3.0.6]{Su23}}]\label{def:brauer-class}
Let $F$ be a field of characteristic zero with algebraic closure $\overline F$, and let $r\colon\Gamma\to\GL(V)$ be an irreducible representation of a group $\Gamma$ on a nonzero finite-dimensional $\overline F$-vector space $V$, with character values in $F$. The \emph{image algebra} $A_r$ is the $F$-linear span of the operators $r(\gamma)$, $\gamma\in\Gamma$. It is an $F$-algebra because $r(\gamma)r(\gamma')=r(\gamma\gamma')$, and it is central simple of degree $\dim r$ by Lemma~\ref{lem:trace-span}. The \emph{Brauer class} of $r$ is
\[
 o_F(r)=[A_r]\in\Br(F).
\]
By Wedderburn's theorem \cite[IV, Theorem~1.17 and Corollary~1.22]{Mil20}, $A_r\cong M_t(\mathcal D_r)$ for an integer $t$ and a central division algebra $\mathcal D_r$ over $F$, unique up to isomorphism, called the \emph{Schur algebra} of $r$. The \emph{Schur index} of $r$ is $\ind(A_r)$, the degree of $\mathcal D_r$.
\end{definition}

We apply this definition with $F=E$, and with $\C$ in place of $\overline F$, to the complex representations $\rho$ and $\tau$. For such a representation $r$ and a field $L$ with $E\subseteq L\subseteq\C$, the representation $r$ has a model over $L$ if and only if $A_r\otimes_EL$ is split (Lemma~\ref{lem:models} and \cite[Lemma~3.0.7]{Su23}). Thus $o_E(r)$ measures the obstruction to a model over $E$, and the Schur algebra, a name we take from Prasad and Ramakrishnan \cite[\S6]{PR12}, is the division algebra in this class.

\subsection{History}\label{subsec:intro-history}

For representations of general linear groups, the field of rationality is a field of definition in the following two settings. Waldspurger proved that an irreducible admissible representation of $\GL_2$ over a nonarchimedean local field has a model over every subfield of $\C$ containing its field of rationality, using the uniqueness of newvectors \cite[Lemma~I.2.1]{Wal85}. Suzuki proved that a square-integrable representation of $\GL_n(K)$ whose central character has finite order has a model over its field of rationality, using the uniqueness of Whittaker models \cite[Lemma~4.0.4]{Su23}. Thus, in the finite-order case, the remaining obstructions are those of $\rho$ and $\tau$.

Prasad and Ramakrishnan compared $\rho$ and $\varphi_\pi$ through signs. A self-dual irreducible representation is orthogonal or symplectic, and a self-dual irreducible admissible unitary representation of a $p$-adic group is orthogonal if and only if it has a model over $\R$ \cite[Lemma~2.1]{PR12}. Prasad and Ramakrishnan proved that a self-dual irreducible representation of $D^\times$ is orthogonal when $n$ is odd, and that for $n$ even it is orthogonal if and only if its Langlands parameter is symplectic \cite[Theorem~A and Corollary~B]{PR12}.

Prasad and Ramakrishnan then asked how the Schur algebras of $\rho$ and of the normalized parameter are related \cite[Question~6.1]{PR12}. Their Conjecture~6.1 describes the common centre of the two Schur algebras as a cyclotomic field, and its part~(1) asserts that the two algebras are equal when $\rho$ is not self-dual \cite[Conjecture~6.1(1)]{PR12}. They also noted that the realization of both correspondences in $\ell$-adic cohomology, by Harris and Taylor, determines the relation between the two Schur algebras at the places of $E$ not above $p$ or $\infty$ \cite[\S6]{PR12}.

For finite groups of Lie type, Oh used the local Langlands correspondence to study Brauer obstructions at the places not above $p$. Oh states that a depth-zero supercuspidal representation of $\GL_n(K)$ whose central character takes a rational value at a fixed uniformizer has no $\ell$-adic Brauer obstruction for $\ell\ne p$ \cite[Theorem~7]{Oh16}. Through the passage from cuspidal representations of $\GL_n(\F_q)$ to depth-zero supercuspidal representations of $\GL_n(K)$, Oh deduces from it the classical fact that cuspidal representations of $\GL_n(\F_q)$ have no $\ell$-adic Brauer obstruction \cite[Theorem~8]{Oh16}. In depth zero we prove that the sum $o_E(\rho)+o_E(\tau)$ has invariant zero at every finite place away from $p$ for every algebraic-valued central character (Theorem~\ref{thm:arbitrary-central-value}).

At the places above $p$, Sherman uses the $p$-adic Langlands program to recover the vanishing for cuspidal representations of $\GL_2(\F_p)$ and $\SL_2(\F_p)$ \cite[Theorem~1.1.1]{She18}. Sherman also compares Brauer obstructions under Fontaine's construction \cite[Lemma~4.2.5]{She18}.

Suzuki phrased the problem through the classes in $\Br(E)$ and observed that the natural quantity is their sum \cite[Remark~4.0.8]{Su23}. For central characters of finite order, Suzuki states that the local invariants of the Schur algebras of $\JL(\pi)$ and of the twisted parameter $\varphi_\pi(-\frac{n-1}2)$ add up to zero at every finite place of $\Q(\pi)$ not above $p$ \cite[Theorem~4.0.5]{Su23}. The proof uses the quotient of the Lubin--Tate tower by $\varpi^{\Z}$, for a uniformizer $\varpi$ of $K$. Since $\varpi$ acts trivially on the cohomology of this quotient, the argument requires $\omega_\pi(\varpi)=1$. Under the hypothesis of Remark~\ref{rem:suzuki-step}, the quotient by $\varpi^{d\Z}$ gives the vanishing for every finite-order central character, where $d$ is the order of $\omega_\pi(\varpi)$. Even for $d=1$, the passage to our notation assumes that Mieda's description of the cohomology \cite[Theorem~3.1]{Mie16} holds in our normalization and that Suzuki's twisted parameter is $\tau$. In even rank, Suzuki combines the vanishing away from $p$ and $\infty$ with the real invariants of Prasad and Ramakrishnan to determine the total above $p$ \cite[(4.0.6)]{Su23}. Our conjecture predicts the invariant of the sum at each place above $p$.

\subsection{The conjecture and a first example}\label{subsec:intro-conjecture}

In every case treated in this paper, the two classes satisfy the identity \eqref{eq:intro-sum} of the opening paragraph, which we now state as a conjecture.

\begin{conjecture}\label{conj:intro-sum}
Let $\pi$ be an irreducible supercuspidal representation of $\GL_n(K)$ whose central character takes algebraic values. Let $\rho=\JL(\pi)$, let $\tau$ be as in Definition~\ref{def:normalized-parameter}, and let $E$ be their common field of rationality as in \S\ref{subsec:intro-classes}. Then \eqref{eq:intro-sum} holds.
\end{conjecture}

By the theorem of Albert, Brauer, Hasse and Noether, a class in $\Br(E)$ is determined by its local invariants \cite[VII, Theorem~7.1]{Mil20}, so the conjecture can be stated place by place. For a place $v$ of $E$ we write $E_v$ for the completion of $E$ at $v$ and $\inv_v$ for the local invariant at $v$, normalized as in \S\ref{subsec:local}, and we use the same notation for the places of other number fields. In terms of invariants, the conjecture asserts that
\[
 \inv_vo_E(\rho)+\inv_vo_E(\tau)=
 \begin{cases}
 \frac{n-1}2\,[E_v:\Q_p] & \text{at each place $v$ above $p$},\\
 \frac{n-1}2 & \text{at each real place $v$},\\
 0 & \text{at every other place $v$}.
 \end{cases}
\]
For $n=1$ the conjecture holds, since $\rho$ and $\tau$ are then characters with values in $E$. 

An example in rank two shows why the field of rationality and the half twist matter.

\begin{example}\label{ex:intro-five}
The half twist can create a Brauer obstruction over the field of rationality. Take $n=2$ and $K=\Q_5$, and let $D$ be the quaternion division algebra. Choose an unramified quadratic extension $K_2\subset D$ and a character $\theta$ of $\F_{25}^\times$ of order $8$. Extend $\theta$ to $K_2^\times(1+\mathfrak p_D)$, trivially on $5$ and on $1+\mathfrak p_D$, and induce to $D^\times$ to obtain $\rho$. Let $\pi$ be its Jacquet--Langlands correspondent and put $\tau=\varphi_\pi\otimes|\cdot|^{-1/2}$.

The common field of rationality is $E=\Q(i)$, and the representations $\rho$ and $\varphi_\pi$ have models over $E$. The image algebra of $\tau$ is
\[
 A_\tau\cong L\oplus Lu,\qquad L=\Q(\zeta_8),\qquad u^2=-5,\qquad ux=\gamma_5(x)u\quad(x\in L),
\]
where $\gamma_5(\zeta_8)=\zeta_8^5$. In the notation of Definition~\ref{def:cyclic}, this is the cyclic algebra $(L/E,\gamma_5,-5)$. This quaternion algebra has invariant $1/2$ at each of the two places of $E$ above $5$ and invariant $0$ at every other place, so it is a division algebra and $\tau$ has no model over $E$. Consequently $o_E(\rho)=0$ and $o_E(\tau)=\res_{E/\Q}[B_{5,\infty}]$. Theorem~\ref{thm:order-eight}, applied with $p=5$, proves these assertions.
\end{example}

Both classes vanish in $\Br(L)$, since $\rho$ and $\tau$ have models over $L$, so the obstruction concerns the descent of a model from $L$ to the field of rationality $E$. To see where it comes from, choose a geometric Frobenius element $s$ corresponding to the uniformizer $5$. The image of inertia spans $L$, and $\tau(s)$ is the element $u$ of the example. The half twist multiplies the square of the Frobenius operator by $5$,
\[
 \varphi_\pi(s)^2=-1=i\,\gamma_5(i)=\Nm_{L/E}(i),\qquad \tau(s)^2=5\,\varphi_\pi(s)^2=-5.
\]
Before the twist the square is a norm from $L$ to $E$, and accordingly $A_{\varphi_\pi}$ is split. At either place of $E$ above $5$, the completion of $L$ is the unramified quadratic extension of $E_v=\Q_5$, whose norms have even valuation, so $-5$ is not a norm there. Thus the half twist changes the split algebra $A_{\varphi_\pi}$ into the division algebra $A_\tau$. 

\subsection{Results}\label{subsec:intro-results}

Conjecture~\ref{conj:intro-sum} holds for every depth-zero supercuspidal representation with algebraic-valued central character (Corollary~\ref{cor:depth-zero-relation}). Our first theorem treats depth zero for $n\ge2$. Proposition~\ref{prop:image-algebras} computes both classes for every algebraic value, as classes of cyclic algebras.

The representations come from characters of an unramified extension of degree $n$. Let $K_n\subset D$ be an unramified extension of $K$ of degree $n$, choose a geometric Frobenius element $s$ corresponding to a uniformizer $\varpi$, and let $U^1_D=1+\mathfrak p_D$. Let $\theta$ be a character of $\F_{q^n}^\times$ of order $m$ whose orbit under $\theta\mapsto\theta^q$ has $n$ elements. We view $\theta$ as a character of $\OO_{K_n}^\times$ through the residue map, and as a character of the inertia group $I_K$ through local class field theory (\S\ref{subsec:depth-zero-reps}). For $z\in\Qbar^\times$, let $\chi_z$ be the character of $K_n^\times U^1_D$ that is trivial on $U^1_D$ and $\mu_z$ the character of $W_{K_n}$ such that
\[
\begin{aligned}
 \chi_z&=\theta\ \text{on }\OO_{K_n}^\times, & \chi_z(\varpi)&=z,\\
 \mu_z&=\theta\ \text{on }I_K, & \mu_z(s^n)&=(-1)^{n-1}z,
\end{aligned}
\]
and put
\[
 \rho_z=\Ind_{K_n^\times U^1_D}^{D^\times}\chi_z,\qquad \tau_z=\Ind_{W_{K_n}}^{W_K}(\mu_z)\otimes|\cdot|^{-(n-1)/2}.
\]

By a theorem of Bushnell and Henniart on the level-zero correspondences \cite[6.3, Corollary~1]{BH11L}, which proves a formula proposed by Prasad and Ramakrishnan \cite[Theorem~4.4]{PR12}, the representation $\rho_z$ is the Jacquet--Langlands transfer of a depth-zero supercuspidal representation $\pi_z$ of $\GL_n(K)$ whose central character sends $\varpi$ to $z$, and $\tau_z$ is the normalized Langlands parameter of $\pi_z$ (Proposition~\ref{prop:depth-zero-parameter}). Every depth-zero supercuspidal representation of $\GL_n(K)$ with algebraic-valued central character is one of the representations $\pi_z$, for a suitable choice of $\theta$ and $z$ (Proposition~\ref{prop:depth-zero-exhaustion}).

Finally, put
\[
 L=\Q(\zeta_m),\qquad \gamma\colon\zeta_m\mapsto\zeta_m^q,\qquad E=L^{\langle\gamma\rangle},\qquad E_z=E(z),
\]
so that $E_1=E$. Example~\ref{ex:intro-five} is the case $n=2$, $K=\Q_5$, $m=8$ and $z=1$ (see \S\ref{subsec:order-eight}), and there $L=\Q(\zeta_8)$ and $E=\Q(i)$.

\begin{introtheorem}\label{thm:intro-depth-zero}
Let $n\ge2$ and $z\in\Qbar^\times$.
\begin{enumerate}
\item[(i)] The representations $\rho_z$ and $\tau_z$ have field of rationality $E_z$, and
\[
 o_{E_z}(\rho_z)+o_{E_z}(\tau_z)=(n-1)\res_{E_z/\Q}[\Bpinf] \qquad\text{in }\Br(E_z).
\]
\item[(ii)] If $v$ is a place of $E_z$ above $p$ at which $z$ is a unit, then $\inv_vo_{E_z}(\rho_z)=0$ and $\inv_vo_{E_z}(\tau_z)=\frac{n-1}2[E_{z,v}:\Q_p]$.
\item[(iii)] The representation $\rho_1$ has a model over $E$. The representation $\tau_1$ has a model over $E$ if and only if $n$ is odd, or $n$ and $f$ are even and $E$ is totally imaginary. Otherwise the Schur index of $\tau_1$ is $2$.
\end{enumerate}
\end{introtheorem}

By Lemma~\ref{lem:common-rationality-field} and part~(i), $\Q(\pi_z)=\Q(\rho_z)=E_z$ for every $z\in\Qbar^\times$, so part~(i) is Conjecture~\ref{conj:intro-sum} for $\pi_z$.

We prove Theorem~\ref{thm:intro-depth-zero} in Section~\ref{sec:depth-zero} by computing the two image algebras. Both representations are sums of $n$ lines, on which $\OO_{K_n}^\times$, respectively the inertia group, acts by the $n$ distinct characters $\theta^{q^i}$. The element $\Pi$ of \eqref{eq:Pi}, whose $n$-th power is $\varpi$, permutes the lines of $\rho_z$ cyclically, and $s$ permutes the lines of $\tau_z$ cyclically. Their $n$-th powers act by $z$ on $\rho_z$ and by $az$ on $\tau_z$, where $a$ is the rational number
\[
 a=(-1)^{n-1}q^{n(n-1)/2}.
\]
The scalar $a$ combines the sign in $\mu_z(s^n)$ with the factor $q^{n(n-1)/2}$ by which the half twist multiplies the action of $s^n$. In the following picture, each node is a line of $\rho_z$ or of $\tau_z$, labelled by the character through which $\OO_{K_n}^\times$, respectively $I_K$, acts on it.
\[
\begin{tikzcd}[column sep=2.2em, row sep=0.8em]
 \rho_z & \theta \arrow[r, "\Pi"] & \theta^{q^{-1}} \arrow[r, "\Pi"]
 & \cdots \arrow[r, "\Pi"]
 & \theta^{q^{-(n-1)}} \arrow[lll, start anchor=north, end anchor=north,
 bend right=30, "\Pi"']
 & \rho_z(\Pi)^n=z\\
 \tau_z & \theta \arrow[r, "s"] & \theta^{q} \arrow[r, "s"]
 & \cdots \arrow[r, "s"]
 & \theta^{q^{n-1}} \arrow[lll, start anchor=south, end anchor=south,
 bend left=30, "s"]
 & \tau_z(s)^n=az
\end{tikzcd}
\]
The two cycles run through the characters $\theta^{q^i}$ in opposite directions, because conjugation by $\Pi$ acts on the values of $\theta$ through $\gamma$ and conjugation by $s$ acts through $\gamma^{-1}$. This is why the cyclic algebras behind the two image algebras have mutually inverse generators. Let $L_z=LE_z$, a cyclic extension of $E_z$, and let $\gamma_z$ be the generator of $\Gal(L_z/E_z)$ that restricts to $\gamma^{[L\cap E_z:E]}$ on $L$ (Definition~\ref{def:coefficient-fields}), so that $\gamma_1=\gamma$. The two image algebras are matrix algebras over the cyclic algebras $(L_z/E_z,\gamma_z,z)$ and $(L_z/E_z,\gamma_z^{-1},az)$ (Lemma~\ref{lem:shape} and Proposition~\ref{prop:image-algebras}), and inverting the generator negates the class, so
\[
 o_{E_z}(\rho_z)+o_{E_z}(\tau_z) = [(L_z/E_z,\gamma_z,z)]-[(L_z/E_z,\gamma_z,az)] = -[(L_z/E_z,\gamma_z,a)].
\]
For $z=1$ we have $o_E(\tau_1)=[(L/E,\gamma^{-1},a)]$, and Theorem~\ref{thm:central-value-one} shows that this class is $(n-1)\res_{E/\Q}[\Bpinf]$, which has order at most two, so that it also equals $[(L/E,\gamma,a)]$. The class $[(L_z/E_z,\gamma_z,a)]$ is the restriction of $[(L/E,\gamma,a)]$ to $E_z$, so the sum is $(n-1)\res_{E_z/\Q}[\Bpinf]$, and it depends on $z$ only through the field $E_z=E(z)$.

The invariants of $o_E(\tau_1)$ come from the valuation of $a$ at $p$ and from the sign of $a$ at the real places. The hard step is at the primes dividing $m$, where $L/E$ can ramify and a unit need not be a norm, and there a congruence for $a$ shows that the invariants vanish (Theorem~\ref{thm:central-value-one} and Lemma~\ref{lem:norm-above-m}).

For central characters of finite order, the classes $o_{E_z}(\rho_z)$ and $o_{E_z}(\tau_z)$ can have order greater than two. For $K=\Q_2$, $n=3$, $\theta$ of order $7$ and $z=\zeta_3$, the two classes are opposite and have order three, and their nonzero invariants occur at the two places above $7$ (Example~\ref{ex:order-three}). A class of order three differs from its opposite, so the sign of the Brauer class, which Section~\ref{sec:brauer-classes} fixes, matters for the sum.

Now one might ask how far apart the two classes can be. Wherever \eqref{eq:intro-sum} holds, the Schur indices of $\rho$ and $\tau$ at each place of $E$, and globally, have the same odd part and differ by at most a factor of two (Propositions~\ref{prop:local-indices} and~\ref{prop:global-indices}). The reason is that the local invariants of the right side of \eqref{eq:intro-sum} lie in $\{0,\frac12\}$, and two elements of $\Q/\Z$ whose sum is $0$ or $\frac12$ have orders with the same odd part and with ratio $\frac12$, $1$ or $2$ (Lemma~\ref{lem:half-shift}). Both extremes occur at places above $p$ (Remark~\ref{rem:cases-occur}). The following family attains the global bound in arbitrarily large ranks.

\begin{introtheorem}\label{thm:intro-family}
Let $p\ge7$ be a prime with $p\equiv3\pmod4$, put $n=p-1$, and let $K$ be the unramified extension of $\Q_p$ of degree $n$. There is a simple supercuspidal representation $\pi$ of $\GL_n(K)$, with central character of finite order, such that the following hold.
\begin{enumerate}
\item[(i)] The representation $\pi$ is not self-dual, and the common field of rationality of $\pi$, $\rho$ and $\tau$ is $E=\Q(\zeta_n)$.
\item[(ii)] The prime $p$ splits completely in $E$, and the places of $E$ above $p$ can be labelled $v_j$, $j\in(\Z/n\Z)^\times$, so that
\[
 \inv_{v_j}o_E(\rho)=\frac jn,\qquad \inv_{v_j}o_E(\tau)=-\frac jn-\frac12.
\]
Both classes have invariant $0$ at every other place of $E$.
\item[(iii)] The Schur indices of $\rho$ and $\tau$ are $n$ and $n/2$, and $o_E(\rho)+o_E(\tau)=\res_{E/\Q}[\Bpinf]$.
\end{enumerate}
\end{introtheorem}

The representation is the simple supercuspidal representation $\pi_{1,\chi,1}$ of Imai and Tsushima \cite[\S5.1]{IT20}, for a suitable character $\chi$ of $\F_q^\times$ of order $n$. In Theorem~\ref{thm:family} we represent $o_E(\rho)$ and $o_E(\tau)$ by cyclic algebras for the extension $E(\zeta_p)/E$, with parameters $\zeta_n^{-1}$ and $-\zeta_n$. On the Weil side the image algebra $A_\tau$ is itself this cyclic algebra. The representation $\rho$ of $D^\times$, however, has dimension $(q^n-1)/(q-1)$, which is larger than $n$, so computing its class takes one further step. Mackey's formula splits the restriction of $\rho$ to $U^1_D$ into lines, and the projection onto $n$ suitable lines lies in $A_\rho$ and cuts out a corner of degree $n$, which is a cyclic algebra with the class of $A_\rho$ (Lemma~\ref{lem:division-side}(iii)). The sign in $-\zeta_n$ comes from the rectifier in the correspondence of Imai and Tsushima \cite[Proposition~5.4]{IT20}, and it produces the quaternion class in the sum (Remark~\ref{rem:sign}).

The Schur indices of $\rho$ and $\tau$ are $1$ and $2$ in Example~\ref{ex:intro-five} and $n$ and $n/2$ in Theorem~\ref{thm:intro-family}. In both cases $\rho$ is not self-dual, its central character has finite order, and the common centre is cyclotomic (Corollary~\ref{cor:pr-counterexamples}). The unequal indices rule out the isomorphism of Schur algebras conjectured in this case by Prasad and Ramakrishnan \cite[Conjecture~6.1(1)]{PR12}. The factor of two allowed by \eqref{eq:intro-sum} is attained in both directions, and in the family both indices are unbounded.

\subsection{Method and open questions}\label{subsec:intro-method}

We prove Theorems~\ref{thm:intro-depth-zero} and~\ref{thm:intro-family} by presenting both classes by cyclic algebras and computing their local invariants with local class field theory \cite{Mil20}. The cyclic presentations come from induced forms of $\rho$ and $\tau$, which we take in depth zero from the theorem of Bushnell and Henniart \cite[6.3, Corollary~1]{BH11L} and for the simple supercuspidal family from the description of Imai and Tsushima \cite[Proposition~5.4]{IT20}, who deduce it from the essentially tame correspondences of Bushnell and Henniart \cite{BH05,BH11J}. 

Our proofs need $\rho$ and $\tau$ to be induced from characters, and the conjecture is open in general. Mieda's form of the realization of both correspondences in the cohomology of the Lubin--Tate tower \cite[Theorem~3.1]{Mie16} enters only in Remark~\ref{rem:suzuki-step} and its hypothesis, which concern the finite places not above $p$. Above $p$, the predicted value has a meaning that we explain in Section~\ref{sec:questions}. Isocrystals are vector spaces with a semilinear Frobenius. The slopes of an isocrystal measure the $p$-adic growth of its Frobenius, and at a place $v$ above $p$ the predicted value is the invariant, after restriction to $E_v$, of the endomorphism algebra of the simple isocrystal of slope $\frac{n-1}2$. In depth zero, suppose that $\rho_z\boxtimes\tau_z$ admits a $p$-adic realization with coefficients in $E_{z,v}$, for a place $v$ of $E_z$ above $p$, as in Hypothesis~\ref{hyp:isocrystal}. If $z$ is a unit at $v$, the half twist then produces the slope $\frac{n-1}2$, and the Frobenius commutant, the algebra of endomorphisms of the realization as an isocrystal that are linear over $E_{z,v}$, has the predicted invariant (Proposition~\ref{prop:conditional-isocrystal}). We do not know whether the $p$-adic cohomology of the Lubin--Tate tower provides such a realization. We then ask whether $p$-adic cohomology can compute the invariant above $p$ in general, what determines each class separately, and what happens for division algebras of other invariants.

\subsection{Organization}\label{subsec:intro-organization}

The weight of the paper lies in Sections~\ref{sec:depth-zero} and~\ref{sec:individual}, which prove Theorems~\ref{thm:intro-depth-zero} and~\ref{thm:intro-family}. Section~\ref{sec:brauer-classes} fixes the sign of the Brauer class and the conventions for Frobenius, cyclic algebras and local invariants, identifies the common field of rationality and contains Remark~\ref{rem:suzuki-step}, and the reader may consult it as needed. Section~\ref{sec:counterexamples} treats the rank-two examples. In \S\ref{subsec:dyadic} we use the normalization hypothesis of Remark~\ref{rem:suzuki-step} to compare their individual invariants above $p$ with the two assertions of Suzuki's proposition \cite[Proposition~4.0.13]{Su23}. Section~\ref{sec:questions} states the open questions. Appendix~\ref{app:section-two-proofs} contains the routine proofs of Lemma~\ref{lem:hom-base-change} and Proposition~\ref{prop:intertwiner-sign}, and Appendix~\ref{app:isocrystal} the conditional computation above $p$.

\section{Brauer classes of representations and sign conventions}
\label{sec:brauer-classes}

This section sets up the two tools with which both theorems of the introduction are proved, the Brauer class of a representation and cyclic algebras with their local invariants. In \S\S\ref{subsec:brauer-class}--\ref{subsec:factor-set} the letter $F$ denotes a general field of characteristic zero, and in \S\ref{subsec:common-rationality-field} we return to the setting of Section~\ref{sec:intro} and show that $\pi$, $\rho$ and $\tau$ have a common field of rationality $E$, a number field, so that $o_E(\rho)$ and $o_E(\tau)$ lie in one Brauer group.

The class $o_F(r)$ of Definition~\ref{def:brauer-class} carries a sign. A central simple algebra and its opposite have the same splitting fields and the same index, so whether $r$ has a model over a given field depends only on the pair of classes $\pm o_F(r)$. Since we add the classes of different representations, however, the sign of $o_F(r)$ matters. Replacing one summand $x$ of a sum by $-x$ subtracts $2x$ from the sum, which changes it as soon as $x$ has order greater than two, as in Example~\ref{ex:order-three}. We use the image algebra, as Suzuki does \cite[Definition~3.0.6 and Lemma~3.0.7]{Su23}. Remark~\ref{rem:sign-module} and \S\ref{subsec:factor-set} describe two other algebras attached to $r$, which represent $-o_F(r)$, and \S\ref{subsec:local} fixes the conventions for Frobenius and for cyclic algebras that determine local invariants.

\subsection{The Brauer class of a representation}
\label{subsec:brauer-class}

Throughout this subsection, $F$ is a field of characteristic zero, $\overline F$ is an algebraically closed field containing $F$, $\Gamma$ is a group, and $r\colon\Gamma\to\GL(V)$ is an irreducible representation of $\Gamma$ on a nonzero $\overline F$-vector space $V$ of finite dimension $d$, with character $\gamma\mapsto\tr r(\gamma)$ taking values in $F$. The group $\Gamma$ may be infinite, as the Weil group $W_K$ is. We use the image algebra $A_r$, the class $o_F(r)$ and the Schur algebra $\mathcal D_r$ of Definition~\ref{def:brauer-class}, with $\overline F$ in place of an algebraic closure of $F$, as Section~\ref{sec:intro} does with $\C$.

\begin{definition}\label{def:model-span}
Let $L$ be a field with $F\subseteq L\subseteq\overline F$. A \emph{model} of $r$ over $L$ is a representation $r_0$ of $\Gamma$ on an $L$-vector space $V_0$ together with an isomorphism of representations $V_0\otimes_L\overline F\cong V$.
\end{definition}

\begin{lemma}\label{lem:trace-span}
The algebra $A_r$ is a central simple $F$-algebra of degree $d$, and the natural map $A_r\otimes_F\overline F\to\End_{\overline F}(V)$ is an isomorphism.
\end{lemma}

\begin{proof}
By Burnside's theorem, which is the density theorem \cite[Theorem~2.5(i)]{EGH11} applied to the group algebra $\overline F[\Gamma]$, some operators $r(\gamma_1),\dots,r(\gamma_{d^2})$ form an $\overline F$-basis of $\End_{\overline F}(V)$. The Gram matrix $\bigl(\tr r(\gamma_i\gamma_j)\bigr)_{i,j}$ of the trace form in this basis is invertible with entries in $F$, and every $a\in A_r$ has its traces $\tr\bigl(a\,r(\gamma_j)\bigr)$ in $F$, so the coordinates of $a$ lie in $F$. Hence $A_r=\bigoplus_i F\,r(\gamma_i)$, the natural map $A_r\otimes_F\overline F\to\End_{\overline F}(V)$ is an isomorphism, and $A_r$ is central simple of degree $d$, since its two-sided ideals and its centre extend to those of $M_d(\overline F)$.
\end{proof}

\begin{lemma}[cf.~{\cite[Lemma~3.0.7]{Su23}}]\label{lem:models}
We have the following: 
\begin{enumerate}
\item[(i)] Let $L$ be a field with $F\subseteq L\subseteq\overline F$. Then $r$ has a model over $L$ if and only if $A_r\otimes_FL$ is split. In particular, $r$ has a model over $F$ if and only if $o_F(r)=0$.
\item[(ii)] Let $W$ be a simple left $A_r$-module. Then $\End_{A_r}(W)\cong\mathcal D_r^{\op}$, and the class of $\End_{A_r}(W)$ in $\Br(F)$ is $-o_F(r)$.
\end{enumerate}
\end{lemma}

\begin{proof}
(i) If $A_r\otimes_FL\cong M_d(L)$ (the degree is $d$ by Lemma~\ref{lem:trace-span}), let $\Gamma$ act on the simple module $L^d$ of $M_d(L)$ through $\gamma\mapsto r(\gamma)\otimes1$. After extension of scalars to $\overline F$, this module has dimension $d$ over $\overline F$ and is a module over $A_r\otimes_F\overline F\cong\End_{\overline F}(V)$, so it is isomorphic to $V$ compatibly with $\Gamma$. Hence $L^d$, with this action of $\Gamma$, is a model of $r$ over $L$. Conversely, a model $r_0$ on $V_0$ gives $A_r\subseteq\End_L(V_0)$, and the induced map $A_r\otimes_FL\to\End_L(V_0)$ is injective because its source is simple, hence bijective by dimension, and $A_r\otimes_FL$ is split. The case $L=F$ is the second assertion.

(ii) Write $A_r=M_t(\mathcal D_r)$. The simple left $A_r$-modules are isomorphic to the space $\mathcal D_r^t$ of columns \cite[IV, Example~1.11 and Corollaries~1.21 and~1.22]{Mil20}, whose $A_r$-linear endomorphisms are the right multiplications by $\mathcal D_r$, so $\End_{A_r}(W)\cong\mathcal D_r^{\op}$. Its class is $-[\mathcal D_r]=-o_F(r)$, since $a\otimes b\mapsto(x\mapsto axb)$ is an isomorphism $\mathcal D_r\otimes_F\mathcal D_r^{\op}\cong\End_F(\mathcal D_r)$.
\end{proof}

\begin{remark}\label{rem:sign-module}
A simple left $A_r$-module is the same as an irreducible representation of $\Gamma$ over $F$ on which the group algebra $F[\Gamma]$ acts through $A_r$, since $A_r$ is spanned by the operators $r(\gamma)$. By Lemma~\ref{lem:models}(ii), the commutant of such a representation represents $-o_F(r)$. So a convention that measures rationality by this commutant has the sign opposite to ours, and Suzuki's convention agrees with ours.
\end{remark}

\subsection{Hom spaces under extension of scalars}\label{subsec:hom-spaces}

Lemma~\ref{lem:models} detects models of a finite-dimensional representation through its Brauer class. In the cohomological arguments described in \S\ref{subsec:intro-history} and Remark~\ref{rem:suzuki-step}, however, the representations $\rho$ and $\tau$ appear in the tensor product $\pi^\vee\boxtimes\rho\boxtimes\tau$, whose factor $\pi^\vee$ is infinite-dimensional when $n\geq2$, and a model of $\rho\boxtimes\tau$ has to be extracted with the help of a model of $\pi^\vee$. Let $P$ be a model of $\pi^\vee$ over a field $F$, and let $Y$ be a representation of $\GL_n(K)\times D^\times\times W_K$ over $F$ that becomes $\pi^\vee\boxtimes\rho\boxtimes\tau$ over an extension of $F$. The natural candidate for a model of $\rho\boxtimes\tau$ over $F$ is then the space $\Hom_{\GL_n(K)}(P,Y)$, on which $D^\times\times W_K$ acts through $Y$, and part~(ii) of the next lemma gives conditions under which such a Hom space is a model.

\begin{lemma}\label{lem:hom-base-change}
Let $F\subset F'$ be fields and let $\Lambda$ be a unital associative $F$-algebra, with $F$ central.
\begin{enumerate}
\item[(i)] For a finitely generated left $\Lambda$-module $X$ and any left $\Lambda$-module $Y$, the natural map
\[
\Hom_\Lambda(X,Y)\otimes_{F}F'\longrightarrow\Hom_{\Lambda\otimes_{F}F'}(X\otimes_{F}F',Y\otimes_{F}F'),
\]
which sends $\xi\otimes c$ to the map $x\otimes c'\mapsto\xi(x)\otimes cc'$, is an isomorphism.
\item[(ii)] Let $\Gamma_1$ and $\Gamma_2$ be groups. Let $U$ be a representation of $\Gamma_1$ over $F$ that is finitely generated over $F[\Gamma_1]$ and such that $U_{F'}=U\otimes_{F}F'$ is absolutely irreducible with $\End_{\Gamma_1}(U_{F'})=F'$. Let $T'$ be a finite-dimensional representation of $\Gamma_2$ over $F'$, and let $Y$ be a representation of $\Gamma_1\times\Gamma_2$ over $F$ with $Y_{F'}\cong U_{F'}\boxtimes T'$. Then $\Hom_{\Gamma_1}(U,Y)$, with $\Gamma_2$ acting through $Y$, satisfies $\Hom_{\Gamma_1}(U,Y)\otimes_{F}F'\cong T'$ as representations of $\Gamma_2$ over $F'$. In particular $\dim_F\Hom_{\Gamma_1}(U,Y)=\dim_{F'}T'$.
\end{enumerate}
\end{lemma}

Appendix~\ref{app:section-two-proofs} proves Lemma~\ref{lem:hom-base-change}, and finite generation enters there because the images of finitely many generators of $X$ can be expanded in one finite family of elements of $F'$ that are linearly independent over $F$. Remark~\ref{rem:suzuki-step} applies part~(i) at the finite places of $E$ not above $p$, and Remark~\ref{rem:hom-lemma} applies part~(ii) twice at a place of $E$ above $p$.

\subsection{Frobenius, cyclic algebras and local invariants}
\label{subsec:local}

The invariant of a cyclic algebra of an unramified extension changes sign when the arithmetic Frobenius is replaced by the geometric one, so we fix the Frobenius on the Weil group and the normalization of invariants together. These conventions hold for the rest of the paper.

Let $K$ be a finite extension of $\Q_p$ whose residue field $\F_q$ has $q=p^f$ elements, and let $W_K$ be the Weil group of $K$. An element $w\in W_K$ acts on an algebraic closure $\overline{\F}_q$ of $\F_q$ as the $\nu(w)$-th power of the arithmetic Frobenius $x\mapsto x^q$, and we call $\nu\colon W_K\to\Z$ the \emph{degree}. We fix a lift $s\in W_K$ of geometric Frobenius, so that $\nu(s)=-1$, and we write $|w|=q^{\nu(w)}$, so that $|s|=q^{-1}$. In a finite quotient of $W_K$ in which wild inertia has trivial image, every element $t$ of the image of inertia satisfies the tame relation
\[
 sts^{-1}=t^{q^{-1}},
\]
where $q^{-1}$ is an inverse of $q$ modulo the order of $t$ (this order is prime to $p$).

On the coefficient side, let $L/F$ be an unramified extension of nonarchimedean local fields. The \emph{arithmetic Frobenius} $\Frob\in\Gal(L/F)$ is the element that induces $x\mapsto x^{q_F}$ on the residue fields, where $q_F$ is the number of elements of the residue field of $F$ \cite[I, \S1]{Mil20}.

\begin{definition}\label{def:cyclic}
Let $L/F$ be a cyclic extension of degree $d$, let $\delta$ be a generator of $\Gal(L/F)$, and let $b\in F^\times$. The \emph{cyclic algebra} $(L/F,\delta,b)$ is the $F$-algebra
\[
 (L/F,\delta,b)=\bigoplus_{i=0}^{d-1}Lu^i,\qquad ux=\delta(x)u\ \ (x\in L), \qquad u^d=b.
\]
\end{definition}

Cyclic algebras go back to Dickson (see \cite[VIII, \S2]{Mil20}), and the local invariant $\inv_F\colon\Br(F)\to\Q/\Z$ of a nonarchimedean local field $F$ is due to Hasse (see the notes to \cite[IV, \S4]{Mil20}). We use the normalization of \cite[IV, \S4]{Mil20}, which part~(ii) of the next proposition records.

\begin{proposition}[see {\cite[IV, Theorem~3.14, Lemma~3.15 and Example~4.2]{Mil20}}]
\label{prop:cyclic}
Let $L/F$ be a cyclic extension of degree $d$, let $\delta$ generate $\Gal(L/F)$, and let $b,b'\in F^\times$. Then $(L/F,\delta,b)$ is a central simple $F$-algebra of degree $d$, and the following hold.
\begin{enumerate}
\item[(i)] The algebras $(L/F,\delta^{-1},b)$ and $(L/F,\delta,b^{-1})$ are isomorphic, and $[(L/F,\delta,bb')]=[(L/F,\delta,b)]+[(L/F,\delta,b')]$ in $\Br(F)$. In particular, replacing $\delta$ by $\delta^{-1}$ negates the class.
\item[(ii)] If $F$ is a nonarchimedean local field, $L/F$ is unramified, $\delta=\Frob$ and $\varpi_F$ is a uniformizer of $F$, then
\[
 \inv_F(L/F,\Frob,\varpi_F)=1/d.
\]
Consequently, for every $k\in\Z$,
\[
 \inv_F(L/F,\Frob^{-1},\varpi_F^{\,k})=-k/d\quad\text{in }\Q/\Z.
\]
\end{enumerate}
\end{proposition}

For the isomorphism in (i), let $u$ be the generator of $(L/F,\delta^{-1},b)$ and $u'$ that of $(L/F,\delta,b^{-1})$. The element $u'^{-1}$ satisfies $u'^{-1}x\,u'=\delta^{-1}(x)$ for $x\in L$ and $(u'^{-1})^d=b$, so the map that is the identity on $L$ and sends $u$ to $u'^{-1}$ is a homomorphism $(L/F,\delta^{-1},b)\to(L/F,\delta,b^{-1})$, injective since its source is simple, and bijective by dimension. The last assertion of (i) follows from this isomorphism and additivity, because $(L/F,\delta,1)$ is split (it acts on $L$ through multiplication and $\delta$, which identifies it with $\End_F(L)$), and the consequence in (ii) follows from its first formula by (i).

\begin{theorem}[see {\cite[III, \S2, (29)]{Mil20}} and {\cite[IV, Remark~4.4(b)]{Mil20}}]
\label{thm:local-brauer}
Let $F'/F$ be a finite extension of nonarchimedean local fields.
\begin{enumerate}
\item[(a)] For $x\in\Br(F)$, $\inv_{F'}(\res_{F'/F}x)=[F':F]\inv_F(x)$.
\item[(b)] The order of a class in $\Br(F)$ equals its index. In particular, a central simple $F$-algebra of degree $h$ whose invariant has order $h$ in $\Q/\Z$ is a division algebra.
\end{enumerate}
\end{theorem}

We use local class field theory in two normalizations. For a nonarchimedean
local field $F$, they are the local Artin map of \cite[I, Theorem~1.1]{Mil20}
and the reciprocity isomorphism of Bushnell and Henniart \cite[1.4]{BH11L},
\[
\begin{aligned}
 \Art_F&\colon F^\times\longrightarrow\Gal(F^{\mathrm{ab}}/F), &&\text{uniformizers}\ \mapsto\ \text{arithmetic Frobenius elements},\\
 a_F&\colon W_F^{\mathrm{ab}}\xrightarrow{\ \sim\ }F^\times, &&\text{geometric Frobenius elements}\ \mapsto\ \text{uniformizers}.
\end{aligned}
\]
We use $\Art_F$ to compute local invariants, and it is consistent with the normalization of invariants in Proposition~\ref{prop:cyclic}(ii). We use $a_F$, as Bushnell and Henniart do, to compare with the Langlands correspondence.

We now compare the two maps. The inverse of $a_F$, followed by the map $W_F^{\mathrm{ab}}\to\Gal(F^{\mathrm{ab}}/F)$, sends uniformizers to geometric Frobenius elements and induces isomorphisms $F^\times/\Nm(M^\times)\cong\Gal(M/F)$ for the finite abelian extensions $M$ of $F$. Composed with inversion on $F^\times$, it therefore has the two properties of \cite[I, Theorem~1.1]{Mil20}, and the uniqueness there identifies it with Milne's map $\Art_F$. Hence, for $w\in W_F$ and every finite abelian extension $M$ of $F$,
\begin{equation}\label{eq:weil-galois}
 w|_M=\Art_{M/F}\bigl(a_F(w)\bigr)^{-1}.
\end{equation}

Three facts about cyclic algebras recur in the later sections. The first two say when a cyclic algebra splits, and the third computes local invariants for ramified cyclic extensions as well.

\begin{lemma}\label{lem:cyclic-split-invariant}
Let $L/F$ be a cyclic extension of degree $d$, let $\delta$ generate $\Gal(L/F)$, and let $b\in F^\times$.
\begin{enumerate}
\item[(i)] If $b$ is a norm from $L$, then $(L/F,\delta,b)$ is split.
\item[(ii)] If $F$ is a nonarchimedean local field, $L/F$ is unramified and $b$ is a unit, then $b$ is a norm from $L$, so $(L/F,\delta,b)$ is split.
\item[(iii)] If $F$ is a nonarchimedean local field and $\Art_F(b)$ acts on $L$ as $\delta^t$, then the invariant of $(L/F,\delta,b)$ is $t/d$.
\end{enumerate}
\end{lemma}

\begin{proof}
(i) Write $b=\Nm_{L/F}(c)$ with $c\in L^\times$, and let $u$ be the generator of $(L/F,\delta,b)$. The element $v=c^{-1}u$ satisfies $vxv^{-1}=\delta(x)$ for $x\in L$, and
\[
 v^d=c^{-1}\delta(c^{-1})\cdots\delta^{d-1}(c^{-1})\,u^d =\Nm_{L/F}(c)^{-1}b=1.
\]
Since $v^i$ is a nonzero element of $L$ times $u^i$, the powers $1,v,\dots,v^{d-1}$ form an $L$-basis, so $(L/F,\delta,b)\cong(L/F,\delta,1)$, which is split (see the discussion after Proposition~\ref{prop:cyclic}).

(ii) The unit $b$ is a norm from the unramified extension $L$ \cite[I, 1.9 and Theorem~1.1(b)]{Mil20}, and (i) applies.

(iii) Let $\chi$ be the character of $\Gal(L/F)$ with $\chi(\delta)=1/d$. The image $\partial\chi$ of $\chi$ under the connecting map $H^1(\Gal(L/F),\Q/\Z)\to H^2(\Gal(L/F),\Z)$ is represented by the cocycle whose value at $(\delta^i,\delta^j)$, for $0\le i,j<d$, is $0$ if $i+j<d$ and $1$ otherwise \cite[III, \S1, before Proposition~1.9]{Mil20}. So $b\cup\partial\chi$ is represented by the cocycle with the values $1$ and $b$ in these two cases \cite[III, before Proposition~3.6]{Mil20}. With $e_{\delta^i}=u^i$, the relations of Definition~\ref{def:cyclic} are those of the crossed product of this cocycle, so the invariant of $(L/F,\delta,b)$ is $\inv_F(b\cup\partial\chi)$ \cite[IV, Theorem~3.14, Corollary~3.16 and Remark~4.4(a)]{Mil20}. By \cite[III, Proposition~3.6]{Mil20}, this invariant is $\chi(\Art_F(b)|_L)=\chi(\delta^t)=t/d$.
\end{proof}

For $L/F$ unramified, $\delta=\Frob$ and a uniformizer $b=\varpi_F$, the element $\Art_F(\varpi_F)$ acts on $L$ as $\Frob$, so part~(iii) recovers the invariant $1/d$ of Proposition~\ref{prop:cyclic}(ii).

The proofs of Proposition~\ref{prop:depth-zero-parameter} and Lemma~\ref{lem:weil-side} evaluate a character of a Weil group at an $n$-th power, and the next lemma says which element of the multiplicative group corresponds to it under reciprocity. The proof applies the compatibility of reciprocity with the transfer \cite[III, Lemma~3.2]{Mil20} in finite Galois quotients.

\begin{lemma}\label{lem:weil-power}
Let $K'/K$ be a cyclic extension of degree $n$ inside $\overline K$, and let $g\in W_K$ restrict to a generator of $\Gal(K'/K)$. Then $g^n\in W_{K'}$ and $a_{K'}(g^n)=a_K(g)$.
\end{lemma}

\begin{proof}
The restriction of $g^n$ to $K'$ is trivial, so $g^n\in W_{K'}$. Put $x=a_K(g)$ and $y=a_{K'}(g^n)$. Let $M$ be a finite abelian extension of $K'$, let $N$ be a finite Galois extension of $K$ containing $M$, and let $\Ver\colon\Gal(N/K)^{\mathrm{ab}}\to\Gal(N/K')^{\mathrm{ab}}$ be the transfer. Let $N_0$ and $N_1$ be the maximal abelian subextensions of $N/K$ and of $N/K'$, so that the two abelianized groups are $\Gal(N_0/K)$ and $\Gal(N_1/K')$, and $M\subseteq N_1$. In them we write $g|_N$ and $g^n|_N$ for the images of $g$ and $g^n$, and $\Art_{N/K}(x)$ and $\Art_{N/K'}(x)$ for $\Art_{N_0/K}(x)$ and $\Art_{N_1/K'}(x)$. Then
\[
 \Ver(g|_N)=g^n|_N,\qquad \Ver\bigl(\Art_{N/K}(x)\bigr)=\Art_{N/K'}(x),\qquad g|_N=\Art_{N/K}(x)^{-1}.
\]
The first identity holds because the images of $1,g,\dots,g^{n-1}$ represent the cosets of $\Gal(N/K')$ and $g$ permutes them cyclically, with $g\cdot g^{n-1}=g^n$. The second is the second diagram of \cite[III, Lemma~3.2]{Mil20}, and the third is \eqref{eq:weil-galois} for $F=K$ and the abelian extension $N_0$. Hence $g^n|_N=\Art_{N/K'}(x)^{-1}$, and restricting from $N_1$ to $M$ gives $g^n|_M=\Art_{M/K'}(x)^{-1}$. By \eqref{eq:weil-galois} for $F=K'$, also $g^n|_M=\Art_{M/K'}(y)^{-1}$.

Hence $x/y$ lies in the kernel $\Nm(M^\times)$ of $\Art_{M/K'}$ \cite[I, Theorem~1.1(b)]{Mil20} for every finite abelian extension $M$ of $K'$. By the local existence theorem \cite[I, Theorem~1.4]{Mil20} these norm groups are the open subgroups of finite index of $K'^\times$, whose intersection is trivial, so $x=y$.
\end{proof}

\subsection{The factor set of intertwiners}
\label{subsec:factor-set}

Oh defines the Brauer obstruction of a representation through a crossed product built from intertwiners, and for finite groups and number fields Oh proved that the class of the image algebra is the negative of this obstruction \cite[\S2, Theorem~1]{Oh16}. This subsection makes the comparison for representations of arbitrary groups, and the later sections use only $o_F(r)$. When $r$ is given by matrices over a finite Galois extension $L$ of $F$, isomorphisms from the Galois conjugates of $r$ to $r$ define a factor set, and the class of its crossed product is $-o_F(r)$ (Proposition~\ref{prop:intertwiner-sign}). The sign depends on the direction of the isomorphisms and on the multiplication in the crossed product, and Definition~\ref{def:factor-set} fixes both.

\begin{definition}\label{def:factor-set}
Let $L/F$ be a finite Galois extension with group $G$.
\begin{enumerate}
\item[(i)] Let $c\colon G\times G\to L^\times$ satisfy $c(g,h)\,c(gh,k)=g\bigl(c(h,k)\bigr)\,c(g,hk)$ for $g,h,k\in G$. The \emph{crossed product} $(L/F,c)$ is the $F$-algebra $\bigoplus_{g\in G}Lu_g$ with $u_gx=g(x)u_g$ for $x\in L$ and $u_gu_h=c(g,h)\,u_{gh}$.
\item[(ii)] Let $\Gamma$ be a group and $r\colon\Gamma\to\GL_d(L)$ a representation whose commutant in $M_d(L)$ is $L$. For $g\in G$ let $r^g=g\circ r$, where $g$ acts on matrices entrywise. Given isomorphisms $\iota_g\colon r^g\to r$ for $g\in G$, their \emph{factor set} is
\[
 c(g,h)=\iota_g\,g(\iota_h)\,\iota_{gh}^{-1}\qquad(g,h\in G).
\]
\end{enumerate}
\end{definition}

In (ii), the factor set takes values in $L^\times$ and satisfies the identity of (i), and the class of $(L/F,c)$ in $\Br(F)$ depends only on $r$ (Appendix~\ref{app:section-two-proofs}).

\begin{proposition}\label{prop:intertwiner-sign}
Let $F$ be a field of characteristic zero with algebraic closure $\overline F$, let $L\subseteq\overline F$ be a finite Galois extension of $F$ with group $G$, and let $r\colon\Gamma\to\GL_d(L)$ be a representation of a group $\Gamma$ that is irreducible over $\overline F$ and whose character takes values in $F$. Then the commutant of $r$ in $M_d(L)$ is $L$. If $\iota_g\colon r^g\to r$, for $g\in G$, are isomorphisms with factor set $c$, then
\[
 [(L/F,c)]=-o_F(r)\qquad\text{in }\Br(F).
\]
\end{proposition}

Appendix~\ref{app:section-two-proofs} proves Proposition~\ref{prop:intertwiner-sign} by the double centralizer theorem.

\subsection{The common field of rationality}\label{subsec:common-rationality-field}

We now return to the arithmetic setting of Section~\ref{sec:intro} and show that $\pi$, $\rho$ and $\tau$ have the same field of rationality. An unramified twist reduces the question whether this field is a number field to central characters of finite order. We first record the compatibility of the Jacquet--Langlands correspondence with twists by arbitrary characters, whose values need not have absolute value one. In \S\ref{subsec:order-eight} we use it to identify $\tau$ with the parameter that Prasad and Ramakrishnan attach to $\rho$.

We use the Jacquet--Langlands and Langlands correspondences in the form stated by Dat, who passes between fields of coefficients in \cite[\S\S2.1.1--2.2.5]{Dat07}, so we may take the coefficients to be $\C$. Dat states the Jacquet--Langlands correspondence as a bijection from the set of irreducible representations of $D^\times$ onto the set $\mathrm{Disc}_\C(\GL_n(K))$ of discrete series representations of $\GL_n(K)$, characterized by the equality of characters $\chi_\pi(g)=(-1)^{n-1}\chi_\rho(x)$ when $\pi$ corresponds to $\rho$ and $g$ and $x$ are corresponding elliptic regular elements \cite[Th\'eor\`eme~2.1.2]{Dat07}. As in Section~\ref{sec:intro}, we write $\JL$ for the inverse bijection, so that $\rho=\JL(\pi)$.

We read $\mathrm{Disc}_\C(\GL_n(K))$ as the set of essentially square-integrable representations of $\GL_n(K)$ (for $n=1$ it consists of all characters of $K^\times$), which the Langlands correspondence maps onto the indecomposable Weil--Deligne representations of dimension $n$ \cite[2.2.7]{Dat07}. This set is stable under twisting by every character $\chi$ of $K^\times$, unitary or not, because the Langlands correspondence is compatible with twisting by characters \cite[\S2.2.1]{Dat07} and twisting by a character preserves indecomposability. The proof of the next lemma uses this stability.

\begin{lemma}\label{lem:jl-twist}
Let $\pi$ be an essentially square-integrable representation of $\GL_n(K)$ and $\chi$ a character of $K^\times$. Then $\JL(\pi\otimes(\chi\circ\det))\cong\JL(\pi)\otimes(\chi\circ\Nrd)$.
\end{lemma}

\begin{proof}
Put $\rho=\JL(\pi)$. Let $g\in\GL_n(K)$ and $x\in D^\times$ be elliptic regular elements with the same minimal polynomial $P$ of degree $n$. Then $P$ is the characteristic polynomial of $g$ and the reduced characteristic polynomial of $x$, so $\det g=\Nrd x=(-1)^nP(0)$. The character of $\pi\otimes(\chi\circ\det)$ at $g$ is therefore $(-1)^{n-1}\chi(\Nrd x)\chi_\rho(x)$, which is the value that \cite[Th\'eor\`eme~2.1.2]{Dat07} prescribes for the representation corresponding to $\rho\otimes(\chi\circ\Nrd)$. The representation $\pi\otimes(\chi\circ\det)$ lies in $\mathrm{Disc}_\C(\GL_n(K))$, as shown above, so the characterization in \cite[Th\'eor\`eme~2.1.2]{Dat07} identifies it with the representation corresponding to $\rho\otimes(\chi\circ\Nrd)$.
\end{proof}

\begin{lemma}\label{lem:common-rationality-field}
Let $\pi$ be a supercuspidal representation of $\GL_n(K)$ whose central character takes algebraic values. Then $\Q(\pi)=\Q(\JL(\pi))=\Q(\tau)$, and their common field $E=\Q(\pi)$ is a number field.
\end{lemma}

\begin{proof}
The Jacquet--Langlands correspondence commutes with the automorphisms of $\C$, as Dat remarks after \cite[Th\'eor\`eme~2.1.2]{Dat07}, and so does the map $\pi\mapsto\tau$ \cite[Fait~2.2.3]{Dat07}. Dat normalizes reciprocity, the character $|\cdot|$ and the half twist as in Definition~\ref{def:normalized-parameter} \cite[\S2.2.1]{Dat07}. Both correspondences are injective, so the three stabilizers coincide.

It remains to see that $E$ is a number field. Let $\omega$ be the central character of $\pi$, and choose an unramified character $\chi$ with $\chi(\varpi)^n=\omega(\varpi)^{-1}$ for a uniformizer $\varpi$. The value $\chi(\varpi)$ is algebraic, and $\pi_0=\pi\otimes(\chi\circ\det)$ has central character of finite order. The representation $\JL(\pi_0)$ is smooth and finite-dimensional, and its central character is that of $\pi_0$, so its kernel is open and contains a subgroup of finite index of the centre $K^\times$. As $D^\times/K^\times$ is compact, $\JL(\pi_0)$ factors through a finite quotient of $D^\times$, so its character takes values in a cyclotomic field, and $E_0=\Q(\pi_0)=\Q(\JL(\pi_0))$ is a number field. Moreover, $\pi_0$ has a model over $E_0$ \cite[Lemma~4.0.4]{Su23}, and twisting it gives a model of $\pi$ over the number field $E_0(\chi(\varpi))$. Since every field of definition of $\pi$ contains $\Q(\pi)$, the common field $E$ is a number field.
\end{proof}

Since $\rho$ and $\tau$ have the same field of rationality, their classes can be added in $\Br(E)$, and the sum is itself the Brauer class of a representation.

\begin{lemma}\label{lem:sum-class}
Let $\pi$ be as in Lemma~\ref{lem:common-rationality-field}, let $\rho=\JL(\pi)$, let $\tau$ be its normalized Langlands parameter, and let $E=\Q(\pi)$. Then the image algebra of the representation $\rho\boxtimes\tau$ of $D^\times\times W_K$ is isomorphic to $A_\rho\otimes_EA_\tau$, and
\[
 o_E(\rho\boxtimes\tau)=o_E(\rho)+o_E(\tau).
\]
\end{lemma}

\begin{proof}
Let $V_\rho$ and $V_\tau$ be the spaces of $\rho$ and $\tau$. The map $A_\rho\otimes_EA_\tau\to\End_\C(V_\rho\otimes_\C V_\tau)$, $\alpha\otimes\beta\mapsto\alpha\otimes\beta$, is a homomorphism onto the $E$-span of the operators $\rho(x)\otimes\tau(w)$, for $x\in D^\times$ and $w\in W_K$, and it is injective because its source is simple. By Lemma~\ref{lem:trace-span} the operators of $\rho$ and of $\tau$ span $\End_\C(V_\rho)$ and $\End_\C(V_\tau)$, so those of $\rho\boxtimes\tau$ span $\End_\C(V_\rho\otimes_\C V_\tau)$, and $\rho\boxtimes\tau$ is irreducible, with character values $\tr\rho(x)\,\tr\tau(w)$ in $E$. So the image algebra of $\rho\boxtimes\tau$ is $A_\rho\otimes_EA_\tau$, of class $o_E(\rho)+o_E(\tau)$.
\end{proof}

The following remark shows, for a central character of finite order, that the invariants of this sum vanish at the finite places of $E$ not above $p$, under the hypothesis stated there.

\begin{remark}[Suzuki's argument for central characters of finite order]
\label{rem:suzuki-step}
Let $\pi$ be as in Lemma~\ref{lem:sum-class}, with central character $\omega_\pi$ of finite order, and let $v$ be a finite place of $E$ not above $p$. Assume the hypothesis on normalizations stated below. Then
\[
 \inv_vo_E(\rho)+\inv_vo_E(\tau)=0.
\]

This vanishing is the form in our notation of the theorem that Suzuki  states \cite[Theorem~4.0.5]{Su23}, recalled in Suzuki's terms in \S\ref{subsec:intro-history}. Suzuki's proof uses the quotient of the Lubin--Tate tower by $\varpi^{\Z}$, for a uniformizer $\varpi$ of $K$. The central element $\varpi$ of $D^\times$ acts trivially on the cohomology of this quotient, so $\pi^\vee\boxtimes\JL(\pi)\boxtimes\varphi_\pi(-\frac{n-1}2)$ occurs in it only when $\varpi$ acts trivially on $\rho=\JL(\pi)$. When the central character has finite order, however, $\varpi$ acts on $\rho$ by a root of unity, which can be nontrivial (in Example~\ref{ex:order-three} it is $\zeta_3$). We follow Suzuki's argument with the quotient by the subgroup $\varpi^{d\Z}$ defined below, which is Suzuki's quotient when $\rho(\varpi)=1$, and the hypothesis enters, even in that case, at the step where we identify a component of its cohomology with $\pi^\vee\boxtimes\rho\boxtimes\tau$.

Let $\ell\neq p$ be the prime below $v$, put $k=E_v$, and choose an algebraic closure $\bar k$ of $k$ and a field isomorphism $\C\to\bar k$ extending $E\to k$, through which we view $\pi$, $\rho$ and $\tau$ over $\bar k$. It preserves the image algebras over $E$, and hence $o_E(\rho)$ and $o_E(\tau)$.

Let $\zeta=\omega_\pi(\varpi)$ and let $d$ be its order. The group $\varpi^{d\Z}$ generated by $\varpi^d$ is a discrete cocompact subgroup of $K^\times$, on which the central character of the contragredient $\pi^\vee$ is trivial. Let $H$ be the direct limit over the levels of the compactly supported \'etale cohomology in degree $n-1$, with coefficients in $\Q_\ell$, of the quotient of the Lubin--Tate tower by $\varpi^{d\Z}$ at that level, base changed to the completion of an algebraic closure of $K$ (the geometric generic fibre), with its coefficients then extended from $\Q_\ell$ to $k$ (Dat treats the integral cohomology of the quotient of the tower by the group generated by an element of $K^\times$ of positive valuation \cite[\S1.2]{Dat12}). The space $H\otimes_k\bar k$ is then the same cohomology with coefficients $\bar k$ (see \cite[\S2.5]{Str08} for the coefficients $\Q_\ell$ and their extension).

The groups $G=\GL_n(K)$, $D^\times$ and $W_K$ act on $H$ by commuting $k$-linear actions, and $H\otimes_k\bar k$ is an admissible representation of $G$ on which $\varpi^{d\Z}$ acts trivially (see \cite[\S3.1]{Mie16} for this cohomology and the three actions on it). Mieda's form of the non-abelian Lubin--Tate theory \cite[Theorem~3.1]{Mie16}, which Mieda credits to Harris and Taylor, to Boyer and to Mieda's earlier work, applies to every supercuspidal representation whose central character is trivial on $\varpi^{d\Z}$, and for $\pi^\vee$ it describes the $\pi^\vee$-isotypic component of $H\otimes_k\bar k$.

We assume that, with $W_K$ acting on $H$ as in \cite[\S3.1]{Mie16}, the $\pi^\vee$-isotypic component of $H\otimes_k\bar k$ is isomorphic to $\pi^\vee\boxtimes\rho\boxtimes\tau$ as a representation of $G\times D^\times\times W_K$, where $\rho=\JL(\pi)$ and $\tau=\varphi_\pi\otimes|\cdot|^{-(n-1)/2}$ are attached to $\pi$ in the normalizations of \S\ref{subsec:local} and \S\ref{subsec:common-rationality-field}. We assume moreover that the twisted parameter $\varphi_\pi(-\frac{n-1}2)$ of \cite[Theorem~4.0.5 and Proposition~4.0.13]{Su23} is $\tau$, even when $\rho(\varpi)=1$.

Both correspondences commute with contragredients, and Mieda's twist by $\frac{n-1}2$ is the twist by $|\cdot|^{(n-1)/2}$ in our normalization. So the isomorphism we assume is the one that \cite[Theorem~3.1]{Mie16}, applied to $\pi^\vee$, gives when the Langlands correspondence in that theorem, for this action of $W_K$ on $H$, is the one of this paper, in which uniformizers correspond to geometric Frobenius elements \cite[\S2.2.1]{Dat07}.

The other components of $H\otimes_k\bar k$ in the decomposition of \cite[\S3.1]{Mie16} are sums of copies of other supercuspidal representations or have no supercuspidal subquotient, so they admit no nonzero map from $\pi^\vee$. Since $\End_G(\pi^\vee)=\bar k$ and $\rho\boxtimes\tau$ is finite dimensional, it follows that $t\mapsto(x\mapsto x\otimes t)$ is an isomorphism of representations of $D^\times\times W_K$
\[
 \rho\boxtimes\tau\xrightarrow{\ \sim\ } \Hom_{\bar k[G]}(\pi^\vee,H\otimes_k\bar k).
\]

It remains to descend this Hom space to $k$. By a lemma of Suzuki \cite[Lemma~4.0.4]{Su23}, $\pi$ has a model over $E$, so $\pi^\vee$ has a model $P$ over $k$. Since $P\otimes_k\bar k\cong\pi^\vee$ is irreducible, $P$ is irreducible, hence cyclic, and in particular finitely generated, as a $k[G]$-module. Lemma~\ref{lem:hom-base-change}(i), applied with $\Lambda=k[G]$, $X=P$, $Y=H$ and $F'=\bar k$, identifies $\Hom_{k[G]}(P,H)\otimes_k\bar k$ with $\Hom_{\bar k[G]}(\pi^\vee,H\otimes_k\bar k)$, compatibly with $D^\times$ and $W_K$, which act through $H$. So $\Hom_{k[G]}(P,H)$ is a model of $\rho\boxtimes\tau$ over $k$. By Lemma~\ref{lem:models}(i), applied with $F=E$, $\overline F=\bar k$ and $L=k$, the image algebra of $\rho\boxtimes\tau$ splits over $k$, and its class is $o_E(\rho)+o_E(\tau)$ by Lemma~\ref{lem:sum-class}. Under the hypothesis, this proves the vanishing at $v$. The same quotient reappears above $p$ in Question~\ref{q:above-p}.
\end{remark}

\section{Representations of depth zero}\label{sec:depth-zero}

This section proves Theorem~\ref{thm:intro-depth-zero}. Proposition~\ref{prop:image-algebras} presents the image algebras of $\rho_z$ and $\tau_z$ as matrix algebras over cyclic algebras for one cyclic extension of their field of rationality $E_z$, with parameters $z$ and $az$ and with mutually inverse generators, so that $z$ cancels from the sum of the two classes. The identity for arbitrary $z$ (\S\ref{subsec:arbitrary}) therefore reduces to the computation of $o_E(\tau_1)$ for $z=1$, which local class field theory carries out place by place (\S\ref{subsec:central-value-one}). Finally, for $n\ge2$ the representations $\pi_z$ exhaust the depth-zero supercuspidal representations of $\GL_n(K)$ with algebraic-valued central character (Proposition~\ref{prop:depth-zero-exhaustion}), and the case $n=1$ is immediate (\S\ref{subsec:intro-conjecture}), so Conjecture~\ref{conj:intro-sum} holds in depth zero (Corollary~\ref{cor:depth-zero-relation}).

\subsection{The representations \texorpdfstring{$\rho_z$ and $\tau_z$}{rho\_z and tau\_z}}\label{subsec:depth-zero-reps}

Throughout this section, $K$ is a finite extension of $\Q_p$ with residue field $\F_q$, where $q=p^f$, and $n\geq2$ is an integer. We use the geometric Frobenius element $s\in W_K$, the character $|\cdot|$ and the reciprocity maps $\Art_F$ and $a_F$ of \S\ref{subsec:local}.

Let $\varpi=a_K(s)$, a uniformizer of $K$, and let $D$ be the central division algebra over $K$ of invariant $1/n$, with ring of integers $\OO_D$, maximal ideal $\mathfrak p_D$ and $U_D^1=1+\mathfrak p_D$. We fix an unramified extension $K_n\subset D$ of $K$ of degree $n$ and an element $\Pi\in D$ with
\begin{equation}\label{eq:Pi}
 \Pi^n=\varpi,\qquad \Pi x\Pi^{-1}=\Frob(x)\quad(x\in K_n),
\end{equation}
where $\Frob$ is the arithmetic Frobenius of $K_n/K$. Such a pair exists, for any uniformizer in place of $\varpi$, because $(K_n/K,\Frob,\varpi)$ has invariant $1/n$ (Proposition~\ref{prop:cyclic}(ii)) and is therefore a division algebra isomorphic to $D$ (Theorem~\ref{thm:local-brauer}(b)).

The map $a_{K_n}$ carries the inertia group $I_K=I_{K_n}$ onto $\OO_{K_n}^\times$, and reduction modulo $\mathfrak p_{K_n}$ then defines a surjection $I_K\to\F_{q^n}^\times$, $t\mapsto\bar t$. Through it we regard characters of $\F_{q^n}^\times$ as characters of $I_K$, so that such a character lives on $\OO_{K_n}^\times$ and on $I_K$ through their common quotient $\F_{q^n}^\times$.

The representations $\rho_z$ and $\tau_z$ are induced from the normal subgroups $K_n^\times U_D^1=\OO_D^\times\varpi^{\Z}$ of $D^\times$ and $W_{K_n}$ of $W_K$. Both quotients are cyclic of order $n$, generated by the image of $\Pi$, respectively of $s$, and conjugation by these generators acts by
\[
 \Pi x\Pi^{-1}\equiv x^{q}\pmod{\mathfrak p_{K_n}}\quad(x\in\OO_{K_n}),\qquad \overline{sts^{-1}}=\bar t^{\,q^{-1}}\quad(t\in I_K),
\]
by \eqref{eq:Pi} and by the tame relation of \S\ref{subsec:local}, where $q^{-1}$ is an inverse of $q$ modulo $q^n-1$. So conjugation by $\Pi$ induces the arithmetic Frobenius $x\mapsto x^q$ on the residue field $\F_{q^n}$ of $K_n$, and conjugation by the geometric Frobenius element $s$ induces its inverse on the quotient $\F_{q^n}^\times$ of $I_K$.

Let $\theta\colon\F_{q^n}^\times\to\Qbar^\times$ be a character of order $m$ whose $q$-power orbit has length $n$, that is, such that $\theta,\theta^q,\dots,\theta^{q^{n-1}}$ are distinct. We call such a character primitive. Equivalently, $q$ has order $n$ in $(\Z/m\Z)^\times$. The integer $m$ divides $q^n-1$ and does not divide $q-1$, so $p\nmid m$ and $m\geq3$.

We fix a generator $\bar t$ of $\F_{q^n}^\times$ and put $\zeta_m=\theta(\bar t)$, a primitive $m$-th root of unity. We identify $\Gal(\Q(\zeta_m)/\Q)$ with $(\Z/m\Z)^\times$, the class of $c$ acting by $\zeta_m\mapsto\zeta_m^c$, and, as in Section~\ref{sec:intro}, we put
\[
 L=\Q(\zeta_m),\qquad \gamma\colon\zeta_m\mapsto\zeta_m^q,\qquad E=L^{\langle\gamma\rangle},
\]
so that $\gamma$ is the class of $q$. Since $q$ has order $n$ modulo $m$, the extension $L/E$ is cyclic of degree $n$ and $\gamma$ generates $\Gal(L/E)$.

\begin{definition}\label{def:depth-zero-reps}
Let $z\in\Qbar^\times$.
\begin{enumerate}
\item[(i)] Let $\chi_z$ be the character of $K_n^\times U_D^1$ that is trivial on $U_D^1$, agrees with $\theta$ on $\OO_{K_n}^\times$ (through the residue field) and sends $\varpi$ to $z$. We put $\rho_z=\Ind_{K_n^\times U_D^1}^{D^\times}\chi_z$.
\item[(ii)] Let $\mu_z$ be the character of $W_{K_n}$ that agrees with $\theta$ on $I_K$ and sends $s^n$ to $(-1)^{n-1}z$. We put $\tau_z=\Ind_{W_{K_n}}^{W_K}(\mu_z)\otimes|\cdot|^{-(n-1)/2}$.
\item[(iii)] We let $\pi_z$ be the representation of $\GL_n(K)$ that corresponds to $\rho_z$ under the Jacquet--Langlands correspondence.
\end{enumerate}
\end{definition}

The character $\chi_z$ is well defined because $K_n^\times U_D^1/U_D^1\cong K_n^\times/(1+\mathfrak p_{K_n})$, and $\mu_z$ exists and is unique because $W_{K_n}$ is generated by $I_K$ and $s^n$, and conjugation by $s^n$ acts on tame inertia as the $q^{-n}$-th power, which is the identity on $\F_{q^n}^\times$. The number $z$ is the value at $\varpi$ of the central character of $\rho_z$, and (iii) uses that $\rho_z$ is irreducible (Lemma~\ref{lem:shape}).

As in Section~\ref{sec:intro}, we put
\[
 a=(-1)^{n-1}q^{n(n-1)/2}\in\Q^\times,
\]
and we recall that $B_{p,\infty}$ is the quaternion algebra over $\Q$ ramified exactly at $p$ and at the real place. The next lemma describes the shape of $\rho_z$ and $\tau_z$, which all later arguments use.

\begin{lemma}\label{lem:shape}
Let $z\in\Qbar^\times$.
\begin{enumerate}
\item[(i)] Put $\mu'_z=\mu_z\otimes\bigl(|\cdot|^{-(n-1)/2}|_{W_{K_n}}\bigr)$. Then $\tau_z\cong\Ind_{W_{K_n}}^{W_K}\mu'_z$, the character $\mu'_z$ agrees with $\theta$ on $I_K$, and $\mu'_z(s^n)=az$.
\item[(ii)] The restriction of $\rho_z$ to $\OO_{K_n}^\times$, resp.\ of $\tau_z$ to $I_K$, is the direct sum of the $n$ distinct characters $\theta^{q^i}$, $0\leq i<n$, each on a line, and $\Pi$, resp.\ $s$, permutes these lines cyclically.
\item[(iii)] The element $\varpi=\Pi^n$ acts on $\rho_z$ by $z$, and $s^n$ acts on $\tau_z$ by $az$.
\end{enumerate}
In particular $\rho_z$ and $\tau_z$ are absolutely irreducible of dimension $n$.
\end{lemma}

\begin{proof}
Part (i) is the projection formula, together with $|s^n|^{-(n-1)/2}=q^{n(n-1)/2}$ (the character $|\cdot|$ is trivial on $I_K$). If $v$ spans the line of $\rho_z$ on which $K_n^\times U_D^1$ acts by $\chi_z$, then $\rho_z$ is the direct sum of the lines $\Pi^iv$, $0\leq i<n$, which $\Pi$ permutes cyclically, and $\OO_{K_n}^\times$ acts on $\Pi^iv$ by $\theta^{q^{-i}}$, since conjugation by $\Pi^{-1}$ induces the inverse of $x\mapsto x^q$ on $\F_{q^n}$. These are the $n$ distinct characters $\theta^{q^j}$, because $\theta$ is primitive, and the central element $\varpi$ acts by $\chi_z(\varpi)=z$. The same argument with the lines $s^iv$ of $\Ind\mu'_z$ gives the characters $\theta^{q^i}$ on $I_K$, since $s^{-1}ts=t^q$ in the tame inertia quotient, and the value $\mu'_z(s^n)=az$. Since the characters on the lines are distinct, a subspace stable under $\OO_{K_n}^\times$ and $\Pi$, or under $I_K$ and $s$, is zero or everything, and the coefficients lie in $\Qbar$, so both representations are absolutely irreducible.
\end{proof}

The next proposition explains the choice of $\tau_z$. It is a formula that Prasad and Ramakrishnan proposed \cite[Theorem~4.4]{PR12} and that Bushnell and Henniart proved.

\begin{proposition}[Bushnell and Henniart {\cite[6.3, Corollary~1]{BH11L}}]
\label{prop:depth-zero-parameter}
Let $z\in\Qbar^\times$. The Langlands parameter of $\pi_z$ is $\Ind_{W_{K_n}}^{W_K}\mu_z$. Consequently $\tau_z$ is the normalized Langlands parameter of $\pi_z$.
\end{proposition}

\begin{proof}
We translate the statement of Bushnell and Henniart into our notation. They parametrize the level-zero representations by admissible tame pairs $(K'/K,\vartheta)$, where $K'/K$ is unramified and $\vartheta$ is a character of $K'^\times$, trivial on the principal units, whose $\Gal(K'/K)$-conjugates are distinct \cite[1.1]{BH11L}.

To match our data with such a pair, we write $\chi_z$ also for the restriction of $\chi_z$ to $K_n^\times$. On $\OO_{K_n}^\times$ this character is the primitive character $\theta$, and $\Gal(K_n/K)$ fixes $\varpi$, so the conjugates of $\chi_z$ are distinct and $(K_n/K,\chi_z)$ is an admissible tame pair of degree $n$. The representation of $D^\times$ that Bushnell and Henniart attach to this pair is induced from the character of $K_n^\times U_D^1$ that extends $\chi_z$ and is trivial on $U_D^1$ \cite[(1.5.1)]{BH11L}, so it is $\rho_z$.

In the case $e=1$ and $f=n$ of their notation, Corollary~1 of Bushnell and Henniart states that the Langlands parameter of the Jacquet--Langlands correspondent of this representation, which is $\pi_z$, is
\[
 \Ind_{W_{K_n}}^{W_K}\bigl((\eta^{n-1}\chi_z)\circ a_{K_n}\bigr) \qquad\text{\cite[6.3, Corollary~1, (1.4.1) and (1.4.3)]{BH11L}},
\]
where $\eta$ is the unramified quadratic character of $K_n^\times$. Here Bushnell and Henniart view characters of $K_n^\times$ as characters of $W_{K_n}$ through the reciprocity map that takes geometric Frobenius elements to prime elements \cite[1.4]{BH11L}, which is our $a_{K_n}$, and their Langlands correspondence is the one we use, that of Harris and Taylor and of Henniart.

For $e=1$, the parameter in \cite[(1.4.3)]{BH11L} is the induced representation of \cite[(1.4.1)]{BH11L}, because $\mathrm{Sp}_1$ leaves every irreducible representation of $W_K$ unchanged. Here we read the exponent $i+(1-n)/2$ in the definition of $\mathrm{Sp}_e$ \cite[1.2]{BH11L} as $i+(1-e)/2=0$, as \cite[3.1, (3.2.1) and Lemma~13]{BH11L} require.

It remains to show that the inducing character $(\eta^{n-1}\chi_z)\circ a_{K_n}$ is $\mu_z$. It agrees with $\theta$ on $I_K$, because $\eta$ is unramified and because of the way we regard $\theta$ as a character of $I_K$. Its value at $s^n$ is the one point that needs care. The element $s$ restricts to a generator of $\Gal(K_n/K)$, so $a_{K_n}(s^n)=a_K(s)=\varpi$ by Lemma~\ref{lem:weil-power}, and therefore
\[
 (\eta^{n-1}\chi_z)\bigl(a_{K_n}(s^n)\bigr) = \chi_z(\varpi)\,\eta(\varpi)^{n-1}=(-1)^{n-1}z=\mu_z(s^n).
\]
So the character is $\mu_z$, and the second assertion follows from Definition~\ref{def:normalized-parameter}.
\end{proof}

Since $\tau_z$ is irreducible (Lemma~\ref{lem:shape}), so is the Langlands parameter $\Ind_{W_{K_n}}^{W_K}\mu_z$ of $\pi_z$, and therefore $\pi_z$ is supercuspidal. It has depth zero, because $\rho_z$ is trivial on $U_D^1$ and the Jacquet--Langlands correspondence preserves depth \cite[Corollary~2.8]{ABPS16}. The sign $(-1)^{n-1}$ in $\mu_z$ is the one that the central character requires. In the basis $v,sv,\dots,s^{n-1}v$ of $\Ind\mu_z$, where $W_{K_n}$ acts on $v$ by $\mu_z$, the element $s$ acts by the companion matrix of $X^n-\mu_z(s^n)$, so
\[
 \det(\Ind\mu_z)(s)=(-1)^{n-1}\mu_z(s^n)=z,
\]
and the local Langlands correspondence matches this determinant, through $a_K$, with the value $z$ of the central character of $\pi_z$ at $\varpi=a_K(s)$.

\subsection{Fields of rationality and image algebras}\label{subsec:image-algebras}

We now compute the fields of rationality of $\rho_z$ and $\tau_z$ (Definition~\ref{def:rationality}) and their image algebras. The field of rationality is $E(z)$ (Proposition~\ref{prop:image-algebras}). Adjoining $z$ to $E$, however, can enlarge the intersection with $L$, and then the cyclic extension of $E(z)$ generated by $L$ has degree smaller than $n$. The following fields and degrees measure this.

\begin{definition}\label{def:coefficient-fields}
For $z\in\Qbar^\times$ we put
\[
 E_z=E(z),\qquad L_z=LE_z,\qquad h_z=[L\cap E_z:E],\qquad d_z=n/h_z,
\]
where the compositum and the intersection are taken in $\Qbar$. Restriction identifies $\Gal(L_z/E_z)$ with $\Gal(L/L\cap E_z)$, the subgroup of $\Gal(L/E)$ generated by $\gamma^{h_z}$, and we let $\gamma_z\in\Gal(L_z/E_z)$ be the generator that restricts to $\gamma^{h_z}$ on $L$.
\end{definition}

Thus $[L_z:E_z]=d_z$, and for $z=1$ we have $E_1=E$, $L_1=L$, $h_1=1$, $d_1=n$ and $\gamma_1=\gamma$. When $h_z>1$, a suitable element of $\OO_{K_n}^\times$ acts on $\rho_z$, and a suitable element of $I_K$ acts on $\tau_z$, by an operator that generates over $E_z$ an algebra isomorphic to $L\otimes_EE_z$, a product of $h_z$ copies of $L_z$. The image algebra is then a crossed product over this product of fields, and the following lemma reduces such crossed products to cyclic algebras of smaller degree.

\begin{lemma}\label{lem:corners}
Let $\delta$ be a generator of $\Gal(L/E)$, let $F'$ be a field extension of $E$ and let $b\in F'^\times$. Put $T=L\otimes_EF'$, on which $\delta$ acts through the first factor, and let $A$ be the $F'$-algebra with $T$-basis $1,u,\dots,u^{n-1}$ and relations $uxu^{-1}=\delta(x)$ for $x\in T$ and $u^n=b$. Let $e$ be a primitive idempotent of $T$ and $g\geq1$ the least integer with $\delta^g(e)=e$. Then $g$ divides $n$, the field $L'=eT$ is a cyclic extension of $F'$ of degree $n/g$ whose Galois group is generated by the restriction $\delta'$ of $\delta^g$, and
\[
 A\cong M_g\bigl((L'/F',\delta',b)\bigr).
\]
If $b\in E$, then $A=(L/E,\delta,b)\otimes_EF'$, and
\[
 \res_{F'/E}[(L/E,\delta,b)]=[A]=[(L'/F',\delta',b)].
\]
\end{lemma}

The primitive idempotents of $T$ form one orbit under $\delta$. The corner $eAe$ at one of them is a cyclic algebra over $F'$, and the elements $u^ieu^{-j}$ are matrix units that connect the conjugate corners.

\begin{proof}
Since $L/E$ is separable, $T$ is a finite product of fields, and $\delta$ permutes its primitive idempotents. The sum over an orbit is fixed by $\delta$, hence lies in the field $T^\delta=F'$, so $\delta$ acts transitively, the primitive idempotents are $e,\delta(e),\dots,\delta^{g-1}(e)$, and $g$ divides $n$. The group of order $n/g$ generated by $\delta^g$ acts faithfully on $L'=eT$ (an element of $\Gal(L/E)$ trivial on $eT$ is trivial on every $\delta^i(e)T$), with fixed ring $eF'$ (for $x\in eT$ fixed by $\delta^g$, the sum $\sum_{0\leq i<g}\delta^i(x)$ lies in $F'$). By Artin's theorem $L'/F'$ is cyclic of degree $n/g$, with Galois group generated by $\delta'$.

Since $u^ieu^{-i}=\delta^i(e)$, the product $eu^ie=e\delta^i(e)u^i$ vanishes unless $g$ divides $i$. Hence
\[
 eAe=\bigoplus_{0\leq j<n/g}L'v^j,\qquad v=eu^g=u^ge,
\]
with $vxv^{-1}=\delta'(x)$ for $x\in L'$ and $v^{n/g}=be$, so the simple algebra $(L'/F',\delta',b)$ maps isomorphically onto $eAe$, both having dimension $(n/g)^2$ over $F'$. Finally, the elements $e_{ij}=u^ieu^{-j}$, $0\leq i,j<g$, form a system of matrix units with $\sum_ie_{ii}=1$, so $x\mapsto(e_{0i}xe_{j0})_{i,j}$ is an isomorphism from $A$ onto $M_g(eAe)$. When $b\in E$, the presentation of $A$ is that of $(L/E,\delta,b)\otimes_EF'$, which gives the last assertion.
\end{proof}

Over $F'=E_z$, Lemma~\ref{lem:corners} gives the matrix presentations of the image algebras (part~(i) of the next corollary), and over the completions of $E_z$ it describes how these algebras localize (parts~(ii) and~(iii)).

\begin{corollary}\label{cor:localization}
Let $F=E_z$ for some $z\in\Qbar^\times$, let $\delta$ generate $\Gal(L/E)$, let $b\in F^\times$, and let $A$ be the $F$-algebra of Lemma~\ref{lem:corners} for $F'=F$.
\begin{enumerate}
\item[(i)] If $\delta=\gamma$, resp.\ $\delta=\gamma^{-1}$, then
\[
 A\cong M_{h_z}\bigl((L_z/E_z,\gamma_z,b)\bigr),\quad\text{resp.}\quad A\cong M_{h_z}\bigl((L_z/E_z,\gamma_z^{-1},b)\bigr).
\]
\item[(ii)] Let $v$ be a place of $F$, let $L'$ be a factor of the product of fields $L\otimes_EF_v$, and let $g$ be the number of these factors. Then $L'$ is generated over $F_v$ by the image of $L$, the integer $g$ is the least positive integer such that $\delta^g$ maps $L'$ onto itself, the restriction $\delta'$ of $\delta^g$ to $L'$ generates $\Gal(L'/F_v)$, and
\[
 \res_{F_v/F}[A]=[(L'/F_v,\delta',b)]\qquad\text{in }\Br(F_v).
\]
If $F=E$, then $L'$ is the completion $L_{\tilde v}$ of $L$ at a place $\tilde v$ above $v$, and $g$ is the number $g_v$ of places of $L$ above $v$.
\item[(iii)] Let $v$ be a finite place of $F$ above a prime at which $L$ is unramified, and suppose that $b$ is a unit at $v$. Then $\res_{F_v/F}[A]=0$.
\end{enumerate}
\end{corollary}

\begin{proof}
(i) Here $T=L\otimes_EE_z$, and the kernel of the surjective multiplication map $T\to L_z$ is $(1-e)T$ for a primitive idempotent $e$ with $eT\cong L_z$. Since this kernel determines $e$, the element $\gamma^i$ fixes $e$ exactly when it extends to an element of $\Gal(L_z/E_z)$, that is, when $h_z$ divides $i$ (Definition~\ref{def:coefficient-fields}). Hence $g=h_z$ for $\delta=\gamma^{\pm1}$, the restriction $\delta'$ of $\gamma^{\pm h_z}$ to $eT\cong L_z$ is $\gamma_z^{\pm1}$, and Lemma~\ref{lem:corners} yields (i).

(ii) The algebra $A\otimes_FF_v$ is the algebra of Lemma~\ref{lem:corners} for $F'=F_v$, with the same $\delta$ and $b$, and its class is $\res_{F_v/F}[A]$. Each factor of $L\otimes_EF_v$ is generated over $F_v$ by the image of $L$, and $\delta^i$ maps the factor $L'=e(L\otimes_EF_v)$ onto $\delta^i(e)(L\otimes_EF_v)$. Since $\delta$ permutes the primitive idempotents transitively, the least $i\ge1$ with $\delta^i(L')=L'$ is their number $g$, and the lemma gives $A\otimes_FF_v\cong M_g\bigl((L'/F_v,\delta',b)\bigr)$. For $F=E$ the factors are the completions $L_{\tilde v}$ at the places $\tilde v$ above $v$, which $\Gal(L/E)$ permutes as it permutes these places, so $g=g_v$.

(iii) By (ii), $\res_{F_v/F}[A]$ is the class of $(L'/F_v,\delta',b)$, where $L'$ is generated by the image of $L$ and hence unramified over $F_v$, and this algebra splits by Lemma~\ref{lem:cyclic-split-invariant}(ii), since $b$ is a unit at $v$.
\end{proof}

\begin{proposition}\label{prop:image-algebras}
Let $z\in\Qbar^\times$. The representations $\rho_z$ and $\tau_z$ both have field of rationality $E_z$, and
\[
 A_{\rho_z}\cong M_{h_z}\bigl((L_z/E_z,\gamma_z,z)\bigr),\qquad A_{\tau_z}\cong M_{h_z}\bigl((L_z/E_z,\gamma_z^{-1},az)\bigr).
\]
In particular $A_{\rho_1}\cong(L/E,\gamma,1)$ and $A_{\tau_1}\cong(L/E,\gamma^{-1},a)$.
\end{proposition}

We treat $\rho_z$ and $\tau_z$ in parallel. The operator of an element with image $\bar t$ in the common quotient $\F_{q^n}^\times$ generates over $E_z$ a copy $T$ of $L\otimes_EE_z$, and the operator $u$ of $\Pi$ or of $s$, which permutes the lines of Lemma~\ref{lem:shape}(ii) cyclically, normalizes $T$. The image algebra is the crossed product $\bigoplus_{0\leq i<n}Tu^i$, and Corollary~\ref{cor:localization}(i) turns this decomposition into the matrix presentation.

\begin{proof}
Let $\sigma\in\Aut(\C)$. By Lemma~\ref{lem:shape}(i), the representations $\rho_z^\sigma$ and $\tau_z^\sigma$ are induced from $\sigma\circ\chi_z$ and $\sigma\circ\mu'_z$, which agree with $\sigma\circ\theta$ on $\OO_{K_n}^\times$, respectively $I_K$, and take the values $\sigma(z)$ at $\varpi$, respectively $a\sigma(z)$ at $s^n$. By parts (ii) and (iii) of Lemma~\ref{lem:shape}, $\rho_z^\sigma\cong\rho_z$ if and only if $\sigma\circ\theta=\theta^{q^i}$ for some $i$ and $\sigma(z)=z$ (conversely, $\sigma\circ\chi_z$ is then the conjugate of $\chi_z$ by a power of $\Pi$), and the same holds for $\tau_z$, with $s$ in place of $\Pi$. Writing $\sigma(\zeta_m)=\zeta_m^c$, we have $\sigma\circ\theta=\theta^{q^i}$ if and only if $c\equiv q^i\pmod m$, so the two conditions read
\[
 \sigma|_L\in\langle\gamma\rangle=\Gal(L/E),\qquad \sigma(z)=z,
\]
that is, $\sigma$ fixes $E(z)$. So the field of rationality is $E_z$.

Now let $r$ be $\rho_z$ or $\tau_z$. Choose $t\in\OO_{K_n}^\times$ with residue class $\bar t$ if $r=\rho_z$, and $t\in I_K$ with image $\bar t$ in $\F_{q^n}^\times$ if $r=\tau_z$. By Lemma~\ref{lem:shape}(ii), the operator $r(t)$ is diagonalizable with the eigenvalues $\theta^{q^i}(\bar t)=\zeta_m^{q^i}$, $0\leq i<n$, each of multiplicity one. These are the roots of the minimal polynomial $P$ of $\zeta_m$ over $E$, so $P$ is also the minimal polynomial of $r(t)$ over $E_z$, and the $E_z$-algebra $T$ generated by $r(t)$ is isomorphic to $E_z[X]/(P)\cong L\otimes_EE_z$, with $r(t)$ corresponding to $\zeta_m\otimes1$.

Let $u$ be the operator of $\Pi$ on $\rho_z$, or of $s$ on $\tau_z$. Then
\[
 u\,x\,u^{-1}=\delta(x)\quad(x\in T),\qquad u^n=b,\qquad\text{where}\quad (u,b,\delta)=
 \begin{cases}
 (\rho_z(\Pi),\,z,\,\gamma) & \text{if } r=\rho_z,\\
 (\tau_z(s),\,az,\,\gamma^{-1}) & \text{if } r=\tau_z,
 \end{cases}
\]
and $\delta$ acts on $T\cong L\otimes_EE_z$ through the first factor. Indeed, conjugation by $\rho_z(\Pi)$ sends $\rho_z(t)$ to $\rho_z(t)^q$, since conjugation by $\Pi$ induces $x\mapsto x^q$ on $\F_{q^n}$ by \eqref{eq:Pi} and $1+\mathfrak p_{K_n}\subset U_D^1$ acts trivially. Conjugation by $\tau_z(s)$ sends $\tau_z(t)$ to $\tau_z(t)^{q'}$, where $qq'\equiv1\pmod m$, since $s$ acts on tame inertia as the $q^{-1}$-th power, the inverse of the action of $\Pi$ (\S\ref{subsec:depth-zero-reps}). The relation $u^n=b$ is Lemma~\ref{lem:shape}(iii).

In a basis of eigenvectors of $r(t)$ the elements of $T$ are diagonal matrices, and $u^i$ carries each eigenline to the $i$-th following one. So the subspaces $Tu^i$, $0\leq i<n$, lie in distinct cyclic diagonals, and their sum is direct, of dimension $n^2$ over $E_z$. Each $Tu^i$ lies in $A_r$, since $r(t)$ and $u$ do and $A_r$ is an algebra. By Lemma~\ref{lem:trace-span}, applied with $F=E_z$, the field of rationality, the algebra $A_r$ also has dimension $n^2$, and therefore $A_r=\bigoplus_{i<n}Tu^i$. Thus $A_r$ is the algebra $A$ of Lemma~\ref{lem:corners} for $F'=E_z$, the generator $\delta$ and the parameter $b$, and Corollary~\ref{cor:localization}(i) yields both isomorphisms.
\end{proof}

\subsection{The case \texorpdfstring{$z=1$}{z=1}}\label{subsec:central-value-one}

\begin{theorem}\label{thm:central-value-one}
The representation $\rho_1$ has a model over $E$, so that $o_E(\rho_1)=0$. The class $o_E(\tau_1)=[(L/E,\gamma^{-1},a)]$ has order at most $2$, so it also equals $[(L/E,\gamma,a)]$, and it vanishes if $n$ is odd. The local invariants of $o_E(\tau_1)$ are as follows.
\begin{enumerate}
\item[(i)] At each place $v$ of $E$ above $p$,
\[
 \inv_v o_E(\tau_1)=\frac{f(n-1)}2=\frac{n-1}2\,[E_v:\Q_p]\quad\text{in }\Q/\Z.
\]
\item[(ii)] At each finite place not above $p$, including the places above the primes dividing $m$, the invariant is $0$.
\item[(iii)] The field $E$ has a real place if and only if $-1\in\langle q\rangle\subset(\Z/m\Z)^\times$. In that case $n$ is even, and the invariant at each real place is $1/2$.
\item[(iv)] At each complex place the invariant is $0$.
\end{enumerate}
Consequently
\begin{equation}\label{eq:central-value-one}
 o_E(\tau_1)=(n-1)\res_{E/\Q}[B_{p,\infty}]\quad\text{in }\Br(E).
\end{equation}
\end{theorem}

We determine the class of $(L/E,\gamma^{-1},a)$ place by place and then compare it with $(n-1)\res_{E/\Q}[B_{p,\infty}]$. At a place $v$ of $E$, Corollary~\ref{cor:localization}(ii) represents the restriction of the class to $E_v$ by a cyclic algebra for $L_{\tilde v}/E_v$ with parameter $a$, where $L_{\tilde v}$ is the completion of $L$ at a place $\tilde v$ above $v$. At the places above $p$ and at the real places, we read off the invariant of this algebra from the valuation and from the sign of $a$. At the other finite places we show that $a$ is a norm from $L_{\tilde v}$, so that the algebra splits by Lemma~\ref{lem:cyclic-split-invariant}(i). At such a place where $L/E$ is unramified, $a$ is a unit and hence a norm (Lemma~\ref{lem:cyclic-split-invariant}(ii)). The hard step is at the places above the primes dividing $m$, where $L/E$ can ramify and a unit need not be a norm, and there we use a congruence for $a$.

\begin{lemma}\label{lem:norm-above-m}
Let $\ell$ be a prime dividing $m$, let $v$ be a place of $E$ above $\ell$, and let $\tilde v$ be a place of $L$ above $v$. Then $a$ is a norm from $L_{\tilde v}^\times$ to $E_v^\times$.
\end{lemma}

\begin{proof}
If $n$ is odd, then $a=q^{n(n-1)/2}=\Nm_{L/E}(q^{(n-1)/2})$, so we may assume that $n$ is even. We shall show that $\Art_{E_v}(a)$ is trivial on $L_{\tilde v}$, and the point is a congruence for $a$ modulo the $\ell$-part of $m$. Write $m=m_\ell m'$, where $m_\ell$ is a power of $\ell$ and $\ell\nmid m'$. In $\Gal(L/\Q)=(\Z/m\Z)^\times$, let $H=\langle q\rangle=\Gal(L/E)$, and let $J$ be the decomposition group of $\tilde v$. Then $J\cap H$ is the decomposition group of $\tilde v$ in $L/E$, and $J/(J\cap H)$ has order $[E_v:\Q_\ell]$. We claim that
\[
 a^{[E_v:\Q_\ell]}\equiv1\pmod{m_\ell}.
\]
Since $q^n\equiv1\pmod m$, we have $a^2=(q^n)^{n-1}\equiv1\pmod m$, which settles the case where $[E_v:\Q_\ell]$ is even.

The case of odd $[E_v:\Q_\ell]$, however, is the hard one. If $m_\ell=2$, the congruence holds because $a$ is odd.

Otherwise $m_\ell>2$, and we use an element of order $2$ in the inertia group of $\tilde v$. Let $j\in(\Z/m\Z)^\times$ be the class with $j\equiv-1\pmod{m_\ell}$ and $j\equiv1\pmod{m'}$. It has order $2$, and it lies in the inertia group of $\tilde v$, which consists of the classes congruent to $1$ modulo $m'$ (the field $\Q(\zeta_{m'})$ is unramified at $\ell$, and $\Q(\zeta_{m_\ell})$ is totally ramified at $\ell$). So $j\in J$. The image of $j$ in $J/(J\cap H)$, a group of odd order, is trivial, and hence $j\in H$. The cyclic group $H$ has even order $n$, so its unique element of order $2$ is $q^{n/2}$, and therefore $j=q^{n/2}$. Thus
\[
 q^{n/2}\equiv-1\pmod{m_\ell},\qquad a=-(q^{n/2})^{n-1}\equiv1\pmod{m_\ell},
\]
since $n-1$ is odd. This proves the claim.

Now $a\in\Q_\ell$, so $\Nm_{E_v/\Q_\ell}(a)=a^{[E_v:\Q_\ell]}$. The compatibility of the Artin map with norms \cite[III, Lemma~3.2]{Mil20} shows that
\[
 \Art_{E_v}(a)|_{L_{\tilde v}}=\Art_{\Q_\ell}\bigl(a^{[E_v:\Q_\ell]}\bigr)\big|_{L_{\tilde v}}.
\]
As $\ell\neq p$, the number $a^{[E_v:\Q_\ell]}$ is a unit of $\Z_\ell$. So the Artin element of this unit fixes the roots of unity of order prime to $\ell$ and raises those of $\ell$-power order to the power $a^{-[E_v:\Q_\ell]}$ \cite[I, Example~3.13]{Mil20}, which is congruent to $1$ modulo $m_\ell$ by the claim. Since $L_{\tilde v}=\Q_\ell(\zeta_m)$, the element $\Art_{E_v}(a)$ is trivial on $L_{\tilde v}$, and $a$ is a norm from $L_{\tilde v}$ by local class field theory \cite[III, Theorem~3.4(b)]{Mil20}.
\end{proof}

\begin{proof}[Proof of Theorem~\ref{thm:central-value-one}]
By Proposition~\ref{prop:image-algebras}, $A_{\rho_1}\cong(L/E,\gamma,1)$, which is split, so $\rho_1$ has a model over $E$ by Lemma~\ref{lem:models}(i).

By the same proposition, $o_E(\tau_1)=[(L/E,\gamma^{-1},a)]$. Since $q^{n-1}\in E$, we have $a^2=q^{n(n-1)}=\Nm_{L/E}(q^{n-1})$, so $2\,o_E(\tau_1)=[(L/E,\gamma^{-1},a^2)]=0$ by Proposition~\ref{prop:cyclic}(i) and Lemma~\ref{lem:cyclic-split-invariant}(i). Hence $o_E(\tau_1)=-o_E(\tau_1)=[(L/E,\gamma,a)]$, again by Proposition~\ref{prop:cyclic}(i). If $n$ is odd, then $a=\Nm_{L/E}(q^{(n-1)/2})$ and the class vanishes by Lemma~\ref{lem:cyclic-split-invariant}(i).

To compute the local invariants, let $v$ be a place of $E$ and $\tilde v$ a place of $L$ above $v$. By Corollary~\ref{cor:localization}(ii), with $F=E$, $\delta=\gamma^{-1}$ and $b=a$,
\[
 \res_{E_v/E}\,o_E(\tau_1)=\bigl[(L_{\tilde v}/E_v,\gamma^{-g_v},a)\bigr],
\]
where $g_v$ is the number of places of $L$ above $v$ and $\gamma^{-g_v}$ stands for its restriction to $L_{\tilde v}$.

\emph{Places above $p$.} Since $p\nmid m$, the prime $p$ is unramified in $L$, and its decomposition group in $(\Z/m\Z)^\times$ is the subgroup $\langle p\rangle$, of order $\ord_m(p)$, the order of $p$ modulo $m$. This subgroup contains $q=p^f$, hence contains $\langle\gamma\rangle$. So for $v\mid p$ we have $g_v=1$, and the extension $L_{\tilde v}/E_v$ is unramified of degree $n$. Since $p^f$ has order $n$ modulo $m$,
\[
 [E_v:\Q_p]=\frac{\ord_m(p)}n=\gcd\bigl(\ord_m(p),f\bigr),
\]
and so $f=c\,[E_v:\Q_p]$ for an integer $c$ prime to $n$. The arithmetic Frobenius of $L_{\tilde v}/E_v$ is $\Frob=p^{[E_v:\Q_p]}$, because $E_v/\Q_p$ is unramified, and therefore
\[
 \gamma=q=p^{c[E_v:\Q_p]}=\Frob^c\qquad\text{on }L_{\tilde v}.
\]
Since $g_v=1$, the localization displayed above is the class of $(L_{\tilde v}/E_v,\Frob^{-c},a)$.

To apply Proposition~\ref{prop:cyclic}(ii) we pass to the generator $\Frob^{-1}$. If $u$ is the generator of $(L_{\tilde v}/E_v,\Frob^{-c},a)$ and $cc'\equiv1\pmod n$, then $u^{c'}$ acts on $L_{\tilde v}$ by $\Frob^{-1}$ and its powers form an $L_{\tilde v}$-basis, so
\[
 (L_{\tilde v}/E_v,\Frob^{-c},a)\cong(L_{\tilde v}/E_v,\Frob^{-1},a^{c'}).
\]
Now $a^{c'}=\pm p^k$ with $k=c'fn(n-1)/2$, the prime $p$ is a uniformizer of $E_v$, and the unit $\pm1$ contributes nothing (Lemma~\ref{lem:cyclic-split-invariant}(ii)), so, by Proposition~\ref{prop:cyclic},
\[
 \inv_v o_E(\tau_1)=-\frac kn=-\frac{c'f(n-1)}2\quad\text{in }\Q/\Z.
\]
So the invariant at $v$ comes from the $p$-adic valuation of $a$, which the half twist contributes. For odd $n$ this number and the two expressions in (i) are integers. For even $n$ the integers $c$ and $c'$ are odd, so $-c'f(n-1)/2\equiv f(n-1)/2=c[E_v:\Q_p](n-1)/2\equiv[E_v:\Q_p](n-1)/2$ modulo $\Z$. This proves (i).

\emph{Finite places not above $p$.} Let $v$ lie above a prime $\ell\neq p$. If $\ell\nmid m$, then $L/E$ is unramified at $v$ and $a$ is a unit at $v$, so the localization splits by Lemma~\ref{lem:cyclic-split-invariant}(ii). If $\ell\mid m$, then $a$ is a norm from $L_{\tilde v}$ by Lemma~\ref{lem:norm-above-m}, and the localization splits by Lemma~\ref{lem:cyclic-split-invariant}(i). This proves (ii).

\emph{Infinite places.} Complex conjugation acts on $L$ as the class of $-1$, which is not trivial since $m\geq3$. As $E/\Q$ is abelian, $E$ is either totally real or totally imaginary, and it is totally real if and only if complex conjugation fixes $E$, that is, if and only if $-1\in\langle q\rangle$. In that case the cyclic group $\langle\gamma\rangle$ contains the element $-1$ of order $2$, so $n$ is even, $-1=\gamma^{n/2}$ and $a=-q^{n(n-1)/2}<0$.

At a real place $v$ of $E$, the places of $L$ above $v$ are complex, with decomposition group generated by $\gamma^{n/2}$, so $g_v=n/2$, and $\gamma^{-g_v}$ acts on $L_{\tilde v}=\C$ by complex conjugation ($\mathrm{conj}(z)=\bar z$). So, by the formula for $\res_{E_v/E}\,o_E(\tau_1)$ at the start of the proof, the restriction of $o_E(\tau_1)$ to $E_v=\R$ is the class of $(\C/\R,\mathrm{conj},a)$. Since $-a>0$ is a norm from $\C$,
\[
 (\C/\R,\mathrm{conj},a)\cong(\C/\R,\mathrm{conj},-1),
\]
the algebra of Hamilton quaternions, whose invariant is $1/2$. At a complex place the Brauer group vanishes. This proves (iii) and (iv).

Finally, by Theorem~\ref{thm:local-brauer}(a), the class $(n-1)\res_{E/\Q}[B_{p,\infty}]$ has invariant $(n-1)[E_v:\Q_p]/2$ at each place $v$ above $p$, invariant $0$ at the other finite places and at the complex places, and invariant $(n-1)/2$ at each real place, which is $1/2$ modulo $\Z$ because $n$ is even when $E$ has a real place. These are the invariants of $o_E(\tau_1)$, so the two classes in \eqref{eq:central-value-one} are equal, because the localization map $\Br(E)\to\bigoplus_v\Br(E_v)$ is injective \cite[VII, Theorem~7.1, and IV, Theorem~3.14]{Mil20}.
\end{proof}

\begin{corollary}\label{cor:splitting}
The representation $\tau_1$ has a model over $E$ if and only if $n$ is odd, or $n$ and $f$ are both even and $E$ is totally imaginary. Otherwise the Schur index of $\tau_1$ is $2$.
\end{corollary}

\begin{proof}
By Lemma~\ref{lem:models}(i), $\tau_1$ has a model over $E$ exactly when all the invariants in Theorem~\ref{thm:central-value-one} vanish. They do for odd $n$. For even $n$ the invariant $f(n-1)/2\equiv f/2$ above $p$ vanishes exactly when $f$ is even, and each real place contributes $1/2$, which gives the criterion, since $E$ is totally real or totally imaginary. Otherwise $o_E(\tau_1)=\res_{E/\Q}[B_{p,\infty}]$ is the class of the nonsplit quaternion algebra $B_{p,\infty}\otimes_\Q E$, a division algebra, so the Schur index of $\tau_1$ is $2$.
\end{proof}

\begin{corollary}\label{cor:fields-of-definition}
The representations $\rho_1$, $\tau_1$ and $\pi_1$ have models over $L$. The representation $\pi_1$ has field of rationality $E$ and a model over $E$. In particular all three representations have models over every completion $L_{\tilde v}$ of $L$, and the classes $o_E(\rho_1)$ and $o_E(\tau_1)$ vanish in $\Br(L_{\tilde v})$.
\end{corollary}

\begin{proof}
The representations $\rho_1$ and $\tau_1$ are induced from $\chi_1$ and from the character $\mu'_1$ of Lemma~\ref{lem:shape}(i), which take values in $L$, since $\chi_1$ is $\theta$ on $\OO_{K_n}^\times$ and $1$ at $\varpi$, and $\mu'_1$ is $\theta$ on $I_K$ and $a\in\Q^\times$ at $s^n$. So they have models over $L$, and no square root of $q$ is needed, because the half twist enters $\mu'_1$ only through the rational value $q^{n(n-1)/2}$ at $s^n$.

For $\pi_1$ we use Suzuki's results, which concern representations whose central character has finite order \cite[Hypothesis~4.0.1]{Su23}. The central character of $\pi_1$ is that of $\rho_1$, which is trivial on $\varpi$ and factors through $\F_q^\times$ on $\OO_K^\times$, so it has finite order, and $\pi_1$, a Jacquet--Langlands correspondent with unitary central character, is square-integrable. The field of rationality of $\pi_1$ is that of $\rho_1$, because the Jacquet--Langlands correspondence commutes with the automorphisms of $\C$ (Dat's remark after \cite[Th\'eor\`eme~2.1.2]{Dat07}), and this field is $E$ by Proposition~\ref{prop:image-algebras}. So $\pi_1$ has a model over $E$ \cite[Lemma~4.0.4]{Su23}, hence over $L$ and over every $L_{\tilde v}$. Finally, by Lemma~\ref{lem:models}(i), the models of $\rho_1$ and $\tau_1$ over $L$ show that $o_E(\rho_1)$ and $o_E(\tau_1)$ vanish in $\Br(L)$, hence in $\Br(L_{\tilde v})$.
\end{proof}

Corollary~\ref{cor:fields-of-definition} shows why the identity \eqref{eq:intro-sum} is stated over the field of rationality. Over $L$, or over any completion $L_{\tilde v}$, both classes vanish, so restriction to these fields loses the invariants $\frac{n-1}2[E_v:\Q_p]$ of $\tau_1$ at the places above $p$.

\subsection{Arbitrary \texorpdfstring{$z$}{z}}\label{subsec:arbitrary}

For $z=1$ the algebra $A_{\rho_1}$ is split. For arbitrary $z$ both classes depend on $z$, but their sum is $(n-1)\res_{E_z/\Q}[B_{p,\infty}]$, which depends on $z$ only through the field $E_z$.

\begin{theorem}\label{thm:arbitrary-central-value}
Let $z\in\Qbar^\times$. Then
\begin{equation}\label{eq:arbitrary-central-value}
 o_{E_z}(\rho_z)+o_{E_z}(\tau_z)=(n-1)\res_{E_z/\Q}[B_{p,\infty}] \quad\text{in }\Br(E_z).
\end{equation}
Thus the invariant of the sum is $\frac{n-1}2[E_{z,v}:\Q_p]$ at each place $v$ above $p$, $\frac{n-1}2$ at each real place, and $0$ at every other place. If $z$ is a unit at a place $v$ above $p$, then
\[
 \inv_v o_{E_z}(\rho_z)=0,\qquad \inv_v o_{E_z}(\tau_z)=\frac{n-1}2\,[E_{z,v}:\Q_p]\quad\text{in }\Q/\Z.
\]
This applies at every place above $p$ when the central character of $\rho_z$ has finite order, which is the case if and only if $z$ is a root of unity.
\end{theorem}

\begin{proof}
By Proposition~\ref{prop:image-algebras} and Proposition~\ref{prop:cyclic}(i),
\[
 o_{E_z}(\rho_z)+o_{E_z}(\tau_z)=[(L_z/E_z,\gamma_z,z)]-[(L_z/E_z,\gamma_z,az)]=-[(L_z/E_z,\gamma_z,a)].
\]
Here the inverse generator on the Weil side produces the minus sign, and additivity in the parameter cancels $z$. Since $a\in E$, Lemma~\ref{lem:corners}, applied with $F'=E_z$, $\delta=\gamma$ and $b=a$, and Corollary~\ref{cor:localization}(i) show that $[(L_z/E_z,\gamma_z,a)]=\res_{E_z/E}[(L/E,\gamma,a)]$. By Theorem~\ref{thm:central-value-one}, $[(L/E,\gamma,a)]=(n-1)\res_{E/\Q}[B_{p,\infty}]$, a class of order at most $2$, and so the sum equals $(n-1)\res_{E_z/\Q}[B_{p,\infty}]$. At the finite places, the invariants of the sum follow from Theorem~\ref{thm:local-brauer}(a), since $B_{p,\infty}$ has invariant $\frac12$ at $p$ and invariant $0$ at the other primes. At a real place $v$ of $E_z$ we have $E_{z,v}=\R$, and the invariant is that of $(n-1)[B_{p,\infty}\otimes_\Q\R]$, namely $\frac{n-1}2$. At a complex place the Brauer group vanishes.

If $z$ is a unit at a place $v$ of $E_z$ above $p$, then $\inv_v o_{E_z}(\rho_z)=0$ by Proposition~\ref{prop:image-algebras} and Corollary~\ref{cor:localization}(i) and~(iii), since $L$ is unramified at $p\nmid m$, and the invariant of $o_{E_z}(\tau_z)$ at $v$ follows from \eqref{eq:arbitrary-central-value}. The central character of $\rho_z$ factors through $\F_q^\times$ on $\OO_K^\times$ and sends $\varpi$ to $z$, so it has finite order exactly when $z$ is a root of unity, which is a unit at every place.
\end{proof}

\begin{example}[Classes of order three with finite central character]
\label{ex:order-three}
This example shows that the classes $o_{E_z}(\rho_z)$ and $o_{E_z}(\tau_z)$ in Theorem~\ref{thm:arbitrary-central-value} can have order three when the central character has finite order, and then the inverse generator in Proposition~\ref{prop:image-algebras} changes the class. It also shows that the algebra $A_{\rho_z}$, which is split for the choice $z=1$ (Theorem~\ref{thm:central-value-one}), need not be split for other $z$.

Take $K=\Q_2$, $n=3$ and $\theta$ of order $7$, which is primitive since $2$ has order $3$ modulo $7$. Then $L=\Q(\zeta_7)$, $\gamma(\zeta_7)=\zeta_7^2$, the field $E$ is the fixed field of the squares $\{1,2,4\}$ in $(\Z/7\Z)^\times$, namely $\Q(\sqrt{-7})$, and $a=8$. Let $z=\zeta_3$ be a primitive cube root of unity. Then $E_z=\Q(\sqrt{-7},\zeta_3)$ has degree $4$ over $\Q$, and $L_z=\Q(\zeta_7,\zeta_3)$. The field $L\cap E_z$ lies between $E$ and $L$, and $[L:E]=3$ is prime, so it is $E$, because $[L:\Q]=6$ does not divide $[E_z:\Q]=4$. Hence $h_z=1$, $d_z=3$, and $\gamma_z$ sends $\zeta_7$ to $\zeta_7^2$ and fixes $\zeta_3$.

By Proposition~\ref{prop:image-algebras},
\[
 A_{\rho_z}\cong(L_z/E_z,\gamma_z,\zeta_3),\qquad A_{\tau_z}\cong(L_z/E_z,\gamma_z^{-1},8\zeta_3).
\]
Since $8=\Nm_{L_z/E_z}(2)$ is a norm, Proposition~\ref{prop:cyclic}(i) and Lemma~\ref{lem:cyclic-split-invariant}(i) show that
\[
 o_{E_z}(\tau_z)=-[(L_z/E_z,\gamma_z,8\zeta_3)]=-o_{E_z}(\rho_z).
\]
So the two classes are opposite, and their sum is $0=2\res_{E_z/\Q}[B_{p,\infty}]$, as \eqref{eq:arbitrary-central-value} requires.

We compute the invariants of $o_{E_z}(\rho_z)$, first at the places above $7$. The prime $7$ is ramified in $\Q(\sqrt{-7})$ and split in $\Q(\zeta_3)$, so $E_z$ has two places $v$ above $7$, each with $[E_{z,v}:\Q_7]=2$ and residue field $\F_7$, in which $\zeta_3$ reduces to $2$ at one place and to $4$ at the other. The ramification index of $7$ is $6$ in $L_z$ and $2$ in $E_z$, so $L_z/E_z$ is totally ramified of degree $3$ at $v$, and the localization of $A_{\rho_z}$ at $v$ is the cyclic algebra $(L_{z,\tilde v}/E_{z,v},\gamma_z,\zeta_3)$ with $L_{z,\tilde v}=\Q_7(\zeta_7)$.

By Lemma~\ref{lem:cyclic-split-invariant}(iii), the invariant of this algebra is $t/3$ if $\Art_{E_{z,v}}(\zeta_3)$ acts on $L_{z,\tilde v}$ as $\gamma_z^t$. Since $\zeta_3\in\Q_7$ and $[E_{z,v}:\Q_7]=2$, the compatibility with norms \cite[III, Lemma~3.2]{Mil20} shows that
\[
 \Art_{E_{z,v}}(\zeta_3)|_{L_{z,\tilde v}}=\Art_{\Q_7}(\zeta_3^2)|_{L_{z,\tilde v}},
\]
and $\Art_{\Q_7}(\zeta_3^2)$ sends $\zeta_7$ to $\zeta_7^c$ with $c\equiv(\zeta_3^2)^{-1}=\zeta_3$ modulo $7$ \cite[I, Example~3.13]{Mil20}. So $\Art_{E_{z,v}}(\zeta_3)$ acts on $L_{z,\tilde v}$ as $\gamma_z$ at the place where $\zeta_3\equiv2$ and as $\gamma_z^2$ at the place where $\zeta_3\equiv4$. Since $o_{E_z}(\tau_z)=-o_{E_z}(\rho_z)$, the invariants at the two places above $7$ are
\[
\begin{aligned}
 &\text{where }\zeta_3\equiv2\text{,} &\quad \inv_v o_{E_z}(\rho_z)&=\tfrac13,
 &\quad \inv_v o_{E_z}(\tau_z)&=\tfrac23,\\
 &\text{where }\zeta_3\equiv4\text{,} &\quad \inv_v o_{E_z}(\rho_z)&=\tfrac23,
 &\quad \inv_v o_{E_z}(\tau_z)&=\tfrac13.
\end{aligned}
\]
At every other finite place the class $o_{E_z}(\rho_z)=[(L_z/E_z,\gamma_z,\zeta_3)]$ vanishes, by parts (i) and (iii) of Corollary~\ref{cor:localization}, applied with $\delta=\gamma$ and $b=\zeta_3$, since $L=\Q(\zeta_7)$ is unramified away from $7$ and $\zeta_3$ is a unit, and so $o_{E_z}(\tau_z)=-o_{E_z}(\rho_z)$ vanishes there as well. The field $E_z$ has no real place, so the two classes have order three.
\end{example}

It remains to see that the representations $\pi_z$ exhaust depth zero (recall that $n\ge2$ in this section).

\begin{proposition}\label{prop:depth-zero-exhaustion}
Every depth-zero supercuspidal representation of $\GL_n(K)$ whose central character takes algebraic values is isomorphic to $\pi_z$ for some primitive character $\theta$ of $\F_{q^n}^\times$ and some $z\in\Qbar^\times$.
\end{proposition}

\begin{proof}
The Jacquet--Langlands correspondence preserves depth \cite[Corollary~2.8]{ABPS16}. On $D^\times$, depth zero is equivalent to having a nonzero $U_D^1$-fixed vector \cite[\S2.4 and Proposition~2.5]{ABPS16}, applied to $\OO_D$, the only hereditary order of $D$. As $U_D^1$ is normal, an irreducible representation with such a vector is trivial on $U_D^1$.

Prasad and Ramakrishnan classify the irreducible representations of $D^\times$ that are trivial on $U_D^1$ \cite[\S4, pp.~758--759]{PR12}. For a divisor $d$ of $n$, let $K_d\subseteq K_n$ be the unramified extension of $K$ of degree $d$. The representations of dimension $d$ correspond to the $\Gal(K_d/K)$-orbits of characters of $K_d^\times$ that are trivial on $1+\mathfrak p_{K_d}$ and whose $\Gal(K_d/K)$-conjugates are pairwise distinct (Prasad and Ramakrishnan call such characters regular). Since $\Gal(K_d/K)$ fixes $\varpi$, the last condition means that the restriction to $\OO_{K_d}^\times$ comes from a character of $\F_{q^d}^\times$ with a $q$-power orbit of $d$ elements. The Weil--Deligne parameter attached to a representation of dimension $d$ has the factor $\mathrm{sp}_{n/d}$ \cite[Theorem~4.4]{PR12}, and a supercuspidal representation has an irreducible Weil parameter and zero monodromy, so $d=n$.

For $d=n$ the representation attached to such a character $\chi$ of $K_n^\times$ is induced from the character of $K_n^\times U_D^1$ that extends $\chi$ and is trivial on $U_D^1$. Its central character, which is that of the corresponding representation of $\GL_n(K)$, sends $\varpi$ to $\chi(\varpi)$ and takes roots of unity as values on $\OO_K^\times$, so it takes algebraic values if and only if $z=\chi(\varpi)$ lies in $\Qbar^\times$, and then the representation is $\rho_z$ for the primitive character $\theta$ that $\chi$ induces on $\F_{q^n}^\times$.
\end{proof}

We can now prove Conjecture~\ref{conj:intro-sum} in depth zero.

\begin{corollary}\label{cor:depth-zero-relation}
Conjecture~\ref{conj:intro-sum} holds for every depth-zero supercuspidal representation with algebraic-valued central character. For every $z\in\Qbar^\times$, the representations $\pi_z$, $\rho_z$ and $\tau_z$ have common field of rationality $E_z$, and
\[
 o_{E_z}(\rho_z)+o_{E_z}(\tau_z)=(n-1)\res_{E_z/\Q}[B_{p,\infty}].
\]
\end{corollary}

\begin{proof}
For $n=1$ the conjecture holds (\S\ref{subsec:intro-conjecture}). For $n\ge2$, Proposition~\ref{prop:depth-zero-exhaustion} reduces the assertion to $\pi_z$. Proposition~\ref{prop:depth-zero-parameter} identifies its normalized Langlands parameter with $\tau_z$, and Proposition~\ref{prop:image-algebras} shows that $\Q(\rho_z)=\Q(\tau_z)=E_z$. Lemma~\ref{lem:common-rationality-field} applies to $\pi_z$, whose central character takes algebraic values, so $\Q(\pi_z)=\Q(\rho_z)=E_z$ as well. The identity then follows from Theorem~\ref{thm:arbitrary-central-value}.
\end{proof}

Section~\ref{sec:individual} takes up the question of how far apart the two classes can be in general.

\section{Schur indices and a family with unbounded unequal indices}\label{sec:individual}

The identity \eqref{eq:intro-sum}, wherever it holds, forces the Schur indices of $\rho$ and $\tau$ to have the same odd part and to differ by at most a factor of two (\S\ref{subsec:index-relations}),
\[
 \frac{\ind A_\tau}{\ind A_\rho}\in\bigl\{\tfrac12,1,2\bigr\}.
\]
Within these limits, however, the indices can be unequal and unbounded. In \S\ref{subsec:family} we construct simple supercuspidal representations whose two classes have indices $n$ and $n/2$, by representing both classes by cyclic algebras for one cyclic extension of the field of rationality.

Throughout this section, $\pi$ is a supercuspidal representation of $\GL_n(K)$ whose central character takes algebraic values, and $\rho=\JL(\pi)$, the normalized parameter $\tau$ and the field $E$ are as in Section~\ref{sec:intro}, so that $E=\Q(\pi)=\Q(\rho)=\Q(\tau)$ is a number field (Lemma~\ref{lem:common-rationality-field}).

\subsection{Consequences for Schur indices}\label{subsec:index-relations}

Recall that the image algebras $A_\rho$ and $A_\tau$ of Definition~\ref{def:brauer-class} are central simple $E$-algebras, and that $A_\tau$ has degree $n$. For a central simple $E$-algebra $A$ and a place $v$ of $E$ we write $\ind A$ for the index of $A$ and $\ind_v A$ for the index of $A\otimes_E E_v$. Over a local field, including $\R$ and $\C$, and over a number field, the index of a central simple algebra equals the order of its class \cite[Theorems~1.5.34(iv) and~1.5.36(iv)]{Poo17}. Since a class in $\Br(E)$ is determined by its local invariants \cite[VII, Theorem~7.1]{Mil20}, it follows that $\ind A$ is the least common multiple of the indices $\ind_v A$, and $\ind_v A$ is the order of $\inv_v[A]$ in $\Q/\Z$.

The two propositions of this subsection rest on the following elementary fact about orders in $\Q/\Z$.

\begin{lemma}\label{lem:half-shift}
Let $x,y\in\Q/\Z$ have orders $m$ and $m'$, and suppose that $x+y\in\{0,\tfrac12\}$. If $x+y=0$, then $m'=m$. If $x+y=\tfrac12$, then
\[
 m'=\begin{cases}
 2m, & m\text{ odd},\\
 m/2, & m\equiv2\pmod4,\\
 m, & 4\mid m.
 \end{cases}
\]
In every case $m$ and $m'$ have the same odd part, and $m'/m\in\{\tfrac12,1,2\}$.
\end{lemma}

\begin{proof}
If $x+y=0$, then $y=-x$ has the order of $x$. Suppose that $x+y=\frac12$, and write $x=c/m$ with $\gcd(c,m)=1$, so that $y=(m-2c)/(2m)$. If $m$ is odd, then $m-2c$ is odd and prime to $m$, so $y$ has order $2m$. If $m=2k$ with $k$ odd, then $c$ is odd, so $k-c$ is even and $y=\frac{(k-c)/2}{k}$, where $(k-c)/2$ is prime to $k$ because $k-c$ is. Hence $y$ has order $k=m/2$. If $4\mid m$, then $c$ is odd and $y=(m/2-c)/m$ with $m/2-c$ odd, and an odd prime dividing both $m$ and $m/2-c$ would divide $c$. Hence $y$ has order $m$.
\end{proof}

At a place $v$ of $E$ above $p$, restriction multiplies local invariants by the local degree \cite[III, \S2, (29)]{Mil20}, and $\inv_p\Bpinf=\frac12$, so \eqref{eq:intro-sum} implies 
\begin{equation}\label{eq:local-sum}
 \inv_v o_E(\rho)+\inv_v o_E(\tau)=\frac{n-1}{2}\,[E_v:\Q_p] \qquad\text{in }\Q/\Z.
\end{equation}

\begin{proposition}\label{prop:local-indices}
Let $v$ be a place of $E$ above $p$, put $m_\rho=\ind_v A_\rho$ and $m_\tau=\ind_v A_\tau$, and assume \eqref{eq:local-sum} at $v$. Then $m_\rho$ and $m_\tau$ divide $n$. If $(n-1)[E_v:\Q_p]$ is even, then $\inv_v o_E(\tau)=-\inv_v o_E(\rho)$ and $m_\tau=m_\rho$. If $(n-1)[E_v:\Q_p]$ is odd, then $m_\tau$ equals $2m_\rho$, $m_\rho/2$ or $m_\rho$ according as $m_\rho$ is odd, $m_\rho\equiv2\pmod4$ or $4\mid m_\rho$.
\end{proposition}

\begin{proof}
The indices $m_\rho$ and $m_\tau$ are the orders of $x=\inv_v o_E(\rho)$ and $y=\inv_v o_E(\tau)$, and $x+y$ is $0$ or $\frac12$ according as $(n-1)[E_v:\Q_p]$ is even or odd. So Lemma~\ref{lem:half-shift} determines $m_\tau$ from $m_\rho$ as stated. Since $A_\tau\otimes_E E_v$ has degree $n$, its index $m_\tau$ divides $n$. If $(n-1)[E_v:\Q_p]$ is even, then $m_\rho=m_\tau$ divides $n$. If $(n-1)[E_v:\Q_p]$ is odd, then $n$ is even, so $n$ kills $\frac12$ and $y$, hence also $x=\frac12-y$, and $m_\rho$ divides $n$.
\end{proof}

\begin{proposition}\label{prop:global-indices}
Assume \eqref{eq:intro-sum}. Then $\ind A_\rho$ and $\ind A_\tau$ divide $n$ and have the same odd part, and $\ind A_\tau/\ind A_\rho\in\{\tfrac12,1,2\}$. If $n$ is odd, then $o_E(\tau)=-o_E(\rho)$ and $\ind A_\tau=\ind A_\rho$.
\end{proposition}

\begin{proof}
Since $A_\tau$ has degree $n$, its index divides $n$, so $n$ kills $o_E(\tau)$. The right side of \eqref{eq:intro-sum} is also killed by $n$, because $2[\Bpinf]=0$ and $n(n-1)$ is even. Hence $n$ kills $o_E(\rho)$, and since the index of a class in $\Br(E)$ is its order, both indices divide $n$. At every place $v$ of $E$ the local invariant of $(n-1)\res_{E/\Q}[\Bpinf]$ lies in $\{0,\frac12\}$, so by \eqref{eq:intro-sum} and Lemma~\ref{lem:half-shift} the local indices $\ind_v A_\rho$ and $\ind_v A_\tau$ have the same odd part, and their $2$-adic valuations differ by at most one. The global indices are the least common multiples of the local ones, so the same holds for them, which is the second assertion. If $n$ is odd, then $n-1$ is even and the right side of \eqref{eq:intro-sum} vanishes, so $o_E(\tau)=-o_E(\rho)$.
\end{proof}

\begin{remark}\label{rem:cases-occur}
The two cases of Proposition~\ref{prop:local-indices} with unequal indices both occur. If $\pi$ has depth zero and central character of finite order, then $\inv_v o_E(\rho)=0$ at every place $v$ above $p$ (Theorem~\ref{thm:arbitrary-central-value} and Proposition~\ref{prop:depth-zero-exhaustion}), so $m_\rho=1$, and $m_\tau=2$ when $(n-1)[E_v:\Q_p]$ is odd. The family of Theorem~\ref{thm:family} has $(n-1)[E_v:\Q_p]$ odd, $m_\rho=n$ and $m_\tau=n/2$ at every place above $p$, with $n=p-1$ unbounded, and the global indices of its two classes are also $n$ and $n/2$.
\end{remark}

\subsection{A family with unbounded unequal indices}\label{subsec:family}

In this subsection $p\ge7$ is a prime with $p\equiv3\pmod 4$, and
\[
 n=p-1,\qquad q=p^n.
\]
Let $K$ be the unramified extension of $\Q_p$ of degree $n$, with residue field $\F_q$, fix a uniformizer $\varpi$ of $K$, and let $D$ be the central division algebra over $K$ of invariant $1/n$. Fix a generator $\kappa$ of $\F_p^\times$ and a generator $b$ of $\F_q^\times$ with $b^{(q-1)/n}=\kappa$ (the map $x\mapsto x^{(q-1)/n}$ carries the generators of $\F_q^\times$ onto those of $\F_p^\times$), and identify $\kappa$ and $b$ with their Teichm\"uller representatives in $\mu_{q-1}\subset K$. Let $\chi\colon \F_q^\times\to\C^\times$ and $\psi\colon \F_q\to\C^\times$ be the characters
\[
 \chi(b)=\zeta_n,\qquad \psi(x)=\zeta_p^{\Tr_{\F_q/\F_p}(x)},
\]
where $\zeta_N$ denotes a fixed primitive $N$-th root of unity in $\C$. Let $\pi=\pi_{1,\chi,1}$ be the simple supercuspidal representation of $\GL_n(K)$ defined by Imai and Tsushima \cite[\S5.1]{IT20} for the uniformizer $\varpi$ and the parameters $\zeta=1$ and $c=1$ of their notation.

We take the residue degree of $K$ to be $n$ for two reasons. It makes $\chi$ trivial on $\F_p^\times$, as we check before Lemma~\ref{lem:weil-side}, and since $n$ is even it makes $q=p^n$ a square, so that the half twist takes rational values and $\tau$ has the same image algebra as $\varphi_\pi$ (end of the proof of Lemma~\ref{lem:weil-side}).

Finally, put
\[
 E=\Q(\zeta_n),\qquad E'=E(\zeta_p).
\]
Since $E$ is unramified at $p\nmid n$ and $\Q(\zeta_p)$ is totally ramified at $p$, the field $E\cap\Q(\zeta_p)$ is $\Q$, and $E'/E$ is cyclic of degree $p-1=n$. Let $\gamma_\kappa$ be the generator of $\Gal(E'/E)$ with $\gamma_\kappa(\zeta_p)=\zeta_p^\kappa$.

\begin{theorem}\label{thm:family}
The representations $\pi$, $\rho$ and $\tau$ have the following properties.
\begin{enumerate}
\item[(i)] The central character of $\pi$ has finite order, $\pi$ is not self-dual, and $\Q(\pi)=\Q(\rho)=\Q(\tau)=E$.
\item[(ii)] The Brauer classes are
\[
 o_E(\rho)=[(E'/E,\gamma_\kappa,\zeta_n^{-1})],\qquad o_E(\tau)=[(E'/E,\gamma_\kappa,-\zeta_n)].
\]
\item[(iii)] The prime $p$ splits completely in $E$. For $j\in(\Z/n\Z)^\times$ let $v_j$ be the place of $E$ above $p$ at which $\zeta_n\equiv \kappa^j$. Then
\[
 \inv_{v_j}o_E(\rho)=\frac jn,\qquad \inv_{v_j}o_E(\tau)=-\frac jn-\frac12,
\]
and both classes have invariant $0$ at every other place of $E$.
\item[(iv)] At every place of $E$ above $p$ the local Schur indices of $\rho$ and $\tau$ are $n$ and $n/2$, and $\ind A_\rho=n$, $\ind A_\tau=n/2$.
\item[(v)] The sum of the two classes is $o_E(\rho)+o_E(\tau)=\res_{E/\Q}[\Bpinf]$, so \eqref{eq:intro-sum} holds.
\end{enumerate}
\end{theorem}

Theorem~\ref{thm:intro-family} follows from parts (i) and (iii)--(v) of Theorem~\ref{thm:family}, whose part~(iii) fixes the labelling $v_j$ of the places above $p$ by the congruence $\zeta_n\equiv\kappa^j$. The Schur indices in Theorem~\ref{thm:family} are unbounded, because there are infinitely many primes $p\equiv3\pmod4$ (Euclid's argument applies to $4p_1\cdots p_r-1$).

\begin{example}\label{ex:rank-six}
The smallest member of the family already has unequal indices, $6$ and $3$, and the invariants of the two classes add up to $\frac12$ at both places above $7$. Take $p=7$, so that $n=6$ and $q=7^6$. Take $\kappa=3$, a generator of $\F_7^\times$, and a generator $b$ of $\F_q^\times$ with $b^{(q-1)/6}=3$. Then $E=\Q(\zeta_6)$ and $E'=\Q(\zeta_6,\zeta_7)$, and $\gamma_3$ fixes $\zeta_6$ and sends $\zeta_7$ to $\zeta_7^3$. By Theorem~\ref{thm:family}(ii), the two classes are $o_E(\rho)=[(E'/E,\gamma_3,\zeta_6^{-1})]$ and $o_E(\tau)=[(E'/E,\gamma_3,-\zeta_6)]$. The minimal polynomial $X^2-X+1$ of $\zeta_6$ has the roots $3=3^1$ and $5=3^5$ modulo $7$, so the prime $7$ splits in $E$ into the places $v_1$, where $\zeta_6\equiv3$, and $v_5$, where $\zeta_6\equiv5$, both with completion $\Q_7$. By Theorem~\ref{thm:family}(iii), the invariants of $o_E(\rho)$ and of $o_E(\tau)$ at $v_1$ and at $v_5$ form the ordered pairs
\[
 \bigl(\inv_{v_1}o_E(\rho),\,\inv_{v_1}o_E(\tau)\bigr)=\bigl(\tfrac16,\,\tfrac13\bigr),\qquad \bigl(\inv_{v_5}o_E(\rho),\,\inv_{v_5}o_E(\tau)\bigr)=\bigl(\tfrac56,\,\tfrac23\bigr).
\]
All other invariants vanish, so the Schur indices are $6$ and $3$. Each pair adds up to $\frac12$ in $\Q/\Z$, so the sum has invariant $\frac12$ at both places above $7$, as \eqref{eq:intro-sum} requires. The half twist multiplies the operator of $s$ by the rational number $q^{5/2}=7^{15}$, so $\tau$ and $\varphi_\pi$ have the same image algebra (see the end of the proof of Lemma~\ref{lem:weil-side}).
\end{example}

We obtain the two cyclic algebras of Theorem~\ref{thm:family}(ii) from the description of $\rho$ and $\varphi_\pi$ by Imai and Tsushima, as the image algebra $A_\tau$ itself on the Weil side (Lemma~\ref{lem:weil-side}) and as a corner of $A_\rho$ on the side of $D^\times$ (Lemma~\ref{lem:division-side}). We recall the data of \cite[\S5.1]{IT20} for $\zeta=c=1$, which we use with complex coefficients. Let $\alpha$ be a root of $X^n-\varpi$ and $K_\alpha=K(\alpha)$. Imai and Tsushima obtain $\varphi_\pi$ by induction from a character of $W_{K_\alpha}$, and $\rho$ by induction from a character of the subgroup $K_\alpha^\times U_D^1$ of $D^\times$, for an embedding of $K_\alpha$ in $D$. They deduce this description from the essentially tame correspondences of Bushnell and Henniart \cite{BH05,BH11J} (their Jacquet--Langlands correspondence goes from $D^\times$ to $\GL_n(K)$, so their identity $\JL(\rho_{1,\chi,1})=\pi_{1,\chi,1}$ is our $\rho=\JL(\pi)$). Since $\mu_n\subset K$ and $p\nmid n$, the extension $K_\alpha/K$ is cyclic, totally and tamely ramified of degree $n$. We write $U^i_{K_\alpha}=1+\alpha^i\mathcal O_{K_\alpha}$, and fix an additive character $\psi_K$ of $K$ with $\psi_K(x)=\psi(\bar x)$ for $x\in\mathcal O_K$.

For $\varphi_\pi$, Imai and Tsushima define characters $\xi$ and $\mu$ of $K_\alpha^\times$ by
\begin{gather*}
 \xi(\alpha)=1,\qquad \xi(x)=\chi(\bar x)\ \ (x\in\mu_{q-1}),\\
 \xi(y)=\psi_K\bigl(\Tr_{K_\alpha/K}(\alpha^{-1}(y-1))\bigr)\ \
 (y\in U^1_{K_\alpha}),\\
 \mu(\alpha)=\lambda_{K_\alpha/K}(\psi_K),\qquad \mu|_{K^\times}=\delta_{K_\alpha/K},\qquad \mu|_{U^1_{K_\alpha}}=1,
\end{gather*}
where $\delta_{K_\alpha/K}$ is the determinant of $\Ind_{K_\alpha/K}1$ and $\lambda_{K_\alpha/K}(\psi_K)$ is the Langlands constant, and they show that
\[
 \varphi_\pi\cong\Ind_{W_{K_\alpha}}^{W_K}(\mu^{-1}\xi)
\]
\cite[Proposition~5.4]{IT20}. Here the character $\mu$ is the rectifier, and Imai and Tsushima view characters of $K_\alpha^\times$ as characters of $W_{K_\alpha}$ through local reciprocity, normalized so that uniformizers correspond to geometric Frobenius elements \cite[\S1.1]{IT20}.

For $\rho$ we fix, as for \eqref{eq:Pi} and with the uniformizer $\varpi$ fixed here, an unramified extension $K_n\subset D$ of $K$ of degree $n$, with residue field $\F_{q^n}$, and $\Pi\in D$ with $\Pi^n=\varpi$ and $\Pi x\Pi^{-1}=\Frob(x)$ for $x\in K_n$. Thus $D=\bigoplus_{i<n}K_n\Pi^i$, and $\Pi x\Pi^{-1}=x^q$ for $x\in\mu_{q^n-1}$. We embed $K_\alpha$ in $D$ by $\alpha\mapsto\Pi$ and write $U_D^i=1+\Pi^i\mathcal O_D$. Imai and Tsushima define the character $\Theta$ of $K_\alpha^\times U_D^1$ by
\begin{gather*}
 \Theta(\Pi)=-1,\qquad \Theta(x)=\chi(\bar x)\ \ (x\in\mu_{q-1}),\\
 \Theta(h)=\psi_K\bigl(\Trd_{D/K}(\Pi^{-1}(h-1))\bigr)\ \ (h\in U_D^1)
\end{gather*}
(their value $(-1)^{n-1}c$ at $\Pi$ is $-1$ here), and they show that
\[
 \rho\cong\Ind_{K_\alpha^\times U_D^1}^{D^\times}\Theta
\]
\cite[Proposition~5.4]{IT20}. This induced representation is independent of the embedding of $K_\alpha$ in $D$ \cite[\S5.1]{IT20}.

Three facts about the residue field recur in the proofs below. Let $\varepsilon$ be the quadratic character of $\F_q^\times$. Then
\[
 \chi|_{\F_p^\times}=1,\qquad \varepsilon|_{\F_p^\times}=1,\qquad \varepsilon(b)=-1.
\]
The character $\chi$ is trivial on $\F_p^\times$ because $\kappa=b^{(q-1)/n}$ generates $\F_p^\times$ and $(q-1)/n=1+p+\dots+p^{n-1}\equiv n\equiv0\pmod n$, so that $\chi(\kappa)=\zeta_n^{(q-1)/n}=1$. The character $\varepsilon$ is trivial on $\F_p^\times$ because every element of $\F_p$ is a square in the extension $\F_q$ of even degree $n$, and this makes the rectifier value $\mu(\alpha)$ a sign (Lemma~\ref{lem:weil-side}). Finally, $\varepsilon(b)=-1$ because $b$ generates $\F_q^\times$.

On the Weil side, $\varphi_\pi$ is induced from the character $\xi'=\mu^{-1}\xi$ of the Weil group $W_{K_\alpha}$ of $K_\alpha$. The operators of $W_{K_\alpha}$ span over $E$ a diagonal copy of $E'$. Let $g\in W_K$ be an element with $a_K(g)=b$. The proof of Lemma~\ref{lem:weil-side} shows that $g$ restricts to a generator of $\Gal(K_\alpha/K)$, that the operator of $g$ normalizes the diagonal copy of $E'$ and acts on it by $\gamma_\kappa$, and that the copy and this operator generate $A_{\varphi_\pi}$, which is the cyclic algebra of Lemma~\ref{lem:weil-side}. Two steps need care here. With the reciprocity convention of \S\ref{subsec:local}, the element $g$ satisfies $g(\alpha)=\kappa\alpha$, and this action on $\alpha$ determines the automorphism $\gamma_\kappa$ by which conjugation by the operator of $g$ acts on the diagonal copy of $E'$. The scalar $\xi'(g^n)$ has a sign that comes from the quadratic character in the rectifier of Imai and Tsushima.

\begin{lemma}\label{lem:weil-side}
The character $\xi'=\mu^{-1}\xi$ takes values in $E'$, and
\begin{gather*}
 \xi'(\alpha)=\pm1,\qquad \xi'(x)=\chi\varepsilon(\bar x)\ \ (x\in\mu_{q-1}),\\
 \xi'(1+x\alpha)=\psi(n\bar x)\ \ (x\in\mathcal O_K),\qquad \xi'|_{U^2_{K_\alpha}}=1.
\end{gather*}
In particular $\xi'(b)=-\zeta_n$. Moreover $\Q(\varphi_\pi)=\Q(\tau)=E$, and
\[
 A_{\varphi_\pi}=A_\tau\cong(E'/E,\gamma_\kappa,-\zeta_n).
\]
\end{lemma}

\begin{proof}
Since $\mu$ is trivial on $U^1_{K_\alpha}$, the character $\xi'$ agrees with $\xi$ there, and $\Tr_{K_\alpha/K}(\alpha^{-1}\cdot x\alpha)=nx$ gives $\xi'(1+x\alpha)=\psi(n\bar x)$ for $x\in\mathcal O_K$. For $y\in U^2_{K_\alpha}$ the element $\alpha^{-1}(y-1)$ lies in $\alpha\mathcal O_{K_\alpha}$, whose traces lie in $\varpi\mathcal O_K$, so $\xi'(y)=1$, and $U^1_{K_\alpha}=(1+\alpha\mathcal O_K)U^2_{K_\alpha}$. On $\mu_{q-1}\subset K^\times$ the character $\mu$ is $\delta_{K_\alpha/K}$, which on $\mathcal O_K^\times$ is $x\mapsto\varepsilon(\bar x)^{n-1}=\varepsilon(\bar x)$ \cite[Lemma~5.3]{IT20}, since $n-1$ is odd. Hence $\xi'(x)=\chi\varepsilon(\bar x)$ there, and $\xi'(b)=\chi(b)\varepsilon(b)=-\zeta_n$.

Since $n$ is even, formula (5.22) in the proof of \cite[Proposition~5.15]{IT20} reads
\[
 \lambda_{K_\alpha/K}(\psi_K)=\Bigl(\frac{n/2}{q}\Bigr)\mathfrak m(\psi_K)^{n-1}, \qquad \mathfrak m(\psi_K)^2=\Bigl(\frac{-1}{q}\Bigr)
\]
\cite[(5.18)]{IT20}, where $(\frac{\cdot}{q})$ is the quadratic residue symbol of $\F_q$, that is, $\varepsilon$. Both $n/2$ and $-1$ lie in $\F_p^\times$, on which $\varepsilon$ is trivial, so $\mathfrak m(\psi_K)=\pm1$ and $\xi'(\alpha)=\lambda_{K_\alpha/K}(\psi_K)^{-1}=\pm1$. Thus $\xi'$ takes values in $E'$, and its values on $1+\alpha\mathcal O_K$ are all the $p$-th roots of unity, since $p\nmid n$.

To describe $\Ind\xi'$, we need an element of $W_K$ that restricts to a generator of $\Gal(K_\alpha/K)$. Let $g\in W_K$ be an element with $a_K(g)=b$. We claim that $g(\alpha)=\kappa\alpha$. Let $P(X)=\varpi X+X^q$ and let $t$ be a nonzero root of $P$. For $\zeta\in\mu_{q-1}$ the series $\zeta X$ commutes with $P$, so it is the endomorphism $[\zeta]_P$ of the Lubin--Tate formal group of $P$. The arithmetic Artin symbol of a unit $x\in\mathcal O_K^\times$ acts on the torsion points of this group by $[x^{-1}]_P$ \cite[I, 1.14, 1.15]{Mil20}, and by \eqref{eq:weil-galois} the element $g$ acts on them as the inverse of the arithmetic Artin symbol of $b$, that is, by $[b]_P$.

It remains to pass from the torsion point $t$ to a root of $X^n-\varpi$. Choose $\eta$ in the maximal unramified extension of $K$ with $\eta^n=-1$. The Artin symbol of a unit acts trivially on that extension \cite[I, 1.8]{Mil20}, and $g$ acts on it as the inverse of the Artin symbol of the unit $b$, so $g(\eta)=\eta$. Hence
\[
\begin{aligned}
 t^{q-1}&=-\varpi, & g(t)&=[b]_P(t)=bt,\\
 \beta&=t^{(q-1)/n}, & \beta^n&=-\varpi,\quad g(\beta)=b^{(q-1)/n}\beta=\kappa\beta,\\
 \alpha'&=\beta/\eta, & \alpha'^n&=\varpi,\quad g(\alpha')=\kappa\alpha'.
\end{aligned}
\]
Every root of $X^n-\varpi$ has the form $\zeta\alpha'$ with $\zeta\in\mu_n\subset K$, so $g(\alpha)=\kappa\alpha$, and $g$ restricts to a generator of $\Gal(K_\alpha/K)$ because $\kappa$ has order $n$.

Local reciprocity is compatible with conjugation, so for $h\in W_{K_\alpha}$ corresponding to $y\in K_\alpha^\times$ the element $ghg^{-1}$ corresponds to $g(y)$. We check that $\xi'(g(y))=\gamma_\kappa(\xi'(y))$ on generators of $K_\alpha^\times$. On $\mu_{q-1}$ the automorphism $g$ acts trivially, and $\xi'$ takes values in $E$ there. For $x\in\mathcal O_K$ we have $g(1+x\alpha)=1+\kappa x\alpha$, and $\psi(n\kappa\bar x)=\gamma_\kappa(\psi(n\bar x))$ because the residue class of $\kappa$ lies in $\F_p$. Finally, $\xi'(g(\alpha))=\xi'(\kappa)\xi'(\alpha)=\xi'(\alpha)$, because $\xi'(\kappa)=\chi\varepsilon(\kappa)=1$ (both characters are trivial on $\F_p^\times$), and $\xi'(\alpha)=\pm1$ is fixed by $\gamma_\kappa$. Hence
\begin{equation}\label{eq:conjugation-rule}
 \xi'(ghg^{-1})=\gamma_\kappa\bigl(\xi'(h)\bigr)\qquad(h\in W_{K_\alpha}).
\end{equation}

The elements $1,g,\dots,g^{n-1}$ represent $W_K/W_{K_\alpha}$. In the basis $e_i=g^i\otimes1$ of $\Ind\xi'$, an element $h\in W_{K_\alpha}$ acts by $\iota(\xi'(h))$, where
\[
 \iota(z)=\diag\bigl(z,\gamma_\kappa^{-1}z,\dots,\gamma_\kappa^{-(n-1)}z\bigr)\qquad (z\in E')
\]
by \eqref{eq:conjugation-rule}, and $u=\varphi_\pi(g)$ maps $e_i$ to $e_{i+1}$ for $i<n-1$ and $e_{n-1}$ to $\xi'(g^n)e_0$. Since $g$ restricts to a generator of $\Gal(K_\alpha/K)$, Lemma~\ref{lem:weil-power} shows that $a_{K_\alpha}(g^n)=a_K(g)=b$, and therefore
\[
 u\,\iota(z)\,u^{-1}=\iota(\gamma_\kappa z),\qquad u^n=\xi'(g^n)=\xi'(b)=\chi(b)\,\varepsilon(b)=-\zeta_n.
\]

Let $A=\bigoplus_{i<n}\iota(E')u^i$. As in the proof of Proposition~\ref{prop:image-algebras}, the summand $\iota(E')u^i$ maps each line $\C e_j$ into $\C e_{j+i}$ (indices modulo $n$), so the sum is direct and $\dim_E A=n^2$, and the relations above show that $A\cong(E'/E,\gamma_\kappa,-\zeta_n)$. Every operator $\varphi_\pi(g^ih)=u^i\iota(\xi'(h))$ lies in $A$. Conversely, $\iota(E')$ is the $E$-span of the matrices $\iota(\xi'(h))$, because the values of $\xi'$ span $E'$ over $E$, and $u$ lies in the image of $\varphi_\pi$. Hence $A_{\varphi_\pi}=A$.

If $\sigma\in\mathrm{Aut}(\C)$ fixes $E$, it acts on $E'$ as some $\gamma_\kappa^k$, and $\sigma\circ\varphi_\pi\cong\Ind(\gamma_\kappa^k\circ\xi')\cong\varphi_\pi$ by \eqref{eq:conjugation-rule}. Conversely, the element of $W_{K_\alpha}$ corresponding to $b\in K^\times$ acts on $\varphi_\pi$ by the scalar $-\zeta_n$, which every $\sigma$ with $\sigma\circ\varphi_\pi\cong\varphi_\pi$ must fix, so $\Q(\varphi_\pi)=\Q(\zeta_n)=E$. Finally, since $q=p^n$ with $n$ even, $|w|^{-(n-1)/2}=p^{-\nu(w)n(n-1)/2}$ is rational for $w\in W_K$, so $\tau(w)$ is a nonzero rational multiple of $\varphi_\pi(w)$. Hence $A_\tau=A_{\varphi_\pi}$ and $\Q(\tau)=\Q(\varphi_\pi)=E$.
\end{proof}

On the side of $D^\times$, the image algebra $A_\rho$ has degree $(q^n-1)/(q-1)>n$, so the class of $\rho$ needs one further idea. We find an idempotent $e_0\in A_\rho$ of rank $n$ whose corner algebra $e_0A_\rho e_0$ is a cyclic algebra of degree $n$. Mackey's formula splits the restriction of $\rho$ to $U_D^1$ into lines $V_c$, indexed by the kernel $N^1$ of a norm, on which the additive group $U_D^1/U_D^2$ acts by distinct characters. The projection $e_0$ onto the $n$ lines with $c\in\F_p^\times$ lies in $A_\rho$, because the sum over $\F_p^\times$ turns the roots of unity in the Fourier expression for $e_0$ into integers. This is the hard step. The corner then contains a copy $\iota'(E')$ of $E'$ and a unit $u'$ that cycles the $n$ lines and represents the cyclic generator. Finally, the corner and $A_\rho$ are matrix algebras over the same division algebra, so the corner has the Brauer class of $A_\rho$.

\begin{lemma}\label{lem:division-side}
Let $V$ be the space of $\rho$, let $N^1\subset \F_{q^n}^\times$ be the kernel of the norm $\F_{q^n}^\times\to \F_q^\times$, and put $\Psi(x)=\zeta_p^{\Tr_{\F_{q^n}/\F_p}(x)}$ for $x\in \F_{q^n}$.
\begin{enumerate}
\item[(i)] The space $V$ is the direct sum of lines $V_c$, $c\in N^1$, on which $U_D^2$ acts trivially and $1+x\Pi$ acts by $\Psi(c\bar x)$ for $x\in\mathcal O_{K_n}$.
\item[(ii)] The central character of $\rho$ has order $n$, and $\Q(\rho)=E$.
\item[(iii)] Let $e_0$ be the projection of $V$ onto $\bigoplus_{c\in\F_p^\times}V_c$. Then $e_0\in A_\rho$, and $e_0A_\rho e_0\cong(E'/E,\gamma_\kappa,\zeta_n^{-1})$. In particular $o_E(\rho)=[(E'/E,\gamma_\kappa,\zeta_n^{-1})]$.
\end{enumerate}
\end{lemma}

\begin{proof}
(i) For $x\in\mathcal O_{K_n}$ the element $\Pi^{-1}x\Pi$ lies in $K_n$ and has the same trace over $K$ as $x$, so $\Theta(1+x\Pi)=\psi_K(\Tr_{K_n/K}x)=\Psi(\bar x)$. If $h\in U_D^2$, then $\Pi^{-1}(h-1)$ lies in $\Pi\mathcal O_D=\bigoplus_{1\le i\le n}\mathcal O_{K_n}\Pi^i$. The reduced trace of $y\Pi^i$ vanishes for $y\in K_n$ and $0<i<n$, and $\Trd(y\varpi)\in\varpi\mathcal O_K$ for $y\in\mathcal O_{K_n}$. So $\Theta$ is trivial on $U_D^2$ and induces $\Psi$ on $U_D^1/U_D^2\cong \F_{q^n}$, $1+x\Pi\mapsto\bar x$.

The group $D^\times$ is the disjoint union of the cosets $dK_\alpha^\times U_D^1$, for $d$ in a set of representatives of $\mu_{q^n-1}/\mu_{q-1}$, and $U_D^1$ is normal, so by Mackey's formula $\rho|_{U_D^1}$ is the sum of the characters $h\mapsto\Theta(d^{-1}hd)$. Since $\Pi d=d^q\Pi$, these characters are given by
\[
 d^{-1}(1+x\Pi)d=1+d^{q-1}x\Pi,\qquad \Theta\bigl(d^{-1}(1+x\Pi)d\bigr)=\Psi(d^{q-1}\bar x).
\]
The map $d\mapsto d^{q-1}$ induces a bijection
\[
 \mu_{q^n-1}/\mu_{q-1}\xrightarrow{\ \sim\ }N^1
\]
(its image lies in $N^1$ because $d^{(q-1)(1+q+\dots+q^{n-1})}=1$, and both groups have order $(q^n-1)/(q-1)$), and the characters $x\mapsto\Psi(cx)$ are distinct for distinct $c$ because the trace form of $\F_{q^n}/\F_p$ is nondegenerate. So the lines of $\rho|_{U_D^1}$ are indexed by $c\in N^1$, and $1+x\Pi$ acts on the line $V_c$ by $\Psi(c\bar x)$. This proves (i).

(ii) The centre $K^\times$ of $D^\times$ is generated by $\varpi=\Pi^n$, $\mu_{q-1}$ and $U_K^1\subset U_D^2$. Since $n$ is even, $\rho(\varpi)=\Theta(\Pi)^n=1$, and $\rho$ acts on $\mu_{q-1}$ by $\chi$, so the central character of $\rho$ has the order $n$ of $\chi$, and $\Q(\rho)$ contains $\Q(\zeta_n)=E$. Conversely, let $\sigma\in\mathrm{Aut}(\C)$ fix $E$. It acts on $E'$, which contains the values of $\Theta$, as $\gamma_c$ for some $c\in\F_p^\times$, where $\gamma_c(\zeta_p)=\zeta_p^{c}$, so $\sigma\circ\Theta$ agrees with $\Theta$ on $\Pi$ and on $\mu_{q-1}$ and sends $1+x\Pi$ to $\Psi(c\bar x)$. Since $c^n=1$, there is $d\in\mu_{q^n-1}$ with $d^{q-1}=c$, and $d$ normalizes $K_\alpha^\times U_D^1$ because $d\Pi d^{-1}=c^{-1}\Pi$. The conjugate character $\Theta^d(h)=\Theta(d^{-1}hd)$ also sends $1+x\Pi$ to $\Psi(c\bar x)$, by the computation in (i), agrees with $\Theta$ on $\mu_{q-1}$, and takes the value $\Theta^d(\Pi)=\Theta(c\Pi)=\chi(c)\Theta(\Pi)=\Theta(\Pi)$ at $\Pi$, because $\chi$ is trivial on $\F_p^\times$. Hence $\sigma\circ\Theta=\Theta^d$, $\sigma\circ\rho\cong\rho$, and $\Q(\rho)=E$.

(iii) Let $\mathrm{pr}_c$ be the projection of $V$ onto $V_c$. Since $U_D^2$ acts trivially,
\[
 \mathrm{pr}_c=q^{-n}\sum_{x\in \F_{q^n}}\Psi(cx)^{-1}\rho(1+[x]\Pi),
\]
where $[x]\in\mathcal O_{K_n}$ is the Teichm\"uller representative of $x$. The norm of $c\in\F_p^\times$ is $c^n=1$, so $\F_p^\times\subset N^1$, and
\[
 e_0=\sum_{c\in\F_p^\times}\mathrm{pr}_c=q^{-n}\sum_{x\in \F_{q^n}}\Bigl(\sum_{c\in\F_p^\times} \zeta_p^{-c\Tr_{\F_{q^n}/\F_p}(x)}\Bigr)\rho(1+[x]\Pi).
\]
The inner sums are rational integers ($p-1$ or $-1$), so $e_0\in A_\rho$, and $e_0$ has rank $n$.

Next we construct a copy of $E'$ in the corner $e_0A_\rho e_0$. For $z\in E'$ put
\[
 \iota'(z)=\sum_{c\in\F_p^\times}\gamma_c(z)\,\mathrm{pr}_c,\qquad \text{so that}\qquad e_0\,\rho(1+[x]\Pi)=\iota'\bigl(\Psi(x)\bigr)\quad(x\in\F_{q^n}),
\]
since, by (i), the operator $\rho(1+[x]\Pi)$ acts on every line $V_c$, $c\in N^1$, by the scalar $\Psi(cx)$, and $\Psi(cx)=\gamma_c(\Psi(x))$ for $c\in\F_p^\times$. The operators $e_0\rho(1+[x]\Pi)$ lie in $e_0A_\rho e_0$, because $e_0\in A_\rho$ and $\rho(1+[x]\Pi)$ preserves every line $V_c$, $c\in N^1$, and so commutes with $e_0$. The map $\iota'$ is $E$-linear, and the values of $\Psi$ span $E'$ over $E$, so $\iota'$ maps $E'$ into $e_0A_\rho e_0$, and it is an injective homomorphism of $E$-algebras $E'\to e_0A_\rho e_0$.

Choose $d_0\in\mu_{q^n-1}$ with $d_0^{q-1}=\kappa^{-1}$, which exists because $\kappa^{-1}\in N^1$. As in (i), $d_0^{-1}(1+x\Pi)d_0=1+\kappa^{-1}x\Pi$, so $\rho(d_0)$ maps $V_c$ onto $V_{c\kappa^{-1}}$. In particular $\rho(d_0)$ commutes with $e_0$, and $u'=\rho(d_0)e_0$ is a unit of $e_0A_\rho e_0$ that maps each line $V_c$, $c\in\F_p^\times$, onto $V_{c\kappa^{-1}}$. Since $\iota'(z)$ acts on $V_c$ by $\gamma_c(z)$ and $\gamma_{c\kappa^{-1}}\gamma_\kappa=\gamma_c$, one checks on each line $V_{c\kappa^{-1}}$ that $u'\,\iota'(z)\,u'^{-1}$ and $\iota'(\gamma_\kappa z)$ both act by $\gamma_c(z)$.

Next, $(d_0^n)^{q-1}=\kappa^{-n}=1$, so $d_0^n$ lies in $\mu_{q-1}$, which is central. Write $d_0^n=b^t$. Then $\kappa^t=(d_0^n)^{(q-1)/n}=d_0^{q-1}=\kappa^{-1}$, so $t\equiv-1\pmod n$ and $\rho(d_0^n)=\chi(b)^t=\zeta_n^{-1}$. So $u'$ and $\iota'$ satisfy the relations of a cyclic algebra,
\[
 u'\,\iota'(z)\,u'^{-1}=\iota'(\gamma_\kappa z)\quad(z\in E'),\qquad u'^n=\rho(d_0^n)\,e_0=\zeta_n^{-1}e_0.
\]

Let $B=\bigoplus_{i<n}\iota'(E')u'^i\subset e_0A_\rho e_0$. Since $u'^i$ maps $V_c$ onto $V_{c\kappa^{-i}}$, the summands occupy disjoint blocks of $\End(e_0V)$, so the sum is direct, $\dim_EB=n^2$, and $B\cong(E'/E,\gamma_\kappa,\zeta_n^{-1})$. Since $\Q(\rho)=E$ by (ii), Lemma~\ref{lem:trace-span}, applied with $F=E$, shows that the natural map $A_\rho\otimes_E\C\to\End_\C(V)$ is an isomorphism. It carries $e_0A_\rho e_0\otimes_E\C$ onto $\End_\C(e_0V)$, which has dimension $n^2$ because $e_0$ has rank $n$. So $\dim_Ee_0A_\rho e_0=n^2$, and $B=e_0A_\rho e_0$.

Finally, write $A_\rho\cong\End_\Delta(W)$ with $\Delta$ a central division algebra over $E$. The corner at the nonzero idempotent $e_0$ is $\End_\Delta(e_0W)$, a matrix algebra over the same $\Delta$, so it has the class $o_E(\rho)$ of $A_\rho$.
\end{proof}

\begin{proof}[Proof of Theorem~\ref{thm:family}]
(i) By Lemma~\ref{lem:division-side}(ii) the central character of $\rho$ has order $n\ge6$. The Jacquet--Langlands correspondence preserves central characters, so the same holds for $\pi$, and since a self-dual representation has central character of order at most two, neither $\pi$ nor $\rho$ is self-dual. By Lemma~\ref{lem:common-rationality-field}, $\Q(\pi)=\Q(\rho)=\Q(\tau)$, and this field is $E$ by Lemmas~\ref{lem:weil-side} and~\ref{lem:division-side}.

(ii) This is Lemma~\ref{lem:division-side}(iii) for $\rho$ and Lemma~\ref{lem:weil-side} for $\tau$.

(iii) Since $p\equiv1\pmod n$, the prime $p$ splits completely in $E=\Q(\zeta_n)$, and the places of $E$ above $p$ correspond to the images of $\zeta_n$ in $\F_p$, which are the generators $\kappa^j$, $j\in(\Z/n\Z)^\times$, of $\F_p^\times$. At $v=v_j$ we have $E_v=\Q_p$, and $E'\otimes_EE_v=\Q_p(\zeta_p)$ is totally ramified of degree $n$ over $\Q_p$, with $\gamma_\kappa$ acting by $\zeta_p\mapsto\zeta_p^\kappa$.

By Lemma~\ref{lem:cyclic-split-invariant}(iii), it suffices to know how the arithmetic Artin symbols of the two parameters act on $\Q_p(\zeta_p)$. For a unit $y\in\Z_p^\times$, the arithmetic Artin symbol of $y$ acts on $\zeta_p$ by $\zeta_p\mapsto\zeta_p^{\bar y^{-1}}$ \cite[I, Example~3.13]{Mil20}, that is, as $\gamma_\kappa^t$ with $\kappa^t=\bar y^{-1}$. At $v_j$ we have $\zeta_n\equiv\kappa^j$ and $-1\equiv\kappa^{n/2}$, so the parameters $y=\zeta_n^{-1}$ of $o_E(\rho)$ and $y=-\zeta_n$ of $o_E(\tau)$ have $\bar y^{-1}=\kappa^j$ and $\bar y^{-1}=-\kappa^{-j}=\kappa^{n/2-j}$, and therefore, in $\Q/\Z$,
\[
 \inv_{v_j}o_E(\rho)=\frac jn,\qquad \inv_{v_j}o_E(\tau)=\frac{n/2-j}n=-\frac jn-\frac12.
\]

At a finite place $v$ of $E$ not above $p$, the extension $E'/E$ is unramified and $\zeta_n^{-1}$ and $-\zeta_n$ are units. The description of Lemma~\ref{lem:cyclic-split-invariant}(iii), applied to the restriction to the decomposition group at $v$ of the character of $\Gal(E'/E)$ that sends $\gamma_\kappa$ to $1/n$, shows that both invariants vanish, because the Artin symbol of a unit lies in the inertia group, which acts trivially on $E'\otimes_EE_v$. The field $E$ is totally imaginary, so the invariants at the infinite places vanish as well.

(iv) For $j$ prime to $n$, the invariant $j/n$ has order $n$, and by (iii) the two invariants at $v_j$ add up to $\frac12$ in $\Q/\Z$. Since $p\equiv3\pmod4$, we have $n=p-1\equiv2\pmod4$, so Lemma~\ref{lem:half-shift} shows that $-j/n-\frac12$ has order $n/2$. These orders are the local Schur indices at $v_j$, every other place has local index $1$, and the global index is the least common multiple of the local indices (\S\ref{subsec:index-relations}), so the global indices are $n$ and $n/2$.

(v) The local invariants of $o_E(\rho)+o_E(\tau)$ are $\frac12$ at each $v_j$ and $0$ elsewhere. Those of $\res_{E/\Q}[\Bpinf]$ are $[E_v:\Q_p]/2=\frac12$ at the places above $p$ \cite[III, \S2, (29)]{Mil20} and $0$ elsewhere, since $E$ has no real place. So the two classes are equal \cite[VII, Theorem~7.1]{Mil20}. Since $n-1$ is odd and $2[\Bpinf]=0$, this is \eqref{eq:intro-sum}.
\end{proof}

\begin{remark}\label{rem:sign}
The sign in $\xi'(b)=-\zeta_n$ produces the quaternion class $\res_{E/\Q}[\Bpinf]$ of Theorem~\ref{thm:family}(v), since
\[
 o_E(\rho)+o_E(\tau)=[(E'/E,\gamma_\kappa,\zeta_n^{-1})]+[(E'/E,\gamma_\kappa,-\zeta_n)]=[(E'/E,\gamma_\kappa,-1)].
\]
Both classes use the generator $\gamma_\kappa$, so additivity in the parameter (Proposition~\ref{prop:cyclic}(i)) leaves the product $\zeta_n^{-1}\cdot(-\zeta_n)=-1$, the sign that the rectifier contributes, and the cyclic algebra $(E'/E,\gamma_\kappa,-1)$ of degree $n$ represents the quaternion class. The central character also forces this sign. Let $g$ and $u=\varphi_\pi(g)$ be as in the proof of Lemma~\ref{lem:weil-side}. Since $u$ maps $e_i$ to $e_{i+1}$ for $i<n-1$ and $e_{n-1}$ to $\xi'(g^n)e_0$, it permutes the lines $\C e_i$ cyclically, and its determinant is $(-1)^{n-1}\xi'(g^n)$. The local Langlands correspondence matches $\det\varphi_\pi$ with the central character of $\pi$, which is $\chi$ on $\mu_{q-1}$, and $a_K(g)=b$, so
\[
 (-1)^{n-1}\xi'(g^n)=\det u=\chi(b)=\zeta_n,
\]
and hence $\xi'(g^n)=-\zeta_n$, as $n$ is even.

So, given the induced forms of $\rho$ and of $\varphi_\pi$, the central character fixes the parameters of both cyclic presentations for $E'/E$, namely $\zeta_n^{-1}=\rho(d_0^n)$ for $\rho$, where $d_0^n\in\mu_{q-1}\subset K^\times$ is central (proof of Lemma~\ref{lem:division-side}(iii)), and $-\zeta_n=\xi'(g^n)$ for $\varphi_\pi$. The proof of Theorem~\ref{thm:family}(iii) and~(iv) computes the indices $n$ and $n/2$ from these two parameters, so the indices do not depend on the Langlands constant $\lambda_{K_\alpha/K}(\psi_K)$ or on the choice of reciprocity map used to describe $\varphi_\pi$. Replacing one of the two algebras by its opposite negates the invariants of that algebra, and reversing the sign convention for local invariants negates the invariants of both. Either change preserves the orders of the invariants, and hence the indices. Under the opposite normalization of Proposition~\ref{prop:cyclic}(ii), the algebra $D$ of \eqref{eq:Pi} would be called the division algebra of invariant $-1/n$, and $\rho$, $\tau$, $A_\rho$ and $A_\tau$ would stay the same.
\end{remark}

\begin{remark}\label{rem:family-hypotheses}
The congruence $p\equiv3\pmod4$ is used only in part~(iv) of the proof of Theorem~\ref{thm:family}. For $j$ prime to $n$, the order of $-j/n-\frac12=-(j+n/2)/n$ is $n/\gcd(n,j+n/2)$, and this gcd is $2$ if and only if $n/2$ is odd. When $n/2$ is even, the numbers $j/n$ and $-j/n-\frac12$ both have order $n$, and for $n=4$ they are equal (both are $1/4$ for $j=1$, and $3/4$ for $j=3$).

The condition $p\ge7$ makes the family non-self-dual and $E$ totally imaginary. For $p=3$, where $n=2$ and $\chi$ is the quadratic character of $\F_9^\times$, the argument of Lemma~\ref{lem:division-side}(ii) shows that $\Q(\rho)=\Q(\zeta_2)=\Q$, and since $\rho$ factors through a finite quotient of $D^\times$, it follows that $\rho^\vee\cong\bar\rho\cong\rho$, so $\rho$ is self-dual.
\end{remark}

The representations $\rho$ of this family are not self-dual, and their field of rationality $E=\Q(\zeta_n)$ is cyclotomic. This is the setting of the conjecture of Prasad and Ramakrishnan \cite[Conjecture~6.1(1)]{PR12}, to which Section~\ref{sec:counterexamples} turns.

\section{The half twist and the conjecture of Prasad and Ramakrishnan}\label{sec:counterexamples}

The smallest counterexamples to the conjecture of Prasad and Ramakrishnan have rank two and depth zero. Example~\ref{ex:intro-five} already exhibits a representation $\rho$ of $D^\times$ that is not self-dual, with $o_E(\rho)=0$ and $o_E(\tau)\neq0$, and in \S\ref{subsec:order-eight} we construct such pairs for every $p\equiv5\pmod8$, with the cyclotomic centre $E=\Q(i)$, as the case $n=2$ of Section~\ref{sec:depth-zero}. The same construction for $p=3$ gives such a pair with centre $\Q(\sqrt{-2})$, which is abelian over $\Q$ and is not a cyclotomic field, so this pair contradicts the conjecture in the second reading of Remark~\ref{rem:cyclotomic-centre}, which requires the centre only to be abelian over $\Q$. In these pairs $\rho$ and $\varphi_\pi$ both have models over $E$, so the half twist produces the difference between the two Schur algebras. Corollary~\ref{cor:pr-counterexamples} collects them together with the family of \S\ref{subsec:family}, whose indices are unbounded.

In \S\ref{subsec:dyadic} we compare these examples with a proposition of Suzuki on the local invariants above $p$ \cite[Proposition~4.0.13]{Su23}. Read in our notation under the hypothesis of Remark~\ref{rem:suzuki-step}, the proposition makes two assertions, and the order-eight examples contradict one of them and an example over $\Q_2$ the other. All these examples satisfy \eqref{eq:intro-sum}.

\subsection{The order-eight examples}\label{subsec:order-eight}

We first recall the normalization under which Prasad and Ramakrishnan compare the two Schur algebras. Prasad and Ramakrishnan multiply the local Langlands correspondence by the character $x\mapsto|x|^{(n-1)/2}$ of $K^\times$, which makes it commute with the automorphisms of $\C$ \cite[\S6]{PR12}, and we read this, with the reciprocity conventions of \S\ref{subsec:local}, as a twist of the representation of $D^\times$. Put $\nu_D=|\Nrd_{D/K}|$, and for an irreducible $n$-dimensional parameter $r$ let $\pi_r$ be the representation of $\GL_n(K)$ with Langlands parameter $r$. The representation of $D^\times$ attached to $r$ is then $\JL(\pi_r)\otimes\nu_D^{(n-1)/2}$. Since the Langlands and Jacquet--Langlands correspondences are compatible with twisting by powers of $|\det|$, $|\cdot|$ and $\nu_D$ (\cite[\S2.2.1]{Dat07} and Lemma~\ref{lem:jl-twist}), the parameter attached to $\rho=\JL(\pi)$ is $\varphi_\pi\otimes|\cdot|^{-(n-1)/2}=\tau$.

The twist of the opposite sign, $\varphi_\pi\otimes|\cdot|^{(n-1)/2}$, differs from $\tau$ by the rational-valued character $|\cdot|^{n-1}$, so it has the same field of rationality and image algebra as $\tau$, and the comparison below holds for both signs of the twist. The same normalization also attaches $\rho\otimes\nu_D^{(n-1)/2}$ to $\varphi_\pi$, which places the twist on the side of $D^\times$. In the order-eight examples below, the Schur algebras of this pair differ as well, since $\varphi_\pi$ has a model over $E$ and the twist replaces the cyclic parameter $1$ of $A_\rho$ by $p^{-1}$, so that $o_E(\rho\otimes\nu_D^{1/2})$ has invariant $\frac12$ at the places of $E$ above $p$.

Recall that for a finite-dimensional irreducible representation $r$ whose field of rationality $F$ is a number field, the Schur algebra of $r$ is the central division algebra over $F$ in the class $o_F(r)$, so that $A_r$ is a matrix algebra over it \cite[Lemma~3.0.7]{Su23}. Prasad and Ramakrishnan conjectured the following \cite[Conjecture~6.1]{PR12}. Let $D$ be a central division algebra of index $d$ over $K$, let $\pi'$ be an irreducible discrete series representation of $\GL_m(D)$ with $md=n$, and let $\sigma$ be its $n$-dimensional parameter, normalized in the same way. The conjecture concerns the Schur algebras of $\pi'$ and $\sigma$, and it describes their common centre as a number field ``(which is a cyclotomic field)''. Part~(1) of the conjecture asserts that the two division algebras are isomorphic when $\pi'$ is not self-dual, and also when $\pi'$ is self-dual and the signs of $\pi'$ and $\sigma$ agree.

We use only the case in which $\pi'$ is not self-dual. The representations of $D^\times$ used in our counterexamples have finite-order central character. They are unitary and square-integrable modulo the centre, since $D^\times/K^\times$ is compact, and hence are discrete series. For $m=1$, the non-self-dual case therefore asserts $o_E(\rho)=o_E(\tau)$ for these representations.

We now construct the smallest examples, which have $n=2$. Let $p=3$ or $p\equiv5\pmod8$, and let $K=\Q_p$, so that $q=p$. We choose $s$ with $a_{\Q_p}(s)=p$, so that the uniformizer $\varpi$ of Section~\ref{sec:depth-zero} is $p$. Let $D$ be the quaternion division algebra over $\Q_p$, let $K_2\subset D$ be an unramified quadratic extension of $\Q_p$, and let $\Pi\in D$ satisfy $\Pi^2=p$ and $\Pi x\Pi^{-1}=\bar x$ for $x\in K_2$, where $x\mapsto\bar x$ is the nontrivial automorphism of $K_2$.

Let $\theta$ be a character of $\F_{p^2}^\times$ of order $8$, viewed as a character of $K_2^\times U_D^1$ that is trivial on $U_D^1$ and on $p$. Put
\[
 \rho_\theta=\Ind_{K_2^\times U_D^1}^{D^\times}\theta,
\]
let $\pi_\theta$ be the representation of $\GL_2(\Q_p)$ with $\JL(\pi_\theta)=\rho_\theta$, let $\varphi_\theta$ be the Langlands parameter of $\pi_\theta$, and let $\tau_\theta=\varphi_\theta\otimes|\cdot|^{-1/2}$ be its normalized parameter. Let $L=\Q(\zeta_8)$, let $\gamma_p\in\Gal(L/\Q)$ be the automorphism with $\gamma_p(\zeta_8)=\zeta_8^p$, and put $E=L^{\langle\gamma_p\rangle}$. In the notation of Section~\ref{sec:depth-zero}, $\rho_\theta$ and $\tau_\theta$ are the depth-zero representations $\rho_1$ and $\tau_1$ for $n=2$, $K=\Q_p$, $m=8$ and $z=1$.

\begin{theorem}\label{thm:order-eight}
Taking $p$, $\theta$, $L$ and $E$ as above, we have the following:
\begin{enumerate}
\item[(i)] The representation $\rho_\theta$ is irreducible and not self-dual, its central character has finite order, and $\Q(\rho_\theta)=\Q(\pi_\theta)=\Q(\tau_\theta)=E$. The field $E$ is $\Q(\sqrt{-2})$ if $p=3$ and $\Q(i)$ if $p\equiv5\pmod8$.
\item[(ii)] The representations $\rho_\theta$ and $\varphi_\theta$ have models over $E$.
\item[(iii)] The algebra $A_{\tau_\theta}$ is isomorphic to $(L/E,\gamma_p,-p)$, a quaternion division algebra over $E$ with invariant $\frac12$ at each of the two places of $E$ above $p$ and invariant $0$ at every other place.
\end{enumerate}
\end{theorem}

The case $n=2$ of Proposition~\ref{prop:image-algebras} and Theorem~\ref{thm:central-value-one} gives the model of $\rho_\theta$ in (ii), part (iii) and the field of rationality of $\rho_\theta$ in (i). The remaining point is a model of the untwisted parameter $\varphi_\theta$ over $E$. We realize $\varphi_\theta$ on the $E$-vector space $L$, with $s$ acting by $\gamma_p$ followed by multiplication by an element of $L$ whose norm to $E$ is $-1$, so that the operator of $s$ has square $-1$. The half twist multiplies this square by $p$, and $-p$ is not a norm at the places of $E$ above $p$, which is why $\tau_\theta$ has no model over $E$.

\begin{proof}
(i) Fix a generator $t$ of $\F_{p^2}^\times$. We may take $\zeta_8=\theta(t)$, and then $\rho_\theta$ is the representation $\rho_1$ of Section~\ref{sec:depth-zero} for $n=2$, $K=\Q_p$ and $m=8$, with $\gamma=\gamma_p$, as noted before the theorem. The character $\theta$ is primitive, since $\theta^p\neq\theta$ because $8\nmid p-1$. So $\rho_\theta$ is irreducible by Lemma~\ref{lem:shape}, $\Q(\rho_\theta)=E$ by Proposition~\ref{prop:image-algebras}, and $\Q(\pi_\theta)=\Q(\tau_\theta)=E$ by Lemma~\ref{lem:common-rationality-field}.

Concretely, $\rho_\theta$ factors through $D^\times/p^{\Z}U_D^1\cong\F_{p^2}^\times\rtimes\langle j\rangle$, where $j$ is the image of $\Pi$, $j^2=1$ and $jxj^{-1}=x^p$ for $x\in\F_{p^2}^\times$, and in a basis of eigenvectors of $t$
\[
 \rho_\theta(t)=\begin{pmatrix}\zeta_8&0\\0&\zeta_8^{p}\end{pmatrix},\qquad \rho_\theta(j)=\begin{pmatrix}0&1\\1&0\end{pmatrix}.
\]
The image of $\rho_\theta$ is finite, so its central character has finite order. The exponents $\{1,p\}$ of $\rho_\theta$ and $\{-1,-p\}$ of its contragredient differ modulo $8$, so $\rho_\theta$ is not self-dual. If $p=3$, the element $c_8=\zeta_8+\zeta_8^3$ is fixed by $\gamma_3$ and $c_8^2=-2$, so $E=\Q(c_8)=\Q(\sqrt{-2})$. If $p\equiv5\pmod8$, then $\gamma_p(\zeta_8)=-\zeta_8$, so $\gamma_p$ fixes $i=\zeta_8^2$ and $E=\Q(i)$.

(ii) and (iii) Here $a=-p$ and $f=1$, so the case $n=2$ of Proposition~\ref{prop:image-algebras} and Theorem~\ref{thm:central-value-one} shows that $\rho_\theta$ has a model over $E$, that $A_{\tau_\theta}\cong(L/E,\gamma_p^{-1},-p)=(L/E,\gamma_p,-p)$, since $\gamma_p$ has order two, and that $o_E(\tau_\theta)$ has invariant $\frac12$ at each place above $p$ and $0$ at the other finite places. The field $E$ has no real place, because $-1\notin\langle p\rangle=\{1,p\}\subset(\Z/8\Z)^\times$, and $p$ splits in $E$ ($-2$ is a square modulo $3$, and $p\equiv1\pmod4$ when $p\equiv5\pmod8$), so $A_{\tau_\theta}$ is a division algebra with nonzero invariants at the two places above $p$.

For $\varphi_\theta$ we use Proposition~\ref{prop:depth-zero-parameter} with $n=2$ and $z=1$. It shows that $\varphi_\theta=\Ind_{W_{K_2}}^{W_K}\mu_1$, where $\mu_1$ agrees with $\theta$ on inertia and $\mu_1(s^2)=-1$. So $\varphi_\theta$ is trivial on wild inertia, its restriction to inertia has the characters $\theta$ and $\theta^p$, $s$ exchanges the two eigenlines, and $\varphi_\theta(s)^2=-1$.

On the $E$-vector space $L$, let the inertia group act through $\theta$, by multiplication, and let $s$ act by the $E$-linear operator
\[
 x\mapsto\zeta_8\,\gamma_3(x)\quad\text{if }p=3,\qquad x\mapsto i\,\gamma_p(x)\quad\text{if }p\equiv5\pmod8,
\]
whose square is multiplication by
\[
 \zeta_8\,\gamma_3(\zeta_8)=\zeta_8^4=-1,\qquad\text{respectively}\qquad i\,\gamma_p(i)=i^2=-1.
\]
Conjugation by $s$ raises the inertia character to the power $p$, which is the tame relation because $p^{-1}\equiv p\pmod 8$. Over $\C$ the eigenlines of the inertia group in $L\otimes_E\C$ have the characters $\theta$ and $\theta^p$ and are exchanged by $s$, so this is a model of $\varphi_\theta$ over $E$. Finally, an automorphism $\zeta_8\mapsto\zeta_8^c$ of $L$ carries $\varphi_\theta$ to the representation whose restriction to inertia has the exponents $\{c,cp\}$, which is isomorphic to $\varphi_\theta$ if and only if $c\in\{1,p\}$. Hence $\Q(\varphi_\theta)=E$.
\end{proof}

\begin{corollary}\label{cor:pr-counterexamples}
We have the following:
\begin{enumerate}
\item[(i)] Let $p\equiv5\pmod8$. The pair $(\rho_\theta,\tau_\theta)$ of Theorem~\ref{thm:order-eight} satisfies the hypotheses of the non-self-dual case of \cite[Conjecture~6.1(1)]{PR12}, with centre the cyclotomic field $\Q(i)$, and its Schur algebras are not isomorphic. For $p=3$ the same holds, with centre $\Q(\sqrt{-2})$, if the centre is only required to be an abelian extension of $\Q$.
\item[(ii)] For every prime $p\ge7$ with $p\equiv3\pmod4$, the pair $(\rho,\tau)$ of Theorem~\ref{thm:family} satisfies the same hypotheses, with centre the cyclotomic field $\Q(\zeta_{p-1})$, and its Schur algebras have indices $p-1$ and $(p-1)/2$. These two indices are unchanged if $\tau$ is replaced by $\varphi_\pi$ or by $\varphi_\pi\otimes|\cdot|^{(n-1)/2}$, if $\rho$ is replaced by $\rho\otimes\nu_D^{(n-1)/2}$, or if either division algebra is replaced by its opposite.
\end{enumerate}
\end{corollary}

\begin{proof}
The representations $\rho_\theta$ and $\rho$ of $D^\times=\GL_1(D)$ are irreducible with finite-order central character, hence unitary discrete series. They are not self-dual, $\rho_\theta$ by Theorem~\ref{thm:order-eight}(i), and $\rho$ because its central character has order $n\ge6$ (Lemma~\ref{lem:division-side}(ii)) and the central character of a self-dual representation has order at most two. The normalization of \cite[\S6]{PR12} attaches $\tau_\theta$ and $\tau$ to them, and the common centre of the Schur algebras is $E$ \cite[Lemma~3.0.7]{Su23}. In (i), $o_E(\rho_\theta)=0$ and $o_E(\tau_\theta)\ne0$ by Theorem~\ref{thm:order-eight}, and the indices in (ii) are those of Theorem~\ref{thm:family}(iv). For the last assertion of (ii), $q$ is a square, so $|\cdot|^{\pm(n-1)/2}$ and $\nu_D^{(n-1)/2}$ take rational values. Hence $\varphi_\pi$ and $\varphi_\pi\otimes|\cdot|^{(n-1)/2}$ have the field of rationality and the image algebra $A_\tau$ of $\tau$, the twist $\rho\otimes\nu_D^{(n-1)/2}$ has those of $\rho$, and an opposite algebra has the same index.
\end{proof}

\begin{remark}\label{rem:cyclotomic-centre}
The conjecture describes the common centre of the two Schur algebras as a cyclotomic field \cite[Conjecture~6.1]{PR12}, and before the conjecture Prasad and Ramakrishnan describe this centre, for discrete series representations whose central character has finite order, as an abelian extension of $\Q$ \cite[\S6]{PR12}. We therefore read the conjecture in two ways. In the first reading, it concerns the pairs whose common centre is a cyclotomic field. In the second, it concerns all pairs in which $\pi'$ has central character of finite order, and the common centre of such a pair is an abelian extension of $\Q$.

The fields $\Q(i)$ and $\Q(\zeta_{p-1})$ are cyclotomic and the central characters have finite order, so the counterexamples of Corollary~\ref{cor:pr-counterexamples}(i) with $p\equiv5\pmod8$ and of Corollary~\ref{cor:pr-counterexamples}(ii) contradict the conjecture under both readings. The case $p=3$ contradicts the second reading.
\end{remark}

\begin{remark}\label{rem:hom-lemma}
At a place $v$ of $E$ above $p$, the representation $\pi_\theta^\vee\boxtimes\rho_\theta\boxtimes\tau_\theta$ of the order-eight examples has no model over $E_v$. Here we view the three representations over an algebraic closure $\overline{E_v}$ of $E_v$, which plays the role of $\overline F$ in \S\ref{subsec:brauer-class} for $F=E$. This is possible because each of them has a model over a number field. Indeed, $\pi_\theta$ has a model over $E$ by Suzuki's result on models \cite[Lemma~4.0.4]{Su23}, and hence so does $\pi_\theta^\vee$. The representation $\rho_\theta$ has a model over $E$ (Theorem~\ref{thm:order-eight}(ii)), and $\tau_\theta$ has one over $L$ (Corollary~\ref{cor:fields-of-definition}).

Suppose that the tensor product had a model over $E_v$. Two applications of Lemma~\ref{lem:hom-base-change}(ii), with $E_v$ and $\overline{E_v}$ in place of $F$ and $F'$, would then extract from it a model of $\tau_\theta$ over $E_v$ through Hom spaces, as Remark~\ref{rem:suzuki-step} extracts a model of $\rho\boxtimes\tau$ at the finite places not above $p$. The first takes for $U$ a model over $E_v$ of $\pi_\theta^\vee$. This model is infinite-dimensional. It is irreducible, since it becomes irreducible over $\overline{E_v}$, and hence cyclic, and in particular finitely generated, over the group algebra (this is why the lemma assumes finite generation). After extension of scalars to $\overline{E_v}$ the model becomes $\pi_\theta^\vee$, whose equivariant endomorphisms are scalars by admissibility. The first application produces a model of $\rho_\theta\boxtimes\tau_\theta$ over $E_v$.

The second application takes for $U$ a model over $E_v$ of the finite-dimensional, absolutely irreducible Jacquet--Langlands transfer $\rho_\theta$, and it produces a model of $\tau_\theta$ over $E_v$. This contradicts $\inv_v o_E(\tau_\theta)=1/2$ by Lemma~\ref{lem:models}(i). So an argument at $p$ that realizes these representations together must account for this failure of linear descent.
\end{remark}

\subsection{A proposition of Suzuki on the local invariants above \texorpdfstring{$p$}{p}}\label{subsec:dyadic}

The examples of this section also bear on a proposition of Suzuki, which bounds the local invariants of the two classes above $p$ \cite[Proposition~4.0.13]{Su23}. Let $\pi$ be a supercuspidal representation of $\GL_n(K)$ whose central character has finite order. Suzuki states that at each place $v$ of $\Q(\pi)$ above $p$, the local invariants of the Schur algebras of $\JL(\pi)$ and of the twisted parameter $\varphi_\pi(-\frac{n-1}2)$ lie in $\frac{1}{\gcd(n,p-1)}\Z/\Z$. This is the first assertion of the proposition. The second is that \cite[Conjecture~4.0.7]{Su23} holds if $\gcd(n,p-1)\le2$, and when $\JL(\pi)$ is not self-dual, this conjecture asserts that the two division algebras agree at every place above $p$. As in \S\ref{subsec:intro-history}, we read these assertions in our notation, with $\JL(\pi)=\rho$ and $\varphi_\pi(-\frac{n-1}2)=\tau$, only under the hypothesis of Remark~\ref{rem:suzuki-step}. The examples themselves are proved unconditionally.

In this reading, Theorem~\ref{thm:order-eight} contradicts the second assertion. There $\pi_\theta$ is supercuspidal, since its parameter $\varphi_\theta$ is irreducible, and its central character has finite order. Moreover $p$ is odd and $n=2$, so $\gcd(n,p-1)=2$. Since $\rho_\theta$ is not self-dual, the second assertion says that the Schur algebras of $\rho_\theta$ and $\tau_\theta$ agree at each place above $p$. By Theorem~\ref{thm:order-eight}, their invariants at such a place are $0$ and $\frac12$, which both lie in $\frac12\Z/\Z$ and are different.

In the same reading, the following example, with $p=n=2$, contradicts the first assertion. Both examples are cases of Theorem~\ref{thm:intro-depth-zero} with $n=2$ (for Example~\ref{ex:dyadic}, see the end of its proof), so they satisfy \eqref{eq:intro-sum}. In both, the invariants of the sum vanish at the finite places not above $p$, and since $n=2$ is even, their total over the places above $p$ is the one that Suzuki's formula gives. So in our reading both examples agree with \cite[Theorem~4.0.5 and~(4.0.6)]{Su23}.

\begin{example}\label{ex:dyadic}
For $x,y\in\Q^\times$, write $(x,y)_{\Q}$ for the quaternion algebra generated by $I,J$ with $I^2=x$, $J^2=y$ and $IJ=-JI$. Let $\zeta=\zeta_3$ and
\[
 T=\begin{pmatrix}\zeta&0\\0&\zeta^{-1}\end{pmatrix},\qquad U_0=\begin{pmatrix}0&-1\\1&0\end{pmatrix}.
\]
Let $\varphi$ be the representation of $W_{\Q_2}$ that is trivial on wild inertia and in which a topological generator $t$ of tame inertia acts by $T$ and the geometric Frobenius element $s$ acts by $U_0$. Then $\varphi$ is the Langlands parameter of a supercuspidal representation $\pi$ of $\GL_2(\Q_2)$ with trivial central character, $\Q(\pi)=\Q$, and the normalized parameter $\tau=\varphi\otimes|\cdot|^{-1/2}$ satisfies
\[
 A_\tau\cong(-3,-2)_{\Q},\qquad \inv_2o_{\Q}(\tau)=\tfrac12.
\]
In the notation of Section~\ref{sec:depth-zero}, the pair $(\JL(\pi),\tau)$ is $(\rho_1,\tau_1)$ for $K=\Q_2$, $n=2$ and a character $\theta$ of order $m=3$, so that $o_\Q(\JL(\pi))=0$ and $o_\Q(\tau)=[B_{2,\infty}]$. Since $\gcd(n,p-1)=\gcd(2,1)=1$, the first assertion of \cite[Proposition~4.0.13]{Su23}, read in our notation under the hypothesis of Remark~\ref{rem:suzuki-step}, would give $\inv_2o_{\Q}(\tau)=0$.
\end{example}

\begin{proof}[Proof of the assertions of Example~\ref{ex:dyadic}]
The matrices satisfy $T^3=1$, $U_0^2=-1$ and $U_0TU_0^{-1}=T^{-1}$. Since $2^{-1}\equiv-1\pmod 3$, the last of these is the tame relation $sts^{-1}=t^{2^{-1}}$, so $\varphi$ is a continuous representation, with image of order $12$. It is irreducible, since $U_0$ exchanges the eigenlines of $T$, and $\det T=\det U_0=1$, so $\varphi$ is the Langlands parameter of a supercuspidal representation $\pi$ of $\GL_2(\Q_2)$ with trivial central character.

Since $|s|^{-1/2}=\sqrt2$, the representation $\tau$ is given by $\tau(t)=T$ and $\tau(s)=U=\sqrt2\,U_0$, with $U^2=-2$. Conjugation by $U_0$ fixes $U$ and replaces $T$ by $T^{-1}$, and conjugation by $\diag(1,-1)$ fixes $T$ and replaces $U$ by $-U$. These conjugations realize the two generators of $\Gal(\Q(\zeta,\sqrt2)/\Q)$ on the matrix entries, so every automorphism of $\C$ preserves the isomorphism class of $\tau$, and $\Q(\tau)=\Q$. By Lemma~\ref{lem:common-rationality-field}, $\Q(\pi)=\Q$ as well.

The $\Q$-span of the image of $\tau$ is
\[
 A_\tau=\Q[T]\oplus\Q[T]U,\qquad \Q[T]=\{\diag(x,\bar x):x\in\Q(\zeta)\}\cong\Q(\zeta),
\]
with $UxU^{-1}=\bar x$ and $U^2=-2$, so $A_\tau\cong(\Q(\zeta)/\Q,x\mapsto\bar x,-2)$. With $I=1+2T$ we have $I^2=-3$ and $UI=-IU$, so $A_\tau\cong(-3,-2)_{\Q}$. The extension $\Q_2(\zeta)/\Q_2$ is unramified of degree $2$, because $X^2+X+1$ is irreducible over $\F_2$, so its norms have even valuation and $-2$ is not a norm. Hence $A_\tau\otimes\Q_2$ is not split, and $\inv_2o_{\Q}(\tau)=\frac12$.

Finally, let $\theta$ be the character of $\F_4^\times$ with $\theta(\bar t)=\zeta$, for the image $\bar t$ of $t$ under the surjection $I_{\Q_2}\to\F_4^\times$ of \S\ref{subsec:depth-zero-reps}, so that $m=3$ and $\theta$ is primitive. The Weil group of the unramified quadratic extension of $\Q_2$ is generated by inertia and by $s^2$, which acts by $U_0^2=-1$, so it acts on the first basis vector by the character $\mu_1$ of Definition~\ref{def:depth-zero-reps}, and $\varphi\cong\Ind\mu_1$ by Frobenius reciprocity. By Proposition~\ref{prop:depth-zero-parameter}, $\pi=\pi_1$ and $\JL(\pi)=\rho_1$, and Theorem~\ref{thm:central-value-one}, with $L=\Q(\zeta)$, $E=\Q$, $a=-2$ and $f=1$, gives $o_\Q(\rho_1)=0$ and $o_\Q(\tau_1)=[B_{2,\infty}]$.
\end{proof}

\section{Questions}\label{sec:questions}

We close with two questions and a conjecture that the computations of this paper leave open. They concern the invariant of the sum above $p$ for a general $\pi$, the individual classes, and the division algebras of other invariants.

\subsection{The general identity}\label{subsec:questions-general}

Conjecture~\ref{conj:intro-sum} holds for $n=1$, and this paper proves it for every depth-zero supercuspidal representation with algebraic-valued central character (Corollary~\ref{cor:depth-zero-relation}) and for the family of Theorem~\ref{thm:intro-family}. Both proofs use that $\rho$ and $\tau$ are induced from characters, which presents both classes by cyclic algebras (for $\rho$ in Theorem~\ref{thm:intro-family}, through the corner of degree $n$ of Lemma~\ref{lem:division-side}(iii)). At the finite places not above $p$, under the hypothesis of Remark~\ref{rem:suzuki-step}, the conjecture holds for every central character of finite order, the case $\omega_\pi(\varpi)=1$ included, and this vanishing is the form in our notation of the theorem that Suzuki states \cite[Theorem~4.0.5]{Su23}. Suzuki's proof, which uses the quotient of the tower by $\varpi^{\Z}$, treats the representations with $\omega_\pi(\varpi)=1$ (\S\ref{subsec:intro-history}). The main open part of the conjecture is the invariant at each place above $p$.

Suzuki's argument does not extend to the places above $p$. In Remark~\ref{rem:suzuki-step}, at a finite place $v$ not above $p$, we take the space of homomorphisms from a model over $E_v$ of $\pi^\vee$ into the cohomology $H$, a representation over $E_v$. Under the hypothesis of that remark, the $\pi^\vee$-isotypic component of $H$, after extension of scalars to an algebraic closure of $E_v$, is $\pi^\vee\boxtimes\rho\boxtimes\tau$, and this identifies the Hom space with a model of $\rho\boxtimes\tau$ over $E_v$. Above $p$, the representation $\pi^\vee\boxtimes\rho\boxtimes\tau$ need not have a model over $E_v$, as the order-eight examples show (Remark~\ref{rem:hom-lemma}), so a realization above $p$ that covers these examples needs more structure than a linear one over $E_v$.

The value that the conjecture predicts above $p$ suggests which structure, and to describe it we recall isocrystals. Let $\Qb$ be the completion of the maximal unramified extension of $\Q_p$, and let $\sigma$ be its Frobenius automorphism, which induces the absolute Frobenius $x\mapsto x^p$ on the residue field $\overline{\F}_p$ of $\Qb$.

\begin{definition*}
An \emph{isocrystal} is a finite-dimensional $\Qb$-vector space $N$ with a bijection $\Phi\colon N\to N$, its \emph{Frobenius}, that is semilinear for $\sigma$, which means that $\Phi$ is additive and $\Phi(cx)=\sigma(c)\Phi(x)$ for $c\in\Qb$ and $x\in N$. The isocrystal is \emph{simple} if $N\neq0$ and no nonzero proper subspace of $N$ is stable under $\Phi$. The \emph{endomorphisms} of the isocrystal are the $\Qb$-linear maps $N\to N$ that commute with $\Phi$.
\end{definition*}

The slopes of an isocrystal are rational numbers that measure the $p$-adic growth of $\Phi$, normalized by the number of iterations. For instance, if $N\neq0$ and, for some integer $k\ge1$, $\Phi^k$ multiplies every vector of some basis of $N$ by one scalar, then the isocrystal has a single slope, the $p$-adic valuation of this scalar divided by $k$.

Dieudonn\'e computed the endomorphism algebra of the simple isocrystal of slope $\lambda$. It is a central division algebra over $\Q_p$, with invariant $-\lambda$ in the normalization that Crew uses \cite[Theorem~1.1.2]{Cre25}. For $\lambda=\frac{n-1}2$ this invariant is $-\frac{n-1}2$, which is $\frac{n-1}2$ modulo $\Z$, and restriction to $E_v$ multiplies it by $[E_v:\Q_p]$ (Theorem~\ref{thm:local-brauer}(a)), which gives
\[
 \frac{n-1}2\,[E_v:\Q_p]\quad\text{in }\Q/\Z.
\]
This is the value that Conjecture~\ref{conj:intro-sum} predicts at a place $v$ of $E$ above $p$.

In depth zero the number $\frac{n-1}2$ comes from the half twist. The element $s^n$ acts on $\tau_z$ by $az$ (Lemma~\ref{lem:shape}(iii)), and the half twist contributes the factor $q^{n(n-1)/2}$ of $a=(-1)^{n-1}q^{n(n-1)/2}$. Since $q=p^f$, and since $s$ acts on the residue field $\overline{\F}_q$ of $\overline K$ as the inverse of $x\mapsto x^q$ (\S\ref{subsec:local}),
\[
 q^{n(n-1)/2}=p^{fn(n-1)/2},\qquad s^n\ \text{acts on }\overline{\F}_q\text{ as }(x\mapsto x^p)^{-fn}.
\]
So the factor has $p$-adic valuation $fn(n-1)/2$, and $s^n$ acts on the residue field as $\sigma^{-fn}$ acts on $\overline{\F}_p$. Dividing the valuation of the factor by this number $fn$ of iterations gives $\frac{n-1}2$. The same valuation of $a$ produces the invariant $\frac{n-1}2[E_v:\Q_p]$ in Theorem~\ref{thm:central-value-one}(i).

On a $p$-adic realization of $\rho_z\boxtimes\tau_z$ as in Hypothesis~\ref{hyp:isocrystal}, Fontaine's construction, which uses the Frobenius to linearize a commuting semilinear Galois action, relates the action of a power of the Frobenius to the scalar $az$. The half twist contributes the term $\frac{n-1}2$ found above to the resulting slope, which equals this term when $z$ is a unit in the $p$-adic field of coefficients. Let $\Q_p^{\mathrm{nr}}\subset\Qb$ be the maximal unramified extension of $\Q_p$, and suppose that a module over $\Q_p^{\mathrm{nr}}$ carries a bijective semilinear Frobenius $\Phi$ and a commuting action $w\mapsto g_w$ of $\Gal(\overline K/K)$ that is semilinear for the Galois action on $\Q_p^{\mathrm{nr}}$, such that an open subgroup fixes a basis. Then Fontaine's linearization
\[
 r(w)=g_w\Phi^{-f\nu(w)}\qquad(w\in W_K)
\]
is a linear action of the Weil group \cite[1.3.1 and~1.3.5]{Fon94}. The element $w$ acts on $\Q_p^{\mathrm{nr}}$ by $\sigma^{f\nu(w)}$, where $\nu$ is the degree of \S\ref{subsec:local}, so $g_w$ is semilinear for $\sigma^{f\nu(w)}$, the factor $\Phi^{-f\nu(w)}$ is semilinear for $\sigma^{-f\nu(w)}$, and their product is linear. In particular $r(s)=g_s\Phi^f$, since $\nu(s)=-1$.

Assume now Hypothesis~\ref{hyp:isocrystal}, under which $\rho_z\boxtimes\tau_z$ arises in this way, with coefficients in $E_{z,v}$ for a place $v$ of $E_z$ above $p$. We do not know whether the $p$-adic cohomology of the Lubin--Tate tower provides such a realization. Under this hypothesis, the operator $r(s)^n=g_{s^n}\Phi^{fn}$ acts on the whole realization by the scalar $az$ by which $s^n$ acts on $\tau_z$ (Proposition~\ref{prop:conditional-isocrystal}(i)). Hence $\Phi^{fn}=az\,g_{s^n}^{-1}$. Some power of $g_{s^n}$ fixes the basis that an open subgroup of the Galois group fixes, so a suitable iterate of $\Phi$ maps the lattice spanned by this basis onto a power of $p$ times itself, and this forces the single slope
\[
 \frac{v_p(az)}{fn}=\frac{n-1}2+\frac{v_p(z)}{fn}
\]
(Proposition~\ref{prop:conditional-isocrystal}(ii)).

If $z$ is a unit at $v$, the slope is $\frac{n-1}2$. Proposition~\ref{prop:conditional-isocrystal}(iii) then shows that the Frobenius commutant with coefficients in $E_{z,v}$, the algebra of endomorphisms of the isocrystal that are linear over $E_{z,v}$, is the image algebra of $\rho_z\boxtimes\tau_z$ over $E_{z,v}$, with invariant $-\frac{n-1}2[E_{z,v}:\Q_p]$, the value predicted by Conjecture~\ref{conj:intro-sum} modulo $\Z$. Appendix~\ref{app:isocrystal} gives the arguments on the fibre and on lattices behind these assertions.

For an arbitrary $\pi$, these observations suggest computing the invariant of $o_E(\rho)+o_E(\tau)$ at a place $v$ of $E$ above $p$ from a realization of $\rho\boxtimes\tau$ on an isocrystal, as the invariant of an algebra of endomorphisms commuting with its Frobenius.

\begin{question}\label{q:above-p}
Let $\pi$ have central character of finite order, let $c\in K^\times$ have positive valuation and satisfy $\rho(c)=1$, and let $v$ be a place of $E$ above $p$. Can the invariant of $o_E(\rho)+o_E(\tau)$ at $v$ be computed from the $p$-adic cohomology of the quotient of the Lubin--Tate tower by the subgroup of $D^\times$ generated by $c$, with its Frobenius and its actions of $\GL_n(K)$, $D^\times$ and $W_K$?
\end{question}

Such $c$ exists for every central character of finite order, since $\varpi$ acts on $\rho$ by a root of unity of some order $d$, and $c=\varpi^d$ gives the quotient that Remark~\ref{rem:suzuki-step} uses at the finite places not above $p$. For algebraic central characters of infinite order, however, descent over the field of rationality is an additional issue, since the reduction of Lemma~\ref{lem:common-rationality-field} passes to an unramified twist of $\pi$ with finite-order central character. For $n\ge2$ this paper proves \eqref{eq:intro-sum} for such central characters only in depth zero, where the cyclic presentations show directly that $z$ cancels from the sum.

\subsection{The class of the transfer above \texorpdfstring{$p$}{p}}\label{subsec:questions-individual}

If Conjecture~\ref{conj:intro-sum} holds for $\pi$, it determines the class of $\tau$ from that of $\rho$. The question is then what determines $o_E(\rho)$ at the places above $p$. In depth zero with central character of finite order, $\inv_vo_E(\rho)=0$ at every place $v$ above $p$ (Theorem~\ref{thm:intro-depth-zero}(ii) and Proposition~\ref{prop:depth-zero-exhaustion}). In Theorem~\ref{thm:intro-family}, however, the invariant $\inv_{v_j}o_E(\rho)=j/n$ varies with the place and takes every value of order $n$ in $\Q/\Z$. Question~\ref{q:above-p} concerns the sum, the class of $\rho\boxtimes\tau$ (Lemma~\ref{lem:sum-class}), which determines neither summand, and in depth zero and for the family the cyclic presentations of Sections~\ref{sec:depth-zero} and~\ref{sec:individual} compute each class separately.

\begin{question}\label{q:individual}
Let $v$ be a place of $E$ above $p$. How does $\inv_vo_E(\rho)$ depend on $\pi$ and on $v$?
\end{question}

\subsection{Division algebras of other invariants}\label{subsec:questions-inner}

Let $t$ be an integer prime to $n$, let $D^{(t)}$ be the central division algebra over $K$ of invariant $t/n$, so that $D^{(1)}=D$, and let $\rho^{(t)}=\JL_{D^{(t)}}(\pi)$ be the Jacquet--Langlands transfer of $\pi$ to $(D^{(t)})^\times$, so that $\rho^{(1)}=\rho$. When the central character of $\pi$ has finite order, $\rho^{(t)}$ again has field of rationality $E$ \cite[Lemma~4.0.3]{Su23}.

If Conjecture~\ref{conj:intro-sum} holds for $\pi$, then subtracting $t$ times \eqref{eq:intro-sum} from the identity of the following conjecture shows that this identity is equivalent to
\[
 o_E(\rho^{(t)})=t\,o_E(\rho).
\]
Thus, together with Conjecture~\ref{conj:intro-sum}, the following conjecture predicts that the Brauer class of the transfer to $D^{(t)}$ is $t$ times the class of the transfer to $D$.

\begin{conjecture}\label{conj:questions-inner-forms}
Let $\pi$ be a supercuspidal representation of $\GL_n(K)$ whose central character has finite order, and let $t$ be an integer prime to $n$. Then
\[
 o_E(\rho^{(t)})+t\,o_E(\tau)=t(n-1)\res_{E/\Q}[\Bpinf]\qquad\text{in }\Br(E).
\]
\end{conjecture}

For $t=1$ this is Conjecture~\ref{conj:intro-sum} for central characters of finite order.

\appendix

\section{Proofs of \texorpdfstring{Lemma~\ref{lem:hom-base-change} and Proposition~\ref{prop:intertwiner-sign}}{Lemma 2.5 and Proposition 2.12}}
\label{app:section-two-proofs}

This appendix proves Lemma~\ref{lem:hom-base-change} and Proposition~\ref{prop:intertwiner-sign}, in the notation of \S\S\ref{subsec:hom-spaces} and~\ref{subsec:factor-set}.

\begin{proof}[Proof of Lemma~\ref{lem:hom-base-change}]
(i) Let $x_1,\dots,x_r$ generate $X$ over $\Lambda$, and let $\xi'\colon X\otimes_FF'\to Y\otimes_FF'$ be $\Lambda\otimes_FF'$-linear. Choose one finite family $c_1,\dots,c_k\in F'$, linearly independent over $F$, and elements $y_{ij}\in Y$ with $\xi'(x_i\otimes1)=\sum_jy_{ij}\otimes c_j$ for $1\le i\le r$. If $\sum_ia_ix_i=0$ with $a_i\in\Lambda$, then $\sum_j\bigl(\sum_ia_iy_{ij}\bigr)\otimes c_j=0$, so $\sum_ia_iy_{ij}=0$ for every $j$. Thus every relation among the generators holds coefficient by coefficient. So, for each $j$, the assignment $x_i\mapsto y_{ij}$ extends to a $\Lambda$-linear map $\xi_j\colon X\to Y$, and $\xi'=\sum_j\xi_j\otimes c_j$, which proves surjectivity. For injectivity, write an element of the kernel as $\sum_l\eta_l\otimes d_l$ with $\eta_l\in\Hom_\Lambda(X,Y)$ and with $d_l\in F'$ linearly independent over $F$. For every $x\in X$, then, $\sum_l\eta_l(x)\otimes d_l=0$, so $\eta_l(x)=0$ for every $l$. Hence every $\eta_l$ vanishes, and so does the element.

(ii) Part (i), with $\Lambda=F[\Gamma_1]$, identifies $\Hom_{\Gamma_1}(U,Y)\otimes_FF'$ with $\Hom_{\Gamma_1}(U_{F'},Y_{F'})$, and this isomorphism commutes with $\Gamma_2$, which acts on both sides through $Y$. Composing with the given isomorphism $Y_{F'}\cong U_{F'}\boxtimes T'$, which also commutes with $\Gamma_2$, identifies $\Hom_{\Gamma_1}(U,Y)\otimes_FF'$ with $\Hom_{\Gamma_1}(U_{F'},U_{F'}\otimes_{F'}T')$ as representations of $\Gamma_2$. Finally, $v\mapsto(u\mapsto u\otimes v)$ identifies $T'$ with the latter, compatibly with $\Gamma_2$, since in a basis of $T'$ a $\Gamma_1$-linear map $U_{F'}\to U_{F'}\otimes_{F'}T'$ is given by elements of $\End_{\Gamma_1}(U_{F'})=F'$.
\end{proof}

We now prove the assertions made after Definition~\ref{def:factor-set}. Since $\iota_g\,g(\iota_h)$ and $\iota_{gh}$ are both isomorphisms $r^{gh}\to r$, each $c(g,h)$ is an automorphism of $r$, hence a scalar in $L^\times$, the commutant of $r$ in $M_d(L)$ being $L$. Expanding $\iota_g\,g(\iota_h)\,(gh)(\iota_k)$ in two ways gives the identity of Definition~\ref{def:factor-set}(i). Any other choice of isomorphisms has the form $\lambda_g\iota_g$ with $\lambda_g\in L^\times$, and its factor set is $c'(g,h)=\lambda_g\,g(\lambda_h)\,\lambda_{gh}^{-1}\,c(g,h)$. The map $(L/F,c')\to(L/F,c)$ that is the identity on $L$ and sends the generator $u'_g$ of $(L/F,c')$ to $\lambda_gu_g$ respects the relations of Definition~\ref{def:factor-set}(i) and carries the $L$-basis $(u'_g)$ to the $L$-basis $(\lambda_gu_g)$, so it is an isomorphism, and the class of $(L/F,c)$ depends only on $r$.

\begin{proof}[Proof of Proposition~\ref{prop:intertwiner-sign}]
By Lemma~\ref{lem:trace-span}, the $F$-span $A_r$ of $r(\Gamma)$ is a central simple $F$-algebra of degree $d$. The induced map $A_r\otimes_FL\to M_d(L)$ is injective, since its source is simple, and hence bijective, since both sides have dimension $d^2$ over $L$. So the $L$-span of $r(\Gamma)$ is $M_d(L)$, and the commutant of $r$ is the centre $L$ of $M_d(L)$.

The isomorphism $\iota_1$ is an automorphism of $r$, hence a scalar, and the class of $(L/F,c)$ depends only on $r$. We may therefore assume $\iota_1=1$. Then $c(1,h)=c(g,1)=1$ for all $g,h\in G$ by the definition of $c$, the element $u_1$ is the identity of $(L/F,c)$, and $(L/F,c)$ is a central simple $F$-algebra \cite[IV, Lemma~3.13]{Mil20}.

Let $V=L^d$, viewed as an $F$-vector space of dimension $d[L:F]$ on which $A_r$ acts by matrix multiplication, and let $C=\End_{A_r}(V)$ be the centralizer of $A_r$ in $\End_F(V)$. By the double centralizer theorem \cite[IV, Theorem~3.1 and Corollary~3.3]{Mil20}, the algebra $C$ is central simple over $F$ of dimension $[L:F]^2$, and multiplication $A_r\otimes_FC\to\End_F(V)$ is an isomorphism. Hence $[C]=-[A_r]=-o_F(r)$.

For $x\in L$ let $m_x$ be multiplication by $x$ on $V$, and for $g\in G$ define $U_g\colon V\to V$ by $U_g(v)=\iota_g\,g(v)$, where $g$ acts on coordinates. Both maps are $F$-linear. Every $m_x$ commutes with $A_r$, since $A_r\subseteq M_d(L)$. For $\delta\in\Gamma$ we have $\iota_g\,g(r(\delta))=r(\delta)\,\iota_g$, because $\iota_g$ is an isomorphism $r^g\to r$, so
\[
 U_g\bigl(r(\delta)v\bigr)=\iota_g\,g\bigl(r(\delta)\bigr)\,g(v)=r(\delta)\,U_g(v),
\]
and $U_g$ commutes with $A_r$.

The operators $m_x$ and $U_g$ satisfy the relations of the crossed product. Indeed, $U_gm_x=m_{g(x)}U_g$, since $g(xv)=g(x)\,g(v)$ and $\iota_g$ is $L$-linear, and $U_gU_h=m_{c(g,h)}U_{gh}$ because $\iota_g\,g(\iota_h)=c(g,h)\,\iota_{gh}$. So $xu_g\mapsto m_xU_g$ defines a homomorphism of $F$-algebras $(L/F,c)\to C$. It is injective because $(L/F,c)$ is simple, and bijective because both algebras have dimension $[L:F]^2$ over $F$. Therefore $[(L/F,c)]=[C]=-o_F(r)$.
\end{proof}

\section{A conditional computation above \texorpdfstring{$p$}{p}}
\label{app:isocrystal}

This appendix proves Proposition~\ref{prop:conditional-isocrystal}, which makes the heuristic of \S\ref{subsec:questions-general} precise. We use the representations $\rho_z$ and $\tau_z$ of Definition~\ref{def:depth-zero-reps} for some $z\in\Qbar^\times$, with $q=p^f$ and $a=(-1)^{n-1}q^{n(n-1)/2}$, so that $s^n$ acts on $\tau_z$ by $az$ (Lemma~\ref{lem:shape}(iii)). The representations $\pi_z$, $\rho_z$ and $\tau_z$ have common field of rationality $E_z$ (Corollary~\ref{cor:depth-zero-relation}). Let $v$ be a place of $E_z$ above $p$, put $k=E_{z,v}$, and normalize the valuation $v_p$ of $k$ by $v_p(p)=1$.

Fix an algebraic closure $\overline K$ of $K$, let $\Q_p^{\mathrm{nr}}$ be the maximal unramified extension of $\Q_p$ in $\overline K$, and let $\Qb$ and $\sigma$ be as in \S\ref{subsec:questions-general}. The realization takes coefficients in the rings
\[
 R_{k,0}=\Q_p^{\mathrm{nr}}\otimes_{\Q_p}k,\qquad R_k=\Qb\otimes_{\Q_p}k,
\]
with $k$ acting through the second factor. The group $G_K=\Gal(\overline K/K)$ acts on $\Q_p^{\mathrm{nr}}$ and continuously on $\Qb$, and an element $w\in W_K$ acts on both as $\sigma^{f\nu(w)}$, where $\nu$ is the degree of \S\ref{subsec:local}.

\begin{hypothesis}\label{hyp:isocrystal}
There is a nonzero free $R_{k,0}$-module $M_0$ of finite rank with
\begin{enumerate}
\item[(a)] a bijective map $\Phi\colon M_0\to M_0$ that is semilinear for $\sigma\otimes1$,
\item[(b)] an action $w\mapsto g_w$ of $G_K$ on $M_0$ by maps that are semilinear for the action of $G_K$ on the first factor of $R_{k,0}$, such that an open subgroup of $G_K$ fixes an $R_{k,0}$-basis of $M_0$, and
\item[(c)] an $R_{k,0}$-linear action of $D^\times$,
\end{enumerate}
such that these three actions commute and the following holds. Extend $\Phi$, $g$ and the action of $D^\times$ to $M=\Qb\otimes_{\Q_p^{\mathrm{nr}}}M_0$, and put $r(w)=g_w\Phi^{-f\nu(w)}$ for $w\in W_K$. Then, for some algebraically closed field $\Omega$ containing $\Qb$ and some $\Q_p$-embedding $k\to\Omega$, the fibre $M\otimes_{R_k}\Omega$ along the homomorphism $R_k\to\Omega$ given by these two embeddings, with $D^\times$ acting through $M$ and $w\in W_K$ acting by $r(w)$, is a representation of $D^\times\times W_K$ isomorphic to $\rho_z\boxtimes\tau_z$.
\end{hypothesis}

The formula for $r$ is Fontaine's linearization of a semilinear Galois action \cite[1.3.1 and~1.3.5]{Fon94}, recalled in \S\ref{subsec:questions-general}. By (a), (b) and the commutation of $\Phi$ with the maps $g_w$, each $r(w)$ is $R_k$-linear and $r(ww')=r(w)r(w')$ (Proposition~\ref{prop:conditional-isocrystal}(i)). The representations $\rho_z$ and $\tau_z$ have models over $L_z$, obtained by inducing characters with values in $L_z$ (Lemma~\ref{lem:shape}(i)), and in the hypothesis we view $\rho_z\boxtimes\tau_z$ over $\Omega$ through an embedding of $L_z$ that extends $E_z\subseteq k\to\Omega$. Two such embeddings differ by an element of $\Gal(L_z/E_z)$, so the resulting representations are isomorphic, $E_z$ being the field of rationality.

\begin{proposition}\label{prop:conditional-isocrystal}
Assume Hypothesis~\ref{hyp:isocrystal}.
\begin{enumerate}
\item[(i)] The maps $r(w)$ form an $R_k$-linear action of $W_K$ on $M$ that commutes with $\Phi$ and with $D^\times$, $r(s)=g_s\Phi^f$, and $r(s)^n=az$ on $M$.
\item[(ii)] The isocrystal $(M,\Phi)$ over $\Qb$ has the single slope
\[
 \lambda=\frac{v_p(az)}{fn}=\frac{n-1}2+\frac{v_p(z)}{fn}.
\]
\item[(iii)] Suppose that $z$ is a unit at $v$. Then the algebra $C=\End_{R_k,\Phi}(M)$ of $R_k$-linear endomorphisms of $M$ that commute with $\Phi$ is a central simple $k$-algebra of degree $n^2$. It is the $k$-span of the products $x\,r(w)$ for $x\in D^\times$ and $w\in W_K$, that is, of the operators of $D^\times\times W_K$ on $M$, and
\[
 [C]=\res_{k/E_z}\bigl(o_{E_z}(\rho_z)+o_{E_z}(\tau_z)\bigr),\qquad \inv_k[C]=-[k:\Q_p]\,\frac{n-1}2.
\]
In particular this invariant is $\frac{n-1}2[k:\Q_p]$ in $\Q/\Z$, the value that Conjecture~\ref{conj:intro-sum} predicts at $v$.
\end{enumerate}
\end{proposition}

We first spread the scalar identity $r(s)^n=az$ from the fibre, where it holds by Lemma~\ref{lem:shape}(iii), to all of $M$. The fibre is formed over one of the field factors of $R_k$, but $\sigma\otimes1$ permutes these factors transitively and $r(s)^n-az$ commutes with $\Phi$, so the identity holds on every factor once it holds on one. In a basis fixed by an open subgroup of $G_K$, the identity becomes a lattice identity for a power of $\Phi$, which determines the slope. Finally, Dieudonn\'e's theorem and the centralizer theorem compute the Frobenius commutant, and Burnside's theorem, applied to the irreducible fibre, identifies it with the image algebra.

\begin{proof}
(i) The map $g_w$ is semilinear for $\sigma^{f\nu(w)}$ on the first factor of $R_k$, and $\Phi^{-f\nu(w)}$ is semilinear for $\sigma^{-f\nu(w)}$, so $r(w)$ is $R_k$-linear. Since $g$ and $\Phi$ commute and $\nu$ is additive, $r(ww')=r(w)r(w')$, and $r$ commutes with $\Phi$ and with $D^\times$. Since $\nu(s)=-1$, we have $r(s)=g_s\Phi^f$.

It remains to show that $r(s)^n=az$ on $M$. Let $B=r(s)^n-az$. It is $R_k$-linear, and it commutes with $\Phi$, because $r(s)$ does and $az\in k$ is fixed by $\sigma\otimes1$. Let $k_0$ be the maximal unramified subextension of $k/\Q_p$. The ring $\Qb\otimes_{\Q_p}k_0$ is a product of copies of $\Qb$, one for each embedding of $k_0$ into $\Qb$, and since $k/k_0$ is totally ramified, each copy stays a field after base change to $k$. So $R_k$ is a finite product of fields, one for each such embedding \cite[Lemma~1.2.1]{Cre25}. The automorphism $\sigma\otimes1$ permutes these factors transitively, as $\sigma$ permutes the embeddings of $k_0$ transitively. Accordingly $M$ is the direct sum of the submodules $eM$ for the primitive idempotents $e$ of $R_k$, and $\Phi$ maps $eM$ onto $(\sigma\otimes1)(e)M$.

We first show that $B$ vanishes on one of these summands. The map $R_k\to\Omega$ given by $\Qb\subseteq\Omega$ and the embedding of $k$ factors through one of the field factors, say $e_0R_k$, and $M\otimes_{R_k}\Omega=e_0M\otimes_{e_0R_k}\Omega$. On this fibre $r(s)^n$ acts as $\tau_z(s)^n=az$ (Lemma~\ref{lem:shape}(iii)), and $e_0R_k\to\Omega$ is a field extension, so $B$ vanishes on $e_0M$. Since $B$ commutes with $\Phi$, it vanishes on every $\Phi^i(e_0M)=(\sigma\otimes1)^i(e_0)M$, and these are all the summands $eM$, because $\sigma\otimes1$ permutes the primitive idempotents transitively. Hence $B$ vanishes on $M$.

(ii) By (i), $g_{s^n}\Phi^{fn}=r(s^n)=az$, so $\Phi^{fn}=az\,g_{s^n}^{-1}$. Let $e_1,\dots,e_N$ be an $R_{k,0}$-basis of $M_0$ fixed by an open subgroup $U$ of $G_K$. The image of $s$ in the finite group $G_K/U'$, where $U'$ is the largest normal subgroup of $G_K$ contained in $U$, has finite order, so $s^{nt}\in U$ for some $t\geq1$, and we may choose $t$ so that $c=t\,v_p(az)$ is an integer. The element $s^{nt}$ acts on $\Qb$ as $\sigma^{-fnt}$ and fixes the $e_j$, so
\[
 \Phi^{fnt}\Bigl(\sum_jx_je_j\Bigr)=(az)^t\sum_j(\sigma^{fnt}\otimes1)(x_j)\,e_j \qquad(x_j\in R_k).
\]
Write $(az)^t=p^c\eta$ with $\eta\in\OO_k^\times$. The lattice $\Lambda=\bigoplus_j\bigl(W(\overline{\F}_p)\otimes_{\Z_p}\OO_k\bigr)e_j$ of the $\Qb$-vector space $M$, where $W(\overline{\F}_p)$ is the ring of integers of $\Qb$, is preserved by $\sigma^{fnt}\otimes1$ and by multiplication by $\eta$, so
\[
 \Phi^{fnt}(\Lambda)=p^c\Lambda.
\]

Iterating this identity forces every slope of $M$ to equal $c/(fnt)$. To see this, we compare $\Lambda$ with the lattices of the Dieudonn\'e--Manin classification \cite[Theorem~1.1.1]{Cre25}, by which $(M,\Phi)$ is the direct sum of its slope components, and a component $S$ of slope $a'/b$, in lowest terms, has a lattice $\Lambda_S$ with $\Phi^b(\Lambda_S)=p^{a'}\Lambda_S$. Since $\Phi(S)=S$, the lattice $\Lambda'_S=\Lambda\cap S$ of $S$ satisfies $\Phi^{fnt}(\Lambda'_S)=p^c\Lambda'_S$. Thus suitable iterates of $\Phi$ multiply $\Lambda_S$ and $\Lambda'_S$ by powers of $p$, at the rates $a'/b$ and $c/(fnt)$ per iteration. Since the two lattices lie within a fixed power of $p$ of each other, and a common iterate of $\Phi$ preserves this relation, these rates must agree.

To make this precise, choose $N'$ with $p^{N'}\Lambda_S\subseteq\Lambda'_S\subseteq p^{-N'}\Lambda_S$. For an integer $h\geq1$, the map $\Phi^{bfnth}$ multiplies $\Lambda_S$ by $p^{a'fnth}$ and $\Lambda'_S$ by $p^{cbh}$, so applying it to these inclusions and multiplying by $p^{-cbh}$ gives
\[
 p^{N'+(a'fnt-cb)h}\Lambda_S\subseteq\Lambda'_S \subseteq p^{-N'+(a'fnt-cb)h}\Lambda_S.
\]

Since $S\neq0$, comparing these inclusions with the ones we started from shows that $|a'fnt-cb|\,h\leq2N'$ for every $h$. Hence $a'fnt=cb$, that is, $a'/b=c/(fnt)=v_p(az)/(fn)$ for every slope $a'/b$ of $M$. Finally $v_p(a)=v_p(q^{n(n-1)/2})=fn(n-1)/2$, which gives the formula for $\lambda$.

(iii) Here $v_p(z)=0$, so $\lambda=\frac{n-1}2$ by (ii). Let $C_0=\End_{\Qb,\Phi}(M)$ be the endomorphism algebra of the isocrystal $(M,\Phi)$. By (ii) and \cite[Theorem~1.1.1]{Cre25}, the isocrystal $(M,\Phi)$ is a direct sum of copies of the simple isocrystal of slope $\lambda$, and by Dieudonn\'e's theorem \cite[Theorem~1.1.2]{Cre25} the endomorphism algebra of the latter is a central division algebra over $\Q_p$ with invariant $-\lambda$. So $C_0$ is a central simple $\Q_p$-algebra with invariant $-\lambda$, of dimension $(\dim_{\Qb}M)^2=([k:\Q_p]N)^2$, where $N$ is the rank of $M$ over $R_k$.

Next we impose $k$-linearity. The coefficient action of $k$ on $M$ is faithful and commutes with $\Phi$, so it makes $k$ a subfield of $C_0$. An element of $C_0$ is $R_k$-linear exactly when it commutes with this subfield, so the centralizer of $k$ in $C_0$ is $C$. By the centralizer theorem \cite[IV, Theorem~3.1]{Mil20}, $C$ is simple, $k$ is the centralizer of $C$, so that the centre of $C$ is $k$, and $[k:\Q_p]\dim_{\Q_p}C=\dim_{\Q_p}C_0$. Hence $C$ is central simple over $k$ with $\dim_kC=N^2$.

To compute the class of $C$, let $C_0\otimes_{\Q_p}k$ act on $C_0$ by left multiplication by $C_0$ and right multiplication by $k$. An endomorphism of this module commutes with left multiplication, so it is a right multiplication by an element of $C_0$, and it commutes with right multiplication by $k$ exactly when this element lies in $C$. Right multiplications compose in the opposite order, so the endomorphism algebra of this module is $C^{\op}$. As in the proof of Lemma~\ref{lem:models}(ii), the endomorphism algebra of a nonzero finitely generated module over a central simple algebra represents the opposite class, so $[C]=[C_0\otimes_{\Q_p}k]=\res_{k/\Q_p}[C_0]$, and $\inv_k[C]=-[k:\Q_p]\lambda$ by Theorem~\ref{thm:local-brauer}(a).

Let $A$ be the $k$-span of the products $x\,r(w)$ for $x\in D^\times$ and $w\in W_K$, that is, of the operators of $D^\times\times W_K$ on $M$. These operators are $R_k$-linear and commute with $\Phi$, so $A$ lies in $C$. The natural map $C\otimes_k\Omega\to\End_\Omega(M\otimes_{R_k}\Omega)$ is injective, since its source is simple, and hence bijective, since both sides have dimension $N^2$. The fibre is the irreducible representation $\rho_z\boxtimes\tau_z$ of $D^\times\times W_K$ (see the proof of Lemma~\ref{lem:sum-class}), so by Burnside's theorem the operators of $D^\times\times W_K$ on the fibre span $\End_\Omega(M\otimes_{R_k}\Omega)$, and the image of $A\otimes_k\Omega$ is all of $\End_\Omega(M\otimes_{R_k}\Omega)$. Hence $\dim_kA\geq N^2=\dim_kC$, and $A=C$.

It remains to identify $C$ with the image algebra of $\rho_z\boxtimes\tau_z$ extended to $k$. Through the injective map $C\to\End_\Omega(M\otimes_{R_k}\Omega)$, the algebra $C$ is the $k$-span of the operators of $\rho_z\boxtimes\tau_z$ on the fibre. This $k$-span is also the image of the natural map from $A_{\rho_z\boxtimes\tau_z}\otimes_{E_z}k$, which is injective. So $[C]=\res_{k/E_z}(o_{E_z}(\rho_z)+o_{E_z}(\tau_z))$ by Lemma~\ref{lem:sum-class}. Since $n=\dim\rho_z=\dim\tau_z$, the degree of $C$ is $N=n^2$. Finally, $-[k:\Q_p]\frac{n-1}2$ and $[k:\Q_p]\frac{n-1}2$ differ by the integer $[k:\Q_p](n-1)$.
\end{proof}

By Theorem~\ref{thm:arbitrary-central-value}, the class $\res_{k/E_z}(o_{E_z}(\rho_z)+o_{E_z}(\tau_z))$ does have the invariant $\frac{n-1}2[k:\Q_p]$, so for $z$ a unit at $v$, Proposition~\ref{prop:conditional-isocrystal} recovers, under Hypothesis~\ref{hyp:isocrystal}, what Section~\ref{sec:depth-zero} proves directly. The proposition is thus a model for one possible computation in Question~\ref{q:above-p}. In the proposition, the invariant is that of the Frobenius commutant over the coefficient field.

\section*{Acknowledgments}

We would like to express our sincere gratitude to Michael Harris, who supervised the senior thesis at Columbia University from which this paper grew, for introducing this problem, which led us to so many beautiful subjects, and we are very grateful for his guidance and for the things we learned from him. We thank Kenta Suzuki for correspondence about \cite[Theorem~4.0.5]{Su23}.

\end{document}